\documentclass[11pt,reqno]{amsart}
\makeatletter
\@namedef{subjclassname@2020}{\textup{2020} Mathematics Subject Classification}
\makeatother
\usepackage{amssymb,amscd,amsthm}
\usepackage[T1]{fontenc}
\usepackage{blkarray}

\usepackage{amsmath,amssymb,graphicx,mathrsfs} 
\usepackage[colorlinks=true,urlcolor=blue,citecolor=red,linkcolor=blue]{hyperref}

\let\frak\mathfrak
\let\Bbb\mathbb

\def\>{\relax\ifmmode\mskip.666667\thinmuskip\relax\else\kern.111111em\fi}
\def\<{\relax\ifmmode\mskip-.333333\thinmuskip\relax\else\kern-.0555556em\fi}
\def\vsk#1>{\vskip#1\baselineskip}
\def\vv#1>{\vadjust{\vsk#1>}\ignorespaces}
\def\vvn#1>{\vadjust{\nobreak\vsk#1>\nobreak}\ignorespaces}

  \let\ssize\scriptstyle
\let\sssize\scriptscriptstyle

\let\Medskip\medskip
\def\medskip{\par\Medskip}
\let\Bigskip\bigskip
\def\bigskip{\par\Bigskip}

\let\Maketitle\maketitle
\def\maketitle{\Maketitle\thispagestyle{empty}\let\maketitle\empty}

\newtheorem{thm}{Theorem}[section]
\newtheorem{cor}[thm]{Corollary}
\newtheorem{lem}[thm]{Lemma}
\newtheorem{prop}[thm]{Proposition}

\numberwithin{equation}{section}

\theoremstyle{definition}
\newtheorem{defn}[thm]{Definition}
\newtheorem{rem}[thm]{Remark}

\newtheorem{example}[thm]{Example}

\let\mc\mathcal
\let\nc\newcommand

\let\al\alpha
\let\bt\beta
\let\dl\delta
\let\Dl\Delta
\let\eps\varepsilon
\let\gm\gamma
\let\Gm\Gamma

\let\la\lambda
\let\La\Lambda

\let\phi\varphi
\let\si\sigma
\let\Si\Sigma

\let\thi\vartheta

\let\om\omega
\let\Om\Omega

\let\der\partial
\let\Hat\widehat

\let\Tilde\widetilde

\let\geq\geqslant

\let\leq\leqslant

\let\on\operatorname
\let\bi\bibitem
\let\bs\boldsymbol

\def\C{{\mathbb C}}
\def\Z{{\mathbb Z}}
\def\Q{{\mathbb Q}}

\def\Pb{{\mathbb P}}

\def\F{{\mc F}}

\def\+#1{^{\{#1\}}}

\def\diag{\on{diag}}
\def\End{\on{End}}

\def\Hom{\on{Hom}}

\def\gl{\mathfrak{gl}}

\def\sln{\mathfrak{sl}_N}

\def\beq{\begin{equation}}
\def\eeq{\end{equation}}
\def\be{\begin{equation*}}
\def\ee{\end{equation*}}

\nc{\bea}{\begin{eqnarray*}}
\nc{\eea}{\end{eqnarray*}}
\nc{\bean}{\begin{eqnarray}}
\nc{\eean}{\end{eqnarray}}
\nc{\bal}{\begin{align*}}
\nc{\eal}{\end{align*}}
\nc{\baln}{\begin{align}}
\nc{\ealn}{\end{align}}

\nc{\Il}{{\mc I_{\bs\la}}}
\nc{\bla}{{\bs\la}}
\nc{\Fla}{\F_\bla}
\nc{\tfl}{{{T}^*\Fla}}
\nc{\GL}{{GL_n(\C)}}
\nc{\GLC}{{GL_n(\C)\times\C^*}}

\let\sd s 

\def\ddk_#1{\kk_{#1}\<\>\frac\der{\der\<\>\kk_{#1}}}

\def\bul{\mathbin{\raise.2ex\hbox{$\sssize\bullet$}}}
\def\intt{\mathchoice
{\mathop{\raise.2ex\rlap{$\,\,\ssize\backslash$}{\intop}}\nolimits}
{\mathop{\raise.3ex\rlap{$\,\sssize\backslash$}{\intop}}\nolimits}
{\mathop{\raise.1ex\rlap{$\sssize\>\backslash$}{\intop}}\nolimits}
{\mathop{\rlap{$\sssize\<\>\backslash$}{\intop}}\nolimits}}

\let\kk q 
\let\cc c

\let\Ko K

\def\GZ/{Gelfand-Zetlin}
\def\KZ/{{\slshape KZ\/}}
\def\qKZ/{{\slshape qKZ\/}}
\def\XXX/{{\slshape XXX\/}}

\nc{\slnl}{{\sln (\lambda)}}
\nc{\PCN}{{   (\C[x])^N   }}
\nc{\di}{\on{Diag}}
\nc{\dio}{\on{Diag}_0}
\nc{\Mm}{{\mc M}}
\nc{\Nn}{{\mc N}}

\nc{\A}{{\mc C}}

\nc{\PCr}{{  P  (\C[x])^n   }}

\nc{\Pk}{{(\bs{P}^1)^k}}

\nc{\N}{{\Bbb N}}

\nc{\Ll}{{\mc L}}

\nc{\ord}{{\on{ord}\,}}

\nc{\Sing}{{\on{Sing}\,}}
\nc{\sing}{{\on{Sing}\,}}

\nc{\Hess}{{\on{Hess}}}

\nc{\R}{{\Bbb R}}

\let\on\operatorname
\nc{\Kk}{{\bs K}}
\nc{\Ap}{{A_\Phi(z)}}
\nc{\ap}{{A_\Phi(z)}}

\nc{\sv}{{\sing V}}
\nc{\cd}{{\C^n-\Delta}}
\nc{\UT}{{U^0}}   
\nc{\ep}{\epsilon}

\usepackage[OT2,T1]{fontenc}
\usepackage[all,cmtip]{xy}
\usepackage{cancel}
\usepackage[vcentermath]{youngtab}
\usepackage{empheq}
\usepackage{changepage}
\usepackage{bm}

\newlanguage\fakelanguage
\newcommand\cyr{\fontencoding{OT2}\fontfamily{wncyr}\selectfont
   \language\fakelanguage}
\DeclareTextFontCommand{\textcyr}{\cyr}

\numberwithin{equation}{section}

\usepackage{calligra}

\DeclareMathOperator{\HOM}{\mathscr{H}\text{\kern -3pt {\calligra\large om}}\,}

\makeatletter
\newsavebox{\@brx}
\newcommand{\llangle}[1][]{\savebox{\@brx}{\(\m@th{#1\langle}\)}%
  \mathopen{\copy\@brx\kern-0.5\wd\@brx\usebox{\@brx}}}
\newcommand{\rrangle}[1][]{\savebox{\@brx}{\(\m@th{#1\rangle}\)}%
  \mathclose{\copy\@brx\kern-0.5\wd\@brx\usebox{\@brx}}}
\makeatother

\newcommand{\bsh}{\begin{shaded}}
\newcommand{\esh}{\end{shaded}}

\newcommand{\ic}{\sqrt{-1}}

\newcommand{\otm}{\otimes}

\newcommand\xqed[1]{%
  \leavevmode\unskip\penalty9999 \hbox{}\nobreak\hfill
  \quad\hbox{#1}}
\newcommand\qetr{\xqed{$\triangle$}}

\begin{document}
\title[Satake correspondence for equivariant qDE and qKZ of Grassmannians ]{ On the Satake correspondence for the equivariant quantum differential equations and qKZ difference equations of Grassmannians}
\author[Giordano Cotti and Alexander Varchenko]{Giordano Cotti$\>^{\circ,\bullet}$ and Alexander Varchenko$\>^\star$}
{\let\thefootnote\relax
\footnotetext{\vskip5pt 
\noindent
$^\circ\>$\textit{ E-mail}:  gcotti@fc.ul.pt
\vskip2pt
\noindent
$^\star\>$\textit{E-mail}: anv@email.unc.edu }
\vskip2mm}
\email{gcotti@fc.ul.pt, anv@email.unc.edu}
\maketitle
\begin{center}
\textit{ 
$^\circ\>$Departamento de Matem\'atica\\ Faculdade de Ci\^encias da Universidade de Lisboa\\ 
Campo Grande Edif\'icio C6, 1749-016 Lisboa, Portugal\/\\
\vskip4pt
$^\bullet\>$Grupo de F\'isica Matem\'atica \\
Departamento de Matem\'atica, Instituto Superior Técnico\\
Av. Rovisco Pais, 1049-001 Lisboa, Portugal\/\\}

\vskip4pt
\textit{ $^{\star}\>$Department of Mathematics, University
of North Carolina at Chapel Hill\\ Chapel Hill, NC 27599-3250, USA\/}

\end{center}

\vskip1cm
\begin{center}
    \textit{In memory of Yuri Ivanovich Manin (1937-2023)}
\end{center}
\vskip1cm

\begin{abstract}
We consider the joint system of equivariant quantum differential equations (qDE) and qKZ difference equations for the Grassmannian $G(k,n)$, which parametrizes $k$-dimensional subspaces of $\C^n$. First, we establish a connection between this joint system for $G(k,n)$ and the corresponding system for the projective space $\Pb^{n-1}$. Specifically, we show that, under suitable {\it Satake identifications} of the equivariant cohomologies of $G(k,n)$ and $\Pb^{n-1}$, the joint system for $G(k,n)$ is gauge equivalent to a differential-difference system on the $k$-th exterior power of the cohomology of $\Pb^{n-1}$.

Secondly, we demonstrate that the \textcyr{B}-theorem for Grassmannians, as stated in \cite{CV21,TV23}, is compatible with the Satake identification. This implies that the \textcyr{B}-theorem for $\Pb^{n-1}$ extends to $G(k,n)$ through the Satake identification. As a consequence, we derive determinantal formulas and new integral representations for multi-dimensional hypergeometric solutions of the joint qDE and qKZ system for $G(k,n)$.

Finally, we analyze the Stokes phenomenon for the joint system of qDE and qKZ equations associated with $G(k,n)$. We prove that the Stokes bases of solutions correspond to explicit $K$-theoretical classes of full exceptional collections in the derived category of equivariant coherent sheaves on $G(k,n)$. Furthermore, we show that the Stokes matrices equal the Gram matrices of the equivariant Euler--Poincaré--Grothendieck pairing with respect to these exceptional $K$-theoretical bases.
\end{abstract}
\vskip0,5cm
\begin{adjustwidth}{35pt}{35pt}
\noindent{\small {\it Key words:} Grassmannian, Gromov--Witten invariants, quantum differential equations, qKZ difference equations, Satake correspondence, derived category, exceptional collections}
\vskip2mm
\noindent{\small {\it 2020 Mathematics Subject Classification:} Primary 14F05, 14N35, 34M40; Secondary: 34M35}
\end{adjustwidth}
\tableofcontents

\section{Introduction}
\noindent 1.1.\,\,{\bf Equivariant quantum cohomology.} Gromov--Witten theory associates to any smooth complex projective variety \( X \) a sequence of numerical invariants, defined as counts of curves on \( X \) with fixed genus and degree that satisfy suitable incidence conditions. Following A.\,Givental \cite{giv1}, when \( X \) is equipped with the action of a complex reductive group \( \mathsf{G} \), these curve counts can be suitably decorated in an 
equivariant manner, and the Gromov--Witten invariants are promoted to classes in the equivariant cohomology algebra of a point. 
See \cite{cox,mirror} for an introduction.
\vskip1.5mm
Remarkably, the genus zero (equivariant) Gromov--Witten invariants of \( X \) can then be used to introduce a family of deformations of the classical (equivariant) cohomology algebra of \( X \): the classical (equivariant) cup product is deformed into a quantum (equivariant) product. This deformation is achieved by adding correction terms given by suitable generating functions of genus zero Gromov--Witten invariants. The resulting family of algebras is called the ({\it equivariant}) {\it quantum cohomology} of \( X \). This is a prototypical example of a {\it formal Frobenius manifold} --in the sense of Yu.I.\,Manin \cite{manin}-- over the (equivariant) cohomology algebra of a point.
\vskip1.5mm
In this paper, we focus on the case in which \( X \) is a Grassmannian \( G(k,n) \), whose points parametrize \( k \)-dimensional subspaces of \( \mathbb{C}^n \), and the group \( \mathsf{G} \) is the \( n \)-dimensional torus \( T = (\mathbb{C}^*)^n \) acting on the ambient space \( \mathbb{C}^n \). See \cite{Ber97,Buc03,BKT03,GRTV13,Kim96,KT03,Mih06,Rie01,ST97}.
\vskip2mm

\noindent 1.2.\,\,{\bf Equivariant qDE and qKZ equations.} This work is devoted to the study of the equivariant quantum differential equations (qDE) and qKZ difference equations for the Grassmannians \( G(k,n) \).
\vskip1.5mm
The equivariant qDE for \( G(k,n) \) consist of a system of partial differential equations given by
\beq\label{intro.eq1}
\kappa\, p_1 \frac{\partial}{\partial p_1} f(\bm{z}; \bm{p}; \kappa) = c_1^T(E_1) *_{\bm{p}, \bm{z}} f(\bm{z}; \bm{p}; \kappa), \qquad \kappa\, p_2 \frac{\partial}{\partial p_2} f(\bm{z}; \bm{p}; \kappa) = c_1^T(E_2) *_{\bm{p}, \bm{z}} f(\bm{z}; \bm{p}; \kappa),
\eeq
where the unknown function \( f(\bm{z}; \bm{p}; \kappa) \) takes values in the equivariant cohomology \( H^\bullet_T(G(k,n), \mathbb{C}) \), and the operators \( c_1^T(E_1) *_{\bm{p}, \bm{z}},\, c_1^T(E_2) *_{\bm{p}, \bm{z}} \colon H^\bullet_T(G(k,n), \mathbb{C}) \to H^\bullet_T(G(k,n), \mathbb{C}) \) denote quantum multiplication by the equivariant first Chern classes of the tautological bundle \( E_1 \) and the quotient bundle \( E_2 \) on \( G(k,n) \), respectively. The qDE depend on the equivariant parameters\footnote{The variables $\bm z$ are the equivariant Chern roots of the trivial bundle $\underline{\C^n}=\C^n\times G(k,n)$. As such, the variables $\bm z$ are elements of $H^\bullet_T(G(k,n),\C)$. It turns out that $H^\bullet_T(G(k,n),\C)$ has a natural structure of $\C[\bm z]$-module.} \( \bm{z} = (z_1, \dots, z_n) \) corresponding to the components of the torus \( T \), as well as on an additional parameter \( \kappa \in \mathbb{C} \).
\vskip1.5mm
The unknown function \( f(\bm{z}; \bm{p}; \kappa) \) must also satisfy a further set of {\it difference} equations, the {\it quantized Knizhnik-Zamolodchikov} (qKZ) equations. These are of the form
\beq\label{intro.eq2}
f(z_1, \dots, z_a + \kappa, \dots, z_n; \bm{p}; \kappa) = K_a(\bm{z}; \bm{p}; \kappa) f(\bm{z}; \bm{p}; \kappa), \quad a = 1, \dots, n,
\eeq
for suitable operators \( K_a \). The qDE and qKZ equations were originally introduced in \cite{TV23} for partial flag varieties, and they define a compatible joint system of differential-difference equations. 

\begin{rem}\label{intro.rem1}Introduce the differential operators $\mathsf D_{i,\bm p,\bm z,\kappa}:=\kappa\,p_i\frac{\der}{\der p_i}-c_1^T(E_i)*_{\bm p,\bm z}$ for $i=1,2$ and the difference operators ${\mathsf Z}_{a,\bm p,\bm z,\kappa}=K_a(\bm z;\bm p;\kappa)^{-1}\on{T}_a$ for $a=1,\dots,n$, where $\on{T}_a=\exp\left(\kappa\frac{\der}{\der z_a}\right)$ acts on functions $f(\bm z;\bm p;\kappa)$ by $f(\bm z;\bm p;\kappa)\mapsto f(z_1,\dots z_a+\kappa,\dots,z_n;\bm p;\kappa)$. The qDE and qKZ equations take then the more compact form $\mathsf D_{i,\bm p,\bm z,\kappa}f=0$ and ${\mathsf Z}_{a,\bm p,\bm z,\kappa}f=f$. The precise meaning of the compatibility of the joint system of qDE and qKZ equations is the following: all the operators $\mathsf D_{i,\bm p,\bm z,\kappa}, {\mathsf Z}_{a,\bm p,\bm z,\kappa}$, with $i=1,2$ and $a=1,\dots,n$, pairwise commute \cite{TV23}.
\end{rem}

\begin{rem}
The qKZ operators discussed above are introduced in \cite{TV23}, with their definition originating from representation theory, as we will explain below. Difference operators in equivariant quantum cohomology, which have a Gromov–Witten theoretical definition, have already appeared in the literature under the names {\it shift} (or {\it Seidel}) {\it operators}, see \cite{OP10,BMO,mo}. Additionally, the interesting paper \cite{Iri17} provides a general definition of these difference operators for the big equivariant quantum cohomology. It would be interesting to establish precise relationships between H.\,Iritani's difference operators in \cite{Iri17} and our qKZ operators $K_a$.
\end{rem}
\vskip1,5mm
\noindent 1.3.\,\,{\bf Previous results.} Before presenting the new results of this paper, let us briefly review the previously known results concerning the joint system of qDE and qKZ equations. In \cite{TV23}, the more general case of qDE and qKZ equations for partial flag varieties was studied, and a basis of solutions for the joint system was constructed in the form of multidimensional hypergeometric functions. In this approach, the solutions are labeled by elements of the \( K \)-theory algebra of the partial flag variety.
\vskip1,5mm
The works \cite{TV21,CV21} focused on the analysis of qDE and qKZ equations for projective spaces. In this particular case, the asymptotic properties of solutions to the joint system were examined, and the Stokes phenomenon was investigated. Initially, in \cite{TV21} the Stokes bases of the qDE at its irregular singular point were described in terms of exceptional bases of the equivariant \( K \)-theory and a suitable braid group action on the set of exceptional bases. Subsequently, in \cite{CV21} the Stokes bases were identified with explicit $T$-full exceptional collections 
in the derived category $\mathsf D^b_{T}(\mathbb P^{n-1})$
of $T$-equivariant coherent sheaves on $\mathbb P^{n-1}$. The resulting exceptional collections are explicit mutations of the A.\,Beilinson's exceptional collection $\left(\mc O,\mc O(1),\dots,\mc O(n-1)\right)$ in $\mathsf D^b_{T}(\mathbb P^{n-1})$, see \cite{Bei78,CV21}.
\vskip1,5mm
Furthermore, in \cite{CV21} relations were established between the monodromy data of the joint qDE and qKZ system for \( \mathbb{P}^{n-1} \) and the characteristic classes of objects in the derived category \( \mathsf{D}^b_T(\mathbb{P}^{n-1}) \) of equivariant coherent sheaves on \( \mathbb{P}^{n-1} \). For example, it was shown that the Stokes matrices connecting different Stokes bases in contiguous Stokes sectors equal the Gram matrices of the Euler--Poincar\'e--Grothendieck pairing with respect to the associated exceptional bases in $K$-theory.
\vskip1,5mm
In addition, in \cite[\textcyr{B}-Theorem]{CV21} a remarkable commutative diagram was obtained: 
\beq\label{intro.diag1}
\xymatrix{
\parbox{4cm}{\begin{center}Equivariant $K$-theory of $\Pb^{n-1}$\end{center}}\ar^{\quad\quad\textcyr{B}_\Pb\quad\quad\quad}[rr]\ar_{\mu_\Pb}[dr]&&
        \parbox{4cm}{\begin{center}Equivariant cohomology of $\Pb^{n-1}$\end{center}}
\ar^{\mathcal S_\Pb}[dl]\\
&\parbox{4cm}{\begin{center}Space of solutions of qDE of $\Pb^{n-1}$\end{center}}&
}
\eeq
The map \(\mu_\Pb\) introduced in \cite{TV21} associates to a \( K \)-theoretical class \( [E] \in K^T_0(\Pb^{n-1}) \) a solution of the joint system of the qDE and qKZ equations for \(\Pb^{n-1}\). On the other hand, the map \textcyr{B}$_\Pb$ assigns to a \( K \)-theoretical class \( [E] \) a distinguished characteristic class, defined as a twisted version of the \(\kappa\)-{\it graded equivariant Chern character} \(\on{Ch}_T^{(\kappa)}(E)\). This is the characteristic class given by
\[
{\rm Ch}^{(\kappa)}_T(E) := \sum_{i=1}^r \exp\left(-\frac{2\pi\sqrt{-1}}{\kappa} \xi_i\right),\quad \text{where \( \xi_1, \dots, \xi_r \) are the equivariant Chern roots of \( E \).}
\]
Introduce the characteristic class \( \Hat{F}_{\Pb} \) of the tangent bundle \( T\Pb^{n-1} \) by the product
\[
\Hat{F}_{\Pb} := \prod_{j=1}^{\binom{n}{k}} F(\delta_j),\quad\text{with \( \delta_1, \dots, \delta_{\binom{n}{k}} \) being the equivariant Chern roots of \( T\Pb^{n-1} \),}
\]and where \( F(t) \) is a formal power series given by the  Taylor expansion at \( t = 0 \) of the function
\begin{multline*}
F(t) = \exp\left(-\frac{\log \kappa}{\kappa} t\right) \Gamma\left(1 - \frac{t}{\kappa}\right)\\ = 1 + \frac{t (\gamma_{\rm EM} - \log \kappa)}{\kappa} + \frac{t^2 \left(6 \log^2 \kappa - 12 \gamma_{\rm EM} \log \kappa + 6 \gamma_{\rm EM}^2 + \pi^2\right)}{12 \kappa^2} + O\left(t^3\right).
\end{multline*}
In the above formula, \( \Gamma(t) \) denotes the Euler Gamma function, and \( \gamma_{\rm EM} := \lim_{n \to \infty} \left(\frac{1}{n} - \sum_{k=1}^n \frac{1}{k}\right) = 0.57721 \dots \) is the Euler–Mascheroni constant.
The specific twist that defines the \textcyr{B}-class \(\text{\textcyr{B}}_\Pb(E)\) is
\[
\text{\textcyr{B}}_\Pb(E) = \Hat{F}_{\Pb} \on{Ch}_T^{(\kappa)}(E). 
\]
Finally, the map \( \mc{S}_\Pb \) is the \textit{topological-enumerative morphism}, which assigns a solution of the qDE (only) to each cohomology class in \( H^\bullet_T(\Pb^{n-1}, \C) \). This solution admits a series expansion in \( q = \frac{p_2}{p_1} \) that serves as a generating function for genus-zero equivariant Gromov--Witten invariants with gravitational descendants, i.e., insertions of \( \psi \)-classes on the moduli space of stable maps. Precise definitions will be provided in Section \ref{sectopenmor}.
\vskip1.5mm
In the subsequent paper \cite{TV23}, the validity of the \textcyr{B}-Theorem was established for all partial flag varieties. In the specific case of Grassmannians, the \textcyr{B}-Theorem represents an equivariant analog of the results in \cite{GGI16} and \cite{CDG1}. Additionally, it can be viewed as a refinement of the Gamma Theorem from \cite{tar-var,TV21}.
\begin{rem}
In the formulation of the \textcyr{B}-Theorem in \cite{TV23}, the topological-enumerative morphism is replaced by the {\it Levelt fundamental solution}, which is uniquely defined by its analytical properties. The same work also conjectures that the Levelt fundamental solution is equal to the topological-enumerative morphism \cite[Conjectures 6.31, 6.32]{TV23}. In this paper, we prove this conjecture in the case of $G(k,n)$, see Theorem \ref{sollev}. The same arguments can be easily extended to the more general case of partial flag varieties. 
\end{rem}

\vskip2mm
\noindent 1.4.\,\,{\bf Main results.} We establish relations between the joint system of qDE and qKZ equations for $G(k,n)$ and the corresponding join system for the projective space $\Pb^{n-1}$, as we are going to explain.
\vskip1,5mm
The joint system \eqref{intro.eq1}, \eqref{intro.eq2} is a set of differential-difference equations for sections $f(\bm z;\bm p;\kappa)$ of the trivial bundle $H_{k,n}$ over $\C^n\times\C^2$, whose fiber over $(\bm z_o;\bm p_o)$ is the specialization $H^\bullet_T(G(k,n),\C)_{\bm z_o}$ of the equivariant cohomology at $\bm z=\bm z_o$,  that is the quotient $H^\bullet_T(G(k,n),\C)/\langle \bm z-\bm z_o\rangle$. 

Any tensor powers of the bundle $H_{k,n}$ naturally inherits an analogous joint system of differential-difference equations, induced by \eqref{intro.eq1}, \eqref{intro.eq2}. Namely, for sections $f(\bm z;\bm p;\kappa)$ of the tensor power $H_{k,n}^{\otimes m}$, we have the induced joint system
\begin{align*}\Tilde{\mathsf D}_{i,\bm p,\bm z,\kappa}f(\bm z;\bm p;\kappa)&=0,&i=1,2,\\
 \Tilde{\mathsf Z}_{a,\bm p,\bm z,\kappa}f(\bm z;\bm p;\kappa)&=f,&a=1,\dots,n,
\end{align*}
where the operators $\Tilde{\mathsf D}_{i,\bm p,\bm z,\kappa},\Tilde{\mathsf Z}_{a,\bm p,\bm z,\kappa}$ act on decomposable sections as follows:
\beq
\Tilde{\mathsf D}_{i,\bm p,\bm z,\kappa}\left[\bigotimes_{j=1}^m f_j(\bm z;\bm p;\kappa)\right]=\sum_{j=1}^mf_1(\bm z;\bm p;\kappa)\otimes\dots\otimes{\mathsf D}_{i,\bm p,\bm z,\kappa}f_j(\bm z;\bm p;\kappa) \otimes\dots\otimes f_m(\bm z;\bm p;\kappa)
\eeq
\beq
\Tilde{\mathsf Z}_{a,\bm p,\bm z,\kappa}\left[\bigotimes_{j=1}^m f_j(\bm z;\bm p;\kappa)\right]=\bigotimes_{j=1}^m {\mathsf Z}_{a,\bm p,\bm z,\kappa}f_j(\bm z;\bm p;\kappa),
\eeq
the operators ${\mathsf D}_{i,\bm p,\bm z,\kappa}$ and ${\mathsf Z}_{a,\bm p,\bm z,\kappa}$ being defined as in Remark \ref{intro.rem1}. Similar equations can be written for symmetric and exterior powers of $H_{k,n}$.
\vskip1,5mm
Our first main result, Theorem \ref{mainthm1}, establishes the gauge equivalence between the joint system of qDE and qKZ equations for sections of \( H_{k,n} \) and the system induced for sections of the exterior power \( \bigwedge^k H_{1,n} \). In particular, this implies that the joint system of qDE and qKZ equations for the Grassmannian \( G(k,n) \) can be explicitly solved by starting from \( k \)-tuples of solutions to the simpler joint system of qDE and qKZ equations for \( \mathbb{P}^{n-1} \). 

This gauge equivalence relies on a suitable identification between the generic fiber of \( H_{k,n} \) (i.e., \( H^\bullet_T(G(k,n), \mathbb{C}) \)) and the \( k \)-th exterior power of the generic fiber of \( H_{1,n} \) (i.e., \( H^\bullet_T(\mathbb{P}^{n-1}, \mathbb{C}) \)). Moreover, it involves identifying the operators of quantum multiplication on \( H^\bullet_T(G(k,n), \mathbb{C}) \) with operators on \( \bigwedge^k_{\C[\bm z]} H^\bullet_T(\mathbb{P}^{n-1}, \mathbb{C}) \) induced by quantum multiplication operators on \( H^\bullet_T(\mathbb{P}^{n-1}, \mathbb{C}) \), Theorem \ref{gSc}. The validity of these identifications is not new; they follow from the {\it geometric Satake correspondence}, in both its classical formulation and its equivariant quantum partial extensions. What is novel and non-trivial is the identification of the qKZ operators for \( G(k,n) \) with the operators on \( \bigwedge^kH_{1,n}\) induced by the qKZ operators for \( \mathbb{P}^{n-1} \).

Our proof proceeds through several steps, relying on constructions from representation theory. Even for the known results, our approach differs from existing ones. 
The strategy can be summarized as follows:

\begin{enumerate}
    \item We begin by considering an ``abstract'' joint system of differential-difference equations, referred to as the {\it dynamical} and {\it qKZ equations}, for sections of the trivial bundle \(
    (\mathbb{C}^2)^{\otimes n} \times \mathbb{C}^n \times \mathbb{C}^2\) over \(\mathbb{C}^n \times \mathbb{C}^2\). This system is inherently tied to representation theory and can be viewed as the limit of a similar system defined using the Yangian \(\mathsf{Y}_h(\mathfrak{gl}_2)\) action on \((\mathbb{C}^2)^{\otimes n}\), in the limit \(h \to \infty\) for the Yangian deformation parameter.
    
    \item Next, by decomposing the space \((\mathbb{C}^2)^{\otimes n} = \bigoplus_{k=0}^n (\mathbb{C}^2)^{\otimes n}_{k, n-k}\), the joint system of dynamical and qKZ equations can be broken down into systems for sections of each vector subbundle \( U_{k,n} = (\mathbb{C}^2)^{\otimes n}_{k, n-k} \times \mathbb{C}^n \times \mathbb{C}^2 \) over \(\mathbb{C}^n \times \mathbb{C}^2\). Then, under a suitable bundle isomorphism \( U_{k,n} \cong \bigwedge^k U_{1,n} \),  we show that the joint system for sections of \( U_{k,n} \) is gauge equivalent to the induced system for sections of \(\bigwedge^k U_{1,n}\), Theorem \ref{abstractthmsat}. This is done by inspection of the explicit formulas defining the dynamical and qKZ operators. 
    
    \item The geometric interpretation in terms of quantum cohomology is established through the introduction of appropriate {\it stable envelope maps} \cite{mo,RTV15,TV23}. These maps define isomorphisms of bundles, \({\rm Stab}_{k,n}^\diamond\colon U_{k,n}\to H_{k,n}\), which intertwine the dynamical operators with the quantum multiplication operators \( c_1^T(E_i) *_{\bm{p}, \bm{z}} \), for \( i = 1, 2 \) (Theorem \ref{stabint}), and the abstract representation-theoretic qKZ operators with those acting on sections of \( H_{k,n} \).
\end{enumerate}
\vskip1,5mm
Our second main result, Theorem \ref{thmcompatibility}, is the compatibility of the \textcyr{B}-theorem for Grassmannians with the Satake correspondence. More precisely, we provide suitable identifications of the equivariant $K$-theory, equivariant cohomology, and the space of solutions of the qDE associated with $G(k,n)$ and the $k$-th exterior powers of the corresponding spaces for $\Pb^{n-1}$. This leads to the commutative diagram
\beq\label{intro.diag2}
\xymatrix @C=1pc @R=.5pc{
\parbox{4cm}{\begin{center}$k$-th exterior power of equivariant $K$-theory of $\Pb^{n-1}$\end{center}}\ar[rr]\ar[dr]\ar[dd]&&\parbox{4cm}{\begin{center}$k$-th exterior power of equivariant cohomology of $\Pb^{n-1}$\end{center}}\ar[dl]\ar[dd]\\
&\parbox{4cm}{\begin{center}$k$-th exterior power of space of solutions of qDE of $\Pb^{n-1}$\end{center}}\ar[dd]&\\
\parbox{4cm}{\begin{center}Equivariant $K$-theory of $G(k,n)$\end{center}}\ar[dr]\ar[rr]&&\parbox{4cm}{\begin{center}Equivariant cohomology of $G(k,n)$\end{center}}\ar[dl]\\
&\parbox{4cm}{\begin{center}Space of solutions of qDE of $G(k,n)$\end{center}}&
}
\eeq
The top face of the diagram is induced by the commutative diagram \eqref{intro.diag1} of the \textcyr{B}-theorem for $\Pb^{n-1}$. The bottom face represents the commutative diagram asserted by the \textcyr{B}-theorem for $G(k,n)$. The vertical arrows denote the suitable identifications mentioned above. The commutativity of diagram \eqref{intro.diag2} is rich of consequences. For example, from it we deduce both a determinantal formula (Theorem \ref{detprop}) and a non-trivial integral representation for solutions of the qDE and qKZ system of $G(k,n)$ (Theorem \ref{intrepthm}).

\vskip1,5mm
Our third main result, Theorem \ref{mainthm3}, establishes a connection between the Stokes bases of solutions for the joint system of qDE and qKZ equations for \( G(k,n) \) and the \( T \)-full exceptional collections in the derived category \(\mathsf{D}^b_T(G(k,n))\) of \( T \)-coherent sheaves on \( G(k,n) \). We demonstrate that these exceptional collections are braid mutations of the twisted M.\,Kapranov's collection \(\left(\mathscr{L} \otimes \Sigma_{\bm{\lambda}} E_1^*\right)_{\bm{\lambda} \subseteq k \times (n-k)}\), where \(\mathscr{L}\) is the line bundle \(\det\bigwedge^2 E_1^*\) on \( G(k,n) \), and \(\Sigma_{\bm{\lambda}}\) denotes the Schur functor associated with the Young diagram \(\bm{\lambda}\), see \cite{Kap88}. Furthermore, the corresponding $K$-theoretical exceptional bases are explicitly described, Theorem \ref{thmqstok}.

Additionally, we show that the Stokes matrices, which relate Stokes bases of solutions, coincide with the Gram matrices of the Euler--Poincaré--Grothendieck product \(\chi^T\) in \( K \)-theory, relative to the corresponding exceptional bases, Theorem \ref{corstokmatr}. Specifically, let \(\mathcal{V}_1\) be a Stokes sector and \(\mathcal{V}_2 = e^{\pi\sqrt{-1}}\mathcal{V}_1\) its complementary sector. Let \(Y_1\) and \(Y_2\) be the Stokes bases on \(\mathcal{V}_1\) and \(\mathcal{V}_2\), respectively, with the base change matrix \(\mathbb{S}\) defined by \(Y_2 = Y_1 \mathbb{S}\). If \(\epsilon_1\) and \(\epsilon_2\) are the exceptional bases corresponding to \(Y_1\) and \(Y_2\), we prove that \(\epsilon_2\) is the \textit{left-dual} of \(\epsilon_1\) (Theorem \ref{thmleftdual}), in the sense of \cite{GK04,CDG1,CV21}. This directly implies that the entries of \(\mathbb{S}\) are the \(\chi^T\)-products of the elements of the \(\epsilon_2\)-basis, specifically \(\chi^T(\epsilon_{2,i}, \epsilon_{2,j})\), for \(i, j = 1, \dots, \binom{n}{k}\).

Finally, it's noteworthy that the identification of Stokes bases with exceptional collections, and Stokes matrices with Gram matrices, is compatible with the Satake correspondence. Under the identification of \(K^T_0(G(k,n))\) with the \(k\)-th exterior power of \(K^T_0(\mathbb{P}^{n-1})\), the exceptional bases that label the Stokes bases for \(G(k,n)\) correspond to \(k\)-fold wedge products of the exceptional bases labeling the Stokes bases for \(\mathbb{P}^{n-1}\), Theorem \ref{satstoksol}. Consequently, the Stokes matrices for \(G(k,n)\) are given by the matrices of \(k\)-minors of the Stokes matrices for \(\mathbb{P}^{n-1}\), Corollary \ref{corstokmat}.
\vskip2mm
\noindent 1.5.\,\,{\bf Relations with Maulik--Okounkov theory.} 
This paper is connected to the work of D.\,Maulik and A.\,Okounkov on the equivariant quantum cohomology of Nakajima quiver varieties \cite{mo}. In their framework, it was proved that the quantum differential equations associated with Nakajima quiver varieties admit a compatible system of difference equations. The coefficients of these difference equations, called {\it shift operators}, are expressed in terms of Gromov--Witten theory \cite[Ch.,8]{mo}. Moreover, it was established that this system of difference equations corresponds to the {\it quantized Knizhnik-Zamolodchikov equations} of I.B.\,Frenkel and N.Y.\,Reshetikhin \cite{FR92}, as discussed in \cite[Ch.,9]{mo}.
\vskip1,5mm

One of the simplest examples of Nakajima quiver varieties is the cotangent bundle $T^*\mc F_{\bm\la}$ of partial flag varieties\footnote{Partial flag varieties themselves are not Nakajima varieties.} $\mc F_{\bm\la}$. For Grassmannians, this corresponds to the most basic quiver, consisting of a single vertex without arrows. In \cite{GRTV13,RTV15}, the combined system of qDE and qKZ equations for $T^*\mc F_{\bm\la}$ was studied. Explicit formulas were obtained, identifying the equivariant qDE and qKZ equations with the dynamical differential equations and qKZ difference equations in representation theory, associated with the evaluation module $\mathbb{C}^N(z_1) \otimes \cdots \otimes \mathbb{C}^N(z_n)$ of the Yangian ${\mathsf Y}_h(\mathfrak{gl}_N)$. 
\vskip 1.5mm
It was expected that taking an appropriate limit as $h \to \infty$ for these equations would produce a system of difference and differential equations related to the equivariant quantum cohomology of the partial flag varieties themselves, rather than their cotangent bundles. In particular, \cite{BMO} showed that the $h \to \infty$ limit of the equivariant qDE for the cotangent bundle of the full flag variety corresponds to the equivariant qDE of the full flag variety itself. Additionally, a suitable $h \to \infty$ limit of the Yangian $R$-matrix, associated with the quantum cohomology of the cotangent bundle of a Grassmannian, was used in \cite{GK,GKS} to compute the quantum multiplication in the cohomology of the Grassmannian.
\vskip1,5mm
This expectation was confirmed in \cite{TV23}, where the ``abstract representation-theoretic'' limiting dynamical differential equation was shown to correspond to the equivariant quantum differential equations for partial flag varieties. Under this identification, the limiting qKZ difference equations become a system of difference equations in the equivariant cohomology of partial flag varieties, forming a new system\footnote{In the case of Grassmannians these are equations \eqref{intro.eq2}.} compatible with the equivariant quantum differential equations. Furthermore, under this identification, the multidimensional hypergeometric solutions to the limiting dynamical and qKZ equations constructed in \cite{TV23} (as limits for $h \to \infty$ of those constructed in \cite{TV94,tar-var}) provide solutions of the equivariant quantum differential equations and qKZ difference equations for partial flag varieties.
\vskip2mm
\noindent 1.6.\,\,{\bf Affine Grassmannian and geometric Satake correspondence.} 
We conclude this introduction with a self-contained explanation of the role played in this paper by the geometric Satake correspondence \cite{BD96,Lus83,Gin90,Gin95,Gin08,MV00,MV07}. This will also clarify the seemingly ``mysterious'' role of the exterior power functors, which enable us to reduce the study of the qDE and qKZ equations for \( G(k,n) \) to the simpler system associated with \( \mathbb{P}^{n-1} \).
\vskip1,5mm
The fundamental reason for the simplifying effect of the exterior powers lies in the geometry of the affine Grassmannian of \( GL(n, \mathbb{C}) \). Given a reductive algebraic group $\mathsf G$ over $\C$, the affine Grassmannian $\mathscr{G}\!r(\mathsf G)$ is an ind-scheme, that is an infinite dimensional space which can be realized as inductive limit of schemes. For sake of simplicity, $\mathscr{G}\!r(\mathsf G)$ can be thought as an algebraic version of a loop group: the space $\mathscr{G}\!r(\mathsf G)$ can indeed be identified with the loop group $\Om\mathsf G_c$, where $\mathsf G_c$ is the maximal compact subgroup of $\mathsf G$. As such, the affine Grassmannian $\mathscr{G}\!r(\mathsf G)$ is equipped with a topological group structure.
By fixing a Borel subgroup \( \mathsf{B} \subset \mathsf{G} \), we can introduce the abstract Cartan torus \( \mathsf{T} \) (the quotient of \( \mathsf{B} \) by its unipotent radical), the group of weights \( \mathbb{X}^\bullet(\mathsf{T}) = \operatorname{Hom}(\mathsf{T}, \mathbb{C}^*) \), and coweights \( \mathbb{X}_\bullet(\mathsf{T}) = \operatorname{Hom}(\mathbb{C}^*, \mathsf{T}) \). The set of roots is denoted by \( \Phi \subset \mathbb{X}^\bullet(\mathsf{T}) \), and the set of coroots by \( \Phi^\vee \subset \mathbb{X}_\bullet(\mathsf{T}) \). The {\it root datum} of \( \mathsf{G} \) is given by the tuple \( (\mathbb{X}^\bullet(\mathsf{T}), \mathbb{X}_\bullet(\mathsf{T}), \Phi, \Phi^\vee) \). The choice of \( \mathsf{B} \) also determines the semigroup of dominant coweights \( \mathbb{X}_\bullet(\mathsf{T})^+ \subset \mathbb{X}_\bullet(\mathsf{T}) \) and the set of positive roots \( \Phi^+ \subset \Phi \).
\vskip1,5mm
The space \( \mathscr{G}\!r(\mathsf{G}) \) admits a decomposition into strata \( S_\lambda \), called {\it (spherical)} {\it Schubert varieties}, which are in one-to-one correspondence with dominant coweights \( \lambda \in \mathbb{X}_\bullet(\mathsf{T})^+ \). Each stratum \( S_\lambda \) is a projective variety, which is typically singular. It turns out that \( S_\lambda \) is smooth if and only if \( \lambda \) is a {\it minuscule} dominant coweight, i.e., such that \( \langle \alpha, \lambda \rangle \leq 1 \) for any positive root \( \alpha \in \Phi^+ \). In such cases, we have \( S_\lambda = \mathsf{G}/P_\lambda \), where \( P_\lambda \) is the parabolic subgroup of \( \mathsf{G} \) associated with \( \lambda \). This class of varieties includes Grassmannians \( G(k,n) \), orthogonal Grassmannians \( OG(n, 2n) \), even-dimensional quadrics \( Q^{2n} \), the Cayley plane \( \mathbb{OP}^2 = E_6/Q_1 \), and the Freudenthal variety \( E_7/Q_7 \). See \cite{LMS79}.
\vskip1,5mm
Intersection cohomology provides a well-behaved cohomology theory for singular and stratified spaces \cite{KW06}. Applying this to \( \mathscr{G}\!r(\mathsf{G}) \), one associates with each Schubert variety \( S_\lambda \) an {\it intersection cohomology sheaf} \( \operatorname{IC}(S_\lambda) \) as an object in a suitable derived category \( \mathcal{D} \) of constructible sheaves on \( \mathscr{G}\!r(\mathsf{G}) \). The complex \( \operatorname{IC}(S_\lambda) \) encodes cohomological information about \( S_\lambda \). Specifically, its hypercohomology gives the intersection cohomology groups of \( S_\lambda \). Furthermore, one can consider the smallest subcategory \( \mathsf{Perv}(\mathsf{G}) \) of \( \mathcal{D} \) generated by objects isomorphic to finite sums of intersection cohomology sheaves \( \operatorname{IC}(S_\lambda) \). Objects in \( \mathsf{Perv}(\mathsf{G}) \) are known as {\it perverse sheaves} on \( \mathscr{G}\!r(\mathsf{G}) \).
The topological group structure on  \( \mathscr{G}\!r(\mathsf{G}) \) induces a convolution product $\circledast$ on the objects of \( \mathsf{Perv}(\mathsf{G}) \).
\vskip1,5mm
The central result of the geometric Satake correspondence is that \((\mathsf{Perv}(\mathsf{G}), \circledast)\) forms a \textit{neutral Tannakian category} \cite{SR72}, meaning it is equivalent to the category of finite-dimensional representations of an affine group scheme \(\mathsf{G}^*\). More precisely, there exists an equivalence of tensor categories
\begin{equation}\label{equivsat}
(\mathsf{Perv}(\mathsf{G}), \circledast) \cong (\operatorname{Rep}_\mathbb{C}(\mathsf{G}^*), \otimes),
\end{equation}
for a suitable affine group scheme \(\mathsf{G}^*\), see \cite{Gin95,MV00,MV07}. It turns out that the Tannakian group \(\mathsf{G}^*\) has a root datum \((\mathbb{X}_\bullet(\mathsf{T}), \mathbb{X}^\bullet(\mathsf{T}), \Phi^\vee, \Phi)\), making it isomorphic to the {\it Langlands dual} group \(\mathsf{G}^\vee\). In this sense, the geometric Satake correspondence provides an intrinsic construction of \(\mathsf{G}^\vee\) that does not rely on the combinatorial data of the root datum but instead uses the classification theorem for reductive groups.

Under the equivalence \eqref{equivsat}, the intersection cohomology sheaf \(\operatorname{IC}(S_\lambda)\) corresponds to the irreducible representation \(V_\lambda\) of \(\mathsf{G}^* = \mathsf{G}^\vee\) with highest weight \(\lambda \in \mathbb{X}_\bullet(\mathsf{T})^+\). Furthermore, the hypercohomology functor \(\mathbb{H}^\bullet\colon \mathsf{Perv}(\mathsf{G}) \to \operatorname{Vect}_\mathbb{C}\) is identified with the forgetful functor \(\operatorname{Rep}_\mathbb{C}(\mathsf{G}^*) \to \operatorname{Vect}_\mathbb{C}\).
\vskip1,5mm
The group \(\mathsf{G} = GL(n, \mathbb{C})\) is self-dual under Langlands duality, and the minuscule Schubert varieties \(S_\lambda\) of the affine Grassmannian \(\mathscr{G}\!r(\mathsf{G})\) are precisely the Grassmannians \(G(k,n)\) for \(k = 1, \ldots, n\). Through the equivalence \eqref{equivsat}, the cohomology \(H^\bullet(S_\lambda, \mathbb{C})\) corresponds to the irreducible representation associated with the minuscule fundamental weight \(\lambda\) of \(\mathsf{G}^\vee\). All fundamental weights of \(GL(n, \mathbb{C})\) are minuscule, and the corresponding irreducible representations are given by the exterior power representations \(\bigwedge^k V\). This justifies the identification \(H^\bullet(G(k,n), \mathbb{C}) \cong \bigwedge^k H^\bullet(\mathbb{P}^{n-1}, \mathbb{C})\) as vector spaces and endows \(H^\bullet(G(k,n), \mathbb{C})\) with a natural module structure over \(\mathfrak{gl}(H^\bullet(\mathbb{P}^{n-1}, \mathbb{C})) \cong \mathfrak{gl}_n(\C)\).

Under this identification, the operators \(X_1, X_2 \in \mathfrak{gl}(H^\bullet(\mathbb{P}^{n-1}, \mathbb{C}))\) corresponding to classical cup multiplication by \(c_1(E_1^{\mathbb{P}})\) and \(c_1(E_2^{\mathbb{P}})\) act on \(H^\bullet(G(k,n), \mathbb{C})\) via classical cup multiplication by \(c_1(E_1)\) and \(c_1(E_2)\), respectively.

The framework described here represents the \textit{classical} geometric Satake correspondence as established in  \cite{Gin95,Gin08,MV00,MV07}. Remarkably, there have been several extensions in both equivariant and quantum directions, thanks to the efforts of many researchers from various perspectives \cite{Gat05,GS09,GS10,Lak08,Lak09,LT09,GM11}. These results originally inspired the present work. For further details, we refer the reader to \cite{CDG1} and references therein.
\vskip2mm
\noindent 1.7.\,\,{\bf Plan of the paper.} In {Section \ref{SEC1}}, we introduce the joint system of dynamical differential and qKZ difference equations on the weight spaces \(\left(\mathbb{C}^2\right)^{\otimes n}_{k,n-k}\) and recall some of its fundamental properties.

In {Section \ref{sec2}}, we define a linear isomorphism, called the Satake identification \(\theta_{k,n}\), between \(\bigwedge^k\left(\mathbb{C}^2\right)^{\otimes n}_{1,n-1}\) and \(\left(\mathbb{C}^2\right)^{\otimes n}_{k,n-k}\). We prove that, under this identification, the joint system of equations is gauge equivalent to the induced system on \(\bigwedge^k\left(\mathbb{C}^2\right)^{\otimes n}_{1,n-1}\).

In {Section \ref{SEC3}}, we introduce the equivariant cohomology of Grassmannians.

In {Section \ref{SEC4}}, we discuss the derived category of equivariant coherent sheaves on Grassmannians and their equivariant \(K\)-theory. We introduce the \textcyr{B}-morphism, a twisted version of the equivariant Chern character, and the concept of exceptional collections and bases.

In {Section \ref{sec5}}, we introduce the equivariant quantum cohomology of Grassmannians and define stable envelope maps, which are used to construct the joint system of qDE and qKZ equations.

In {Section \ref{SEC6}}, we describe several constructions of solutions to the joint system of qDE and qKZ equations, including the topological-enumerative morphism, multidimensional hypergeometric solutions and their relation to \(K\)-theoretic classes. We prove the first main theorem and show how to construct solutions for \(G(k,n)\) from solutions for \(\mathbb{P}^{n-1}\).

In {Section \ref{SEC7}}, we prove the second main theorem, concerning the compatibility of the \textcyr{B}-theorem and the Satake identification, leading to a determinantal formula and new integral representations of the hypergeometric solutions.

In {Section \ref{SEC8}}, we analyze the Stokes phenomenon for the joint system, showing the correspondence between Stokes bases and suitable \(T\)-full exceptional collections, and demonstrating that the Stokes matrices are Gram matrices of the Euler–Poincaré–Grothendieck pairing in \(K\)-theory.

\vskip2mm
\noindent{\bf Acknowledgements.} This research was supported by the FCT Projects UIDB/00208/2020 (DOI: 10.54499/UIDB/00208/2020), UIDP/00208/2020 (DOI: 10.54499/UIDP/00208/2020), 2021.01521.\\CEECIND (DOI:10.13039/501100005765), and the project GENIDE 2022.03702.PTDC (DOI: 10.54\-499/2022.03702.PTDC). The first author is a member of the COST Action CA21109 CaLISTA.

\section{Dynamical and qKZ equations}\label{SEC1}
\subsection{Notations}\label{sec1.1}Let $0\leq k\leq n$. Let $\mc I_{k,n}$ be the set of all double partitions $I=(I_1,I_2)$ of $\{1,\dots,n\}$ with 
\[I_1\coprod I_2=\{1,\dots, n\}, \quad {\rm card\,} I_1=k,\quad {\rm card\,} I_2=n-k.
\]Consider the space $\C^2$ with the standard basis $v_1=(1,0)$ and $v_2=(0,1)$.  The tensor product 
$(\C^2)^{\otimes n}$ has 
basis $(v_I)_I$, indexed by the set $\bigcup_{0\leq k\leq n}\mc I_{k,n}$, where
\[v_I=v_{j_1}\otimes\dots\otimes v_{j_n},\qquad j_i=1\text{ if }i\in I_1,\qquad j_i=2\text{ if }i\in I_2.
\]The basis $(v_I)_I$ will always be lexicographically ordered. We have a weight decomposition
\[(\C^2)^{\otimes n}=\sum_{k=0}^n(\C^2)^{\otimes n}_{k,n-k},\qquad (\C^2)^{\otimes n}_{k,n-k}\text{ is the subspace with basis }(v_I)_{I\in\mc I_{k,n}}.
\]The space $(\C^2)^{\otimes n}$ is a $\gl_2$-module. The Lie algebra $\gl_2$ has basis $e_{ij}=(\dl_{ia}\dl_{jb})_{a,b=1,2}$, for $i,j=1,2$.

\subsection{Difference qKZ equations}Introduce the $R$-matrix $R(u)$, with $u\in\C$, acting on $(\C^2)^{\otimes 2}$:
\[R(u)=P+u\,e_{11}\otimes e_{22},\quad P\text{ permutation of factors of $(\C^2)^{\otimes 2}$}.
\]In the basis $(v_i\otimes v_j)_{i,j}$, the matrix $R$ equals
\[R(u)=\begin{pmatrix}
1&0&0&0\\
0&u&1&0\\
0&1&0&0\\
0&0&0&1
\end{pmatrix}.
\]
Let $\bm p=(p_1,p_2)\in\C^2$. Define the {\it qKZ operators} $K_1,\dots, K_n$ acting on $(\C^2)^{\otimes n}$:
\begin{multline}
\label{qkzop}K_a(\bm z;\bm p;\kappa)=R(z_a-z_{a-1}+\kappa)^{(a,a-1)}\dots R(z_a-z_{1}+\kappa)^{(a,1)}\times\\
\times p_1^{e_{11}^{(a)}}p_2^{e_{22}^{(a)}} R(z_i-z_n)^{(a,n)}\dots R(z_a-z_{a+1})^{(a,a+1)},
\end{multline}for $a=1,\dots,n$. The superscripts indicate the factors of $(\C^2)^{\otm n}$ on which the corresponding operator acts.
The qKZ operators $K_1,\dots, K_n$ preserve the weight decomposition of $(\C^2)^{\otimes n}$.

\begin{example}\label{exqkz}
Let $n=3$. The qKZ operators are
\begin{align*}
K_1&=p_1^{e_{11}^{(1)}}p_2^{e_{22}^{(1)}}R(z_1-z_3)^{(1,3)}R(z_1-z_2)^{(1,2)},\\
K_2&=R(z_2-z_1+\kappa)^{(2,1)}p_1^{e_{11}^{(2)}}p_2^{e_{22}^{(2)}}R(z_2-z_3)^{(2,3)},\\
K_3&=R(z_3-z_2+\kappa)^{(3,2)}R(z_3-z_2+\kappa)^{(3,1)}p_1^{e_{11}^{(3)}}p_2^{e_{2}^{(3)}}.
\end{align*}
We have the weight decomposition of $(\C^2)^{\otimes 3}$:
\[(\C^2)^{\otimes 3}_{3,0}={\rm span}(v_1\otimes v_1\otimes v_1),
\]
\[(\C^2)^{\otimes 3}_{2,1}={\rm span}(v_1\otimes v_1\otimes v_2,\quad v_1\otimes v_2\otimes v_1,\quad v_2\otimes v_1\otimes v_1),
\]
\[
(\C^2)^{\otimes 3}_{1,2}={\rm span}(v_1\otimes v_2\otimes v_2,\quad v_2\otimes v_1\otimes v_2,\quad v_2\otimes v_2\otimes v_1),
\]
\[(\C^2)^{\otimes 3}_{0,3}={\rm span}(v_2\otimes v_2\otimes v_2).
\]
\vskip2mm
With respect to these bases, the matrices of the qKZ operators $K_1,K:2,K_3$ are:
\begin{itemize}
\item on $(\C^2)^{\otm 3}_{3,0}$:
\[K_1=K_2=K_3=p_1;
\]
\item on $(\C^2)^{\otm 3}_{2,1}$:
\beq\label{n3K1.1}
K_1=\left(
\begin{array}{ccc}
 p_1 (z_1-z_3) & p_1 & 0 \\
 0 & p_1 (z_1-z_2) & p_1 \\
 p_2 & 0 & 0 \\
\end{array}
\right),
\eeq
\beq\label{n3K2.1}
K_2=\left(
\begin{array}{ccc}
 p_1 (z_2-z_3) & p_1 & 0 \\
 0 & 0 & p_1 \\
 p_2 & 0 & p_1 (\kappa-z_1+z_2) \\
\end{array}
\right),
\eeq
\beq\label{n3K3.1}
K_3=\left(
\begin{array}{ccc}
 0 & p_1 & 0 \\
 0 & p_1 (\kappa-z_2+z_3) & p_1 \\
 p_2 & 0 & p_1 (\kappa-z_1+z_3) \\
\end{array}
\right);
\eeq
\item on $(\C^2)^{\otm 3}_{1,2}$:
\beq\label{n3K1.2}
K_1=
\left(
\begin{array}{ccc}
 p_1 (z_1-z_2) (z_1-z_3) & p_1 (z_1-z_3) & p_1 \\
 p_2 & 0 & 0 \\
 p_2 (z_1-z_2) & p_2 & 0 \\
\end{array}
\right),
\eeq

\beq\label{n3K2.2}
K_2=
\left(
\begin{array}{ccc}
 0 & p_1 (z_2-z_3) & p_1 \\
 p_2 & p_1 (\kappa-z_1+z_2) (z_2-z_3) & p_1 (\kappa-z_1+z_2) \\
 0 & p_2 & 0 \\
\end{array}
\right),
\eeq

\beq\label{n3K3.2}
K_3=
\left(
\begin{array}{ccc}
 0 & 0 & p_1 \\
 p_2 & 0 & p_1 (\kappa-z_1+z_3) \\
 p_2 (\kappa-z_2+z_3) & p_2 & p_1 (\kappa-z_1+z_3) (\kappa-z_2+z_3) \\
\end{array}
\right);
\eeq
\item on $(\C^2)^{\otm 3}_{0,3}$:
\[K_1=K_2=K_3=p_2.
\]
\end{itemize}
\qetr
\end{example}

\begin{thm}\cite{TV23}
The qKZ operators $K_1,\dots, K_n$ define a discrete flat connection with parameter $\kappa$:
\beq\label{discrconn}
K_a(z_1,\dots,z_b+\kappa,\dots,z_n;\bm p;\kappa)K_b(\bm z;\bm p;\kappa)=K_b(z_1,\dots,z_a+\kappa,\dots,z_n;\bm p;\kappa)K_a(\bm z;\bm p;\kappa),
\eeq
for all $a,b=1,\dots,n$.\qed
\end{thm}

The system of difference equations with step $\kappa$,
\beq\label{qKZ.0}
f(z_1,\dots,z_a+\kappa,\dots,z_n;\bm p)=K_a(\bm z;\bm p;\kappa)f(\bm z;\bm p),\quad a=1,\dots,n,
\eeq
for a $(\C^2)^{\otm n}$-valued function $f(\bm z;\bm p)$ is called system of {\it qKZ equations}.
\subsection{Differential dynamical equations}
For all $i,j=1,2$, $a,b=1,\dots,n$ define the operators $Q^{a,b}_{i,j}$ acting on $(\C^2)^{\otm n}$ by the rule
\begin{align*}
&Q^{a,b}_{i,j}v_I=v_{s_{ab}(I)},\quad &&\text{if $a\in I_i$, $b\in I_j$, and $(j-i)(a-b)\equiv1\,({\rm mod}\,n)$,}\\
&Q^{a,b}_{i,j}v_I=0,\quad&&\text{otherwise.}
\end{align*}
The {\it dynamical operators} $X_1,X_2$ acting on $(\C^2)^{\otm n}$ are given by the formula
\begin{align}
\label{dynop1}
X_1(\bm z;\bm p)&=\sum_{a=1}^nz_a e_{11}^{(a)}+\sum_{1\leq b<a\leq n}Q^{a,b}_{1,2}+\frac{p_2}{p_1}\sum_{1\leq a<b\leq n}Q^{a,b}_{1,2},\\
\label{dynop2}
X_2(\bm z;\bm p)&=\sum_{a=1}^nz_a e_{22}^{(a)}-\frac{p_2}{p_1}\sum_{1\leq b<a\leq n}Q^{a,b}_{2,1}-\sum_{1\leq a<b\leq n}Q^{a,b}_{2,1}.
\end{align}
The operators $X_1,X_2$ preserve the weight decomposition of $(\C^2)^{\otm n}$.
\begin{rem}
The only nonzero operators $Q^{a,b}_{1,2}$ appearing in the operator $X_1$ are $Q^{1,n}_{1,2},\,Q^{a,a-1}_{1,2}$ for $a=2,\dots,n$. That is, we have the simplified expression
\[X_1(\bm z;\bm p)=\sum_{a=1}^nz_a e_{11}^{(a)}+\frac{p_2}{p_1}Q^{1,n}_{1,2}+\sum_{a=2}^nQ^{a,a-1}_{1,2}.
\]

Similarly, the only nonzero operators $Q^{a,b}_{2,1}$ appearing in the operator $X_2$ are $Q^{n,1}_{2,1},\,Q^{a,a+1}_{2,1}$ for $a=1,\dots,n-1$. That is, we have the simplified expression
\[X_2(\bm z;\bm p)=\sum_{a=1}^nz_a e_{22}^{(a)}-\frac{p_2}{p_1}Q^{n,1}_{2,1}-\sum_{a=1}^{n-1}Q^{a,a+1}_{2,1}.
\]
\end{rem}

The system of differential equations with parameter $\kappa$,
\beq\label{dyn.0}
\kappa p_i\frac{\der f}{\der p_i}=X_i(\bm z;\bm p)f,\quad i=1,2,
\eeq
for a $(\C^2)^{\otm n}$-valued function $f(\bm z;\bm p)$ will be called the {\it dynamical equation}.

\begin{thm}\cite{TV23}\label{cjs}
The joint system of dynamical and qKZ equations with the same parameter $\kappa$ is compatible.This means that besides equations \eqref{discrconn}, we have the following compatibility conditions:
\begin{align}
\label{ccond1}
&\kappa p_i\frac{\der X_j}{\der p_i}-\kappa p_j\frac{\der X_i}{\der p_j}=[X_i,X_j],\\
\label{ccond2}
&\kappa p_i\frac{\der K_a(\bm z;\bm p;\kappa)}{\der p_i}=X_i(z_1,\dots,z_a+\kappa,\dots,z_n;\bm p)K_a(\bm z;\bm p;\kappa)-K_a(\bm z;\bm p;\kappa)X_i(\bm z;\bm p),
\end{align}
for any $i,j=1,2$ and any $a=1,\dots,n$.\qed
\end{thm}

\begin{example}\label{ex:ex2}
Let $n=3$. The operators 
\[X_1=\sum_{a=1}^nz_a e_{11}^{(a)}+Q^{2,1}_{1,2}+Q^{3,2}_{1,2}+\frac{p_2}{p_1}Q^{1,3}_{1,2},\qquad X_2=\sum_{a=1}^nz_a e_{22}^{(a)}-Q^{1,2}_{2,1}-Q^{2,3}_{2,1}-\frac{p_2}{p_1}Q^{3,1}_{2,1},
\]preserve the subspaces 
\[(\C^2)^{\otm 3}_{3,0}\text{ spanned by }v_1\otm v_1\otm v_1,
\]
\[(\C^2)^{\otm 2}_{2,1}\text{ spanned by }v_1\otm v_1\otm v_2,\quad v_1\otm v_2\otm v_1,\quad v_2\otm v_1\otm v_1
\]
\[(\C^2)^{\otm 2}_{1,2}\text{ spanned by }v_1\otm v_2\otm v_2,\quad v_2\otm v_1\otm v_2,\quad v_2\otm v_2\otm v_1,
\]
\[(\C^2)^{\otm 3}_{0,3}\text{ spanned by }v_2\otm v_2\otm v_2.
\]

With respect to these bases, the matrices of the dynamical operators $X_1,X_2$ are:
\begin{itemize}
\item on $(\C^2)^{\otm 3}_{3,0}$
\[X_1=z_1+z_2+z_3,\quad X_2=0;
\]
\item on $(\C^2)^{\otm 3}_{2,1}$
\[X_1=\begin{pmatrix}
z_1+z_2&1&0\\
0&z_1+z_3&1\\
\frac{p_2}{p_1}&0&z_2+z_3
\end{pmatrix},\quad 
X_2=\begin{pmatrix}
z_3&-1&0\\
0&z_2&-1\\
-\frac{p_2}{p_1}&0&z_1
\end{pmatrix};
\]
\item on $(\C^2)^{\otm 3}_{1,2}$
\[X_1=\begin{pmatrix}
z_1&1&0\\
0&z_2&1\\
\frac{p_2}{p_1}&0&z_3
\end{pmatrix},\quad 
X_2=\begin{pmatrix}
z_2+z_3&-1&0\\
0&z_1+z_3&-1\\
-\frac{p_2}{p_1}&0&z_1+z_2
\end{pmatrix};
\]
\item on $(\C^2)^{\otm 3}_{0,3}$
\[X_1=0,\quad X_2=z_1+z_2+z_3.
\]
\end{itemize}
The reader can check that the dynamical operators $X_1,X_2$, together with the qKZ operators of Example \ref{exqkz}, satisfy the compatibility conditions \eqref{ccond1} and \eqref{ccond2}.
\qetr
\end{example}

\section{Dynamical and qKZ equations via Satake identifications}\label{sec2}
\subsection{Satake identifications}Consider the first level $(\C^2)^{\otm n}_{1,n-1}$. Its basis elements are labelled by $\mc I_{1,n}$, that is partitions $I=(I_1,I_2)$ of the form 
\[[a]:=\left(\{a\},\,\,\,\{1,\dots,n\}\setminus\{a\}\right),\qquad \text{for $a=1,\dots, n$. }
\]
The basis $(v_{[1]},\dots, v_{[n]})$ induce a basis of the $k$-th exterior power $\bigwedge^k_\C\left(\C^2\right)^{\otm n}_{1,n-1}$, of the form
\[v_{[a_1]}\wedge\dots\wedge v_{[a_k]},\quad 1\leq a_1<a_2<\dots<a_k\leq n.
\]
\begin{prop}
The $\C$-linear map $\theta_{k,n}\colon \bigwedge^k_\C\left(\C^2\right)^{\otm n}_{1,n-1}\to (\C^2)^{\otm n}_{k,n-k}$, defined on the basis elements by
\beq\label{thk}
v_{[a_1]}\wedge\dots\wedge v_{[a_k]}\mapsto v_I,\quad I=(\{a_1,\dots,a_k\},\,\,\,\{1,\dots,n\}\setminus\{a_1,\dots,a_k\}),
\eeq
is an isomorphism of $\C$-vector spaces.\qed
\end{prop}
The isomorphism $\theta_{k,n}$ will be called the {\it $k$-th level Satake identification}.
\subsection{Exterior power of dynamical and qKZ difference equations}
Given a linear endomorphism $f$ of $\left(\C^2\right)^{\otm n}_{1,n-1}$, 
denote by $\bigwedge^kf$ the induced morphism on $\bigwedge^k_\C\left(\C^2\right)^{\otm n}_{1,n-1}$ defined by
\[\bigwedge\nolimits^kf\left(v_{[a_1]}\wedge\dots\wedge v_{[a_k]}\right)=fv_{[a_1]}\wedge\dots\wedge fv_{[a_k]},
\]for any $1\leq a_1<\dots<a_k\leq n$.

For any $a=1,\dots, n$, define the operator $\widetilde K_a^{(k)}(\bm z;\bm p;\kappa)$ acting on 
$\bigwedge^k_\C\left(\C^2\right)^{\otm n}_{1,n-1}$ 
by the formula
\beq\label{tildeka}\widetilde K_a^{(k)}(\bm z;\bm p;\kappa)=
\bigwedge\nolimits^kK_a[\bm z;(-1)^{k-1}p_1,p_2;\kappa].
\eeq Moreover, define two operators $\widetilde X_1^{(k)}(\bm z;\bm p)$ and $\widetilde X_2^{(k)}(\bm z;\bm p)$ on $\bigwedge^k_\C(\C^2)^{\otm n}_{1,n-1}$ by
\beq\label{tildexi}\widetilde X_i^{(k)}(\bm z;\bm p)(v_{[a_1]}\wedge\dots\wedge v_{[a_k]})=\sum_{j=1}^kv_{[a_1]}\wedge\dots\wedge X_i(\bm z;(-1)^{k-1}p_1,p_2)v_{[a_j]}\wedge\dots\wedge v_{[a_k]},\quad i=1,2,
\eeq for any increasing integer sequence $1\leq a_1<\dots<a_k\leq n$. Via the Satake identification $\theta_{k,n}$, the operators $\widetilde K_a^{(k)},\widetilde X_i^{(k)}$ induce operators $\widehat K_a^{(k)},\widehat X_i^{(k)}$ acting on $(\C^2)^{\otm n}_{k,n-k}$:
\beq
\widehat K_a^{(k)}:=\theta_{k,n}\circ\widetilde K_a^{(k)}\circ \theta_{k,n}^{-1},\qquad\widehat X_i^{(k)}=\theta_{k,n}\circ\widetilde X_i^{(k)}\circ \theta_{k,n}^{-1},
\eeq for all $a=1,\dots, n$ and $i=1,2$.
\vskip2mm
Consider then the joint system of difference--differential equations
\begin{align}
\label{qKZ.1}
f(z_1,\dots, z_a+\kappa,\dots,z_n;\bm p;\kappa)&=\widehat K_a^{(k)}(\bm z;\bm p;\kappa)f(\bm z;\bm p),&a=1,\dots,n,\\
\label{dyn.1}
\kappa p_i\frac{\der f}{\der p_i}&=\widehat X_i^{(k)}(\bm z;\bm p) f,&i=1,2,
\end{align}where $f(\bm z;\bm p)$ is $(\C^2)^{\otm n}_{k,n-k}$-valued function. By construction, this is a compatible system of equations. The next result clarifies its relation with the joint system of qKZ and dynamical equations \eqref{qKZ.0} and \eqref{dyn.0} on $(\C^2)^{\otm n}_{k,n-k}$.

\begin{thm}\label{abstractthmsat}
The joint system of equations \eqref{qKZ.1},\eqref{dyn.1} is gauge equivalent with the joint system of qKZ and dynamical equations on $(\C^2)^{\otm n}_{k,n-k}$. More precisely, if $f(\bm z;\bm p)$ is a solution of \eqref{qKZ.1},\eqref{dyn.1}, then the function
\[G(\bm z;\bm p;\kappa)f(\bm z;\bm p),\quad \text{with }G(\bm z;\bm p;\kappa)=p_2^{(1-k)\sum_{a=1}^n z_a/\kappa},
\]is a solution of the qKZ and dynamical equations \eqref{qKZ.0} and \eqref{dyn.0}.
\end{thm}
The claim is equivalent to the following identities among the operators $K_a,\widehat K_a^{(k)}, X_i,\widehat X_i^{(k)}$ acting on $(\C^2)^{\otm n}_{k,n-k}$:
\begin{align}
\nonumber
K_a(\bm z;\bm p;\kappa)&=G(z_1,\dots,z_a+\kappa,\dots z_n;\bm p;\kappa)\widehat K_a^{(k)}(\bm z;\bm p;\kappa) G(\bm z;\bm p;\kappa)^{-1}&\\
\label{gauge1}
&=p_2^{1-k}\widehat K_a^{(k)}(\bm z;\bm p;\kappa),&a=1,\dots,n,\\
\nonumber
X_1(\bm z;\bm p)&=G(\bm z;\bm p;\kappa)\widehat X_1^{(k)}(\bm z;\bm p)G(\bm z;\bm p;\kappa)^{-1}+\kappa p_1\frac{\der G(\bm z;\bm p;\kappa)}{\der p_1}G(\bm z;\bm p;\kappa)^{-1}&\\
\label{gauge2}
&=\widehat X_1^{(k)}(\bm z;\bm p),&\\
\nonumber
X_2(\bm z;\bm p)&=G(\bm z;\bm p;\kappa)\widehat X_2^{(k)}(\bm z;\bm p)G(\bm z;\bm p;\kappa)^{-1}+\kappa p_2\frac{\der G(\bm z;\bm p;\kappa)}{\der p_2}G(\bm z;\bm p;\kappa)^{-1}&\\
\label{gauge3}
&=\widehat X_2^{(k)}(\bm z;\bm p)+(1-k)\sum_{a=1}^nz_a.&
\end{align}
These identities will be proved in the next two Sections \ref{secqkzsat} and \ref{secqkzdyn}. 
\vskip2mm
The following result directly follows from Theorem \ref{abstractthmsat}.
\begin{cor}\label{cor1}
Given $k$ solutions $f_1(\bm z;\bm p),\dots, f_k(\bm z;\bm p)$ of the joint system of dynamical and qKZ equations on $(\C^2)^{\otm n}_{1,n-1}$, the function
\beq\label{extprodsol}
G(\bm z;\bm p;\kappa)\cdot \theta_{k,n}\left(f_1(\bm z;(-1)^{k-1}p_1,p_2)\wedge\dots\wedge f_k(\bm z;(-1)^{k-1}p_1,p_2)\right),\quad G(\bm z;\bm p;\kappa)=p_2^{(1-k)\sum_{a=1}^n z_a/\kappa},
\eeq
is a solution of the joint system of dynamical and qKZ system of equations on $(\C^2)^{\otm n}_{k,n-k}$.\qed
\end{cor}

\begin{rem}
The solutions of the joint systems of dynamical and qKZ equations are, in general, multivalued. Consequently, in \eqref{extprodsol} a determination of the argument of $-1$ has to be fixed.
\end{rem}
\subsection{qKZ operators and Satake identifications}\label{secqkzsat}

In this Section we prove the first set of identities \eqref{gauge1}.
\begin{thm}\label{satqkz}
The isomorphism $\theta_{k,n}$ intertwines the operators $\widetilde K^{(k)}_a(\bm z;\bm p;\kappa)$ 
acting on $\bigwedge^k_\C\left(\C^2\right)^{\otm n}_{1,n-1}$ and the operator 
$p_2^{k-1}K_a(\bm z;\bm p;\kappa)$ acting on $(\C^2)^{\otm n}_{k,n-k}$, that is
\beq
\label{idqkzsat}
\theta_{k,n}\circ \widetilde K^{(k)}_a=p_2^{k-1}K_a\circ \theta_{k,n},
\eeq
for any $a=1,\dots, n$.
\end{thm}

For the proof, we need some preliminary Lemmata.

\begin{lem}\label{lemsatqkz}
Let $u\in\C$, $i,j\in\{1,\dots,n\}$, with $i\neq j$, and $1\leq a_1<\dots<a_k\leq n$.
\begin{enumerate}
\item If $i,j\neq a_1,\dots, a_k$, we have
\[\bigwedge\nolimits^kR(u)^{(i,j)}[v_{[a_1]}\wedge\dots\wedge v_{[a_k]}]=v_{[a_1]}\wedge\dots\wedge v_{[a_k]}.
\]

\item If $j=a_h$, and $i\neq a_1,\dots,a_k$,  
we have
\[\bigwedge\nolimits^kR(u)^{(i,j)}[v_{[a_1]}\wedge\dots\wedge v_{[a_k]}]=v_{[a_1]}\wedge\dots\wedge v_{[a_{h-1}]}\wedge v_{[i]}\wedge v_{[a_{h+1}]}\wedge\dots\wedge v_{[a_k]}.
\]
\item If $i=a_h$, and $j\neq a_1,\dots,a_k$, we have
\begin{multline*}
\bigwedge\nolimits^kR(u)^{(i,j)}[v_{[a_1]}\wedge\dots\wedge v_{[a_k]}]
=v_{[a_1]}\wedge\dots\wedge v_{[a_{h-1}]}\wedge\left(-u v_{[a_h]}+v_{[j]}\right)\wedge v_{[a_{h+1}]}\wedge\dots\wedge v_{[a_k]}.
\end{multline*}
\item If $i=a_h$ and $j=a_\ell$, with $h\neq \ell$, we have
\begin{align*}\bigwedge\nolimits^kR(u)^{(i,j)}[v_{[a_1]}\wedge\dots\wedge v_{[a_k]}]&=v_{[a_1]}\wedge\dots\wedge\underset{h\text{-th}}{v_{[a_\ell]}}\wedge\dots\wedge\underset{\ell\text{-th}}{v_{[a_h]}}\wedge\dots\wedge v_{[a_k]}\\
&=-v_{[a_1]}\wedge\dots\wedge v_{[a_k]}.
\end{align*}
\end{enumerate}
\end{lem}
\proof
The operator $R(u)^{(i,j)}$, with $i\neq j$, acts on $(\C^2)^{\otm n}_{1,n-1}$ as follows:
\begin{align*}
&R(u)^{(i,j)}v_{[a]}=v_{[a]}\quad&\text{if $i,j\neq a$,}\\
&R(u)^{(i,j)}v_{[a]}=v_{[i]}\quad&\text{if $i\neq a,j=a$,}\\
&R(u)^{(i,j)}v_{[a]}=-uv_{[a]}+v_{[j]}\quad&\text{if $i=a,j\neq a$}.
\end{align*}
The Lemma immediately follows.
\endproof

\begin{lem}\label{lem1satqkz}
Let $u_1,\dots,u_{n-1}\in\C$, $i\in\{1,\dots, n\}$, and $1\leq a_1<\dots<a_k\leq n$. Let $h\in\{1,\dots, k\}$ be such that $a_{h-1}<i\leq a_h$. We have
\begin{multline}\label{id1}
\left(\theta_{k,n}\circ\bigwedge\nolimits^k\left[((-1)^{k-h})^{e_{11}^{(i)}}R(u_{n-i})^{(i,n)}R(u_{n-i-1})^{(i,n-1)}\dots R(u_1)^{(i,i+1)}\right]\right)v_{[a_1]}\wedge\dots\wedge v_{[a_k]}\\
=R(u_{n-i})^{(i,n)}R(u_{n-i-1})^{(i,n-1)}\dots R(u_1)^{(i,i+1)}\left[\theta_{k,n}(v_{[a_1]}\wedge\dots\wedge v_{[a_k]})\right],
\end{multline}
\begin{multline}\label{id2}
\left(\theta_{k,n}\circ\bigwedge\nolimits^k\left[R(u_{n-1})^{(i,i-1)}\dots R(u_{n-i+1})^{(i,1)}((-1)^{h-1})^{e_{11}^{(i)}}\right]\right)v_{[a_1]}\wedge\dots\wedge v_{[a_k]}\\
=R(u_{n-1})^{(i,i-1)}\dots R(u_{n-i+1})^{(i,1)}\left[\theta_{k,n}(v_{[a_1]}\wedge\dots\wedge v_{[a_k]})\right].
\end{multline}
\end{lem}

\proof 
 Introduce the sequence $0\leq d_0<d_1<d_2<\dots<d_{k-h}\leq n-i$ such that
\[i+d_0=a_h, \quad i+d_1=a_{h+1},\quad \dots,\quad i+d_{k-h}=a_k.
\]If $d_0>0$, the product of $R$-matrices $\bigwedge^k\left[R(u_{d_0})^{(i,i+d_0)}\dots R(u_2)^{(i,i+2)}R(u_1)^{(i,i+1)}\right]$ acts on $v_{[a_1]}\wedge\dots \wedge v_{[a_k]}$ as follows
\[v_{[a_1]}\wedge\dots \wedge v_{[a_k]}\mapsto v_{[a_1]}\wedge\dots \wedge v_{[a_{h-1}]}\wedge v_{[i]}\wedge v_{[a_{h+1}]}\wedge\dots \wedge v_{[a_k]}.
\]More precisely, the only operator which acts non-identically is the last one $\bigwedge^kR(u_{d_0})^{(i,i+d_0)}$. On the other hand, if $d_0=0$ we immediately start from the vector $$v_{[a_1]}\wedge\dots \wedge v_{[a_k]}=v_{[a_1]}\wedge\dots \wedge v_{[a_{h-1}]}\wedge v_{[i]}\wedge v_{[a_{h+1}]}\wedge\dots \wedge v_{[a_k]}.$$
Let us then consider the action of the product of $R$-matrices $$\bigwedge\nolimits^k\left[R(u_{d_1-1})^{(i,i+d_1-1)}\dots R(u_{d_0+2})^{(i,i+d_0+2)}R(u_{d_0+1})^{(i,i+d_0+1)}\right]$$ on $v_{[a_1]}\wedge\dots \wedge v_{[a_{h-1}]}\wedge v_{[i]}\wedge v_{[a_{h+1}]}\wedge\dots \wedge v_{[a_k]}$.
By a consecutive application of point (3) of Lemma \ref{lemsatqkz}, one sees that the product of $R$-matrices $\bigwedge\nolimits^k\left[R(u_{d_1-1})^{(i,i+d_1-1)}\dots R(u_{1})^{(i,i+1)}\right]$ acts on the initial vector $v_{[a_1]}\wedge\dots\wedge v_{[a_k]}$ as follows:
\begin{multline*}v_{[a_1]}\wedge\dots\wedge v_{[a_k]}\mapsto
\left(\prod_{r=d_0+1}^{d_1-1}(-u_r)\right)v_{[a_1]}\wedge\dots\wedge\underset{h\text{-th}}{v_{[i]}}\wedge\dots\wedge v_{[a_k]}\\+\sum_{j=d_0+1}^{d_1-1}
\left(\prod_{r=d_0+1}^{j-1}(-u_r)\right)v_{[a_1]}\wedge\dots\wedge\underset{h\text{-th}}{v_{[i+j]}}\wedge\dots\wedge v_{[a_k]}.
\end{multline*} 
If we further apply the operator $\bigwedge^kR(u_{d_1})^{(i,i+d_1)}$, we obtain
\begin{multline*}
\left(\prod_{r=d_0+1}^{d_1-1}(-u_r)\right)v_{[a_1]}\wedge\dots\wedge\underset{h\text{-th}}{v_{[i+d_1]}}\wedge\underset{(h+1)\text{-th}}{v_{[i]}}\wedge\dots\wedge v_{[a_k]}\\
+\sum_{j=d_0+1}^{d_1-1}
\left(\prod_{r=d_0+1}^{j-1}(-u_r)\right)v_{[a_1]}\wedge\dots\wedge\underset{h\text{-th}}{v_{[i+j]}}\wedge\underset{(h+1)\text{-th}}{v_{[i]}}\wedge\dots\wedge v_{[a_k]}\\
=\sum_{j=d_0+1}^{d_1}
\left(\prod_{r=d_0+1}^{j-1}(-u_r)\right)v_{[a_1]}\wedge\dots\wedge\underset{h\text{-th}}{v_{[i+j]}}\wedge\underset{(h+1)\text{-th}}{v_{[i]}}\wedge\dots\wedge v_{[a_k]},
\end{multline*}
where we applied points (3) and (4) of Lemma \ref{lemsatqkz}. By proceeding with the same argument, one sees that the product of $R$-matrices $\bigwedge\nolimits^k\left[R(u_{d_{k-h}})^{(i,i+d_{k-h})}\dots R(u_{1})^{(i,i+1)}\right]$ acts on the initial vector $v_{[a_1]}\wedge\dots\wedge v_{[a_k]}$ as follows:
\begin{multline*}v_{[a_1]}\wedge\dots\wedge v_{[a_k]}\mapsto\\
\sum_{j_1=d_0+1}^{d_1}\dots\sum_{j_{k-h}=d_{k-h-1}+1}^{d_{k-h}}
\left(\prod_{b=1}^{k-h}\prod_{r_b=d_{b-1}+1}^{j_b-1}(-u_{r_b})\right)
v_{[a_1]}\wedge\dots\wedge\underset{h\text{-th}}{v_{[i+j_1]}}\wedge v_{[i+j_2]}\wedge\dots\wedge v_{[i+j_{k-h}]}\wedge v_{[i]}.
\end{multline*} 
If we apply the last product $\bigwedge^k[R(u_n)^{(i,n)}\dots R(u_{d_{k-h}+1})^{(i,i+d_{k-h}+1)}]$, we obtain the vector
\begin{multline*}
\sum_{j_1=d_0+1}^{d_1}\dots\sum_{j_{k-h}=d_{k-h-1}+1}^{d_{k-h}}
\left(\prod_{b=1}^{k-h}\prod_{r_b=d_{b-1}+1}^{j_b-1}(-u_{r_b})\right)\left(\prod_{p=d_{k-h}+1}^{n}(-u_p)\right)\\
\cdot v_{[a_1]}\wedge\dots\wedge\underset{h\text{-th}}{v_{[i+j_1]}}\wedge v_{[i+j_2]}\wedge\dots\wedge v_{[i+j_{k-h}]}\wedge v_{[i]}\\
+\sum_{j_1=d_0+1}^{d_1}\dots\sum_{{j_{k-h}=d_{k-h-1}+1}}^{d_{k-h}}\,\,\sum_{{j_{k-h+1}=d_{k-h}+1}}^{n-i}\left(\prod_{b=1}^{k-h+1}\prod_{r_b=d_{b-1}+1}^{j_b-1}(-u_{r_b})\right)\\
\cdot v_{[a_1]}\wedge\dots\wedge\underset{h\text{-th}}{v_{[i+j_1]}}\wedge v_{[i+j_2]}\wedge\dots\wedge v_{[i+j_{k-h}]}\wedge v_{[i+j_{k-h+1}]}.
\end{multline*}
If we apply the operator $\bigwedge^k[((-1)^{k-h})^{e_{11}^{(i)}}]$, each term $v_{[a_1]}\wedge\dots\wedge\underset{h\text{-th}}{v_{[i+j_1]}}\wedge v_{[i+j_2]}\wedge\dots\wedge v_{[i+j_{k-h}]}\wedge v_{[i]}$ is multiplied by $(-1)^{k-h}$, and it can equivalently be rewritten in lexicographical order
\[v_{[a_1]}\wedge\dots\wedge\underset{h\text{-th}}{v_{[i]}}\wedge{v_{[i+j_1]}}\wedge v_{[i+j_2]}\wedge\dots\wedge v_{[i+j_{k-h}]}.
\]If we finally apply $\theta_{k,n}$, and write for short $v_{\{a_1,\dots,a_k\}}$ rather than $v_{\{a_1,\dots,a_k\},\{1,\dots n\}\setminus\{a_1,\dots,a_k\}}$, the l.h.s.\,\,of \eqref{id1} equals
\begin{multline*}
\sum_{j_1=d_0+1}^{d_1}\dots\sum_{j_{k-h}=d_{k-h-1}+1}^{d_{k-h}}
\left(\prod_{b=1}^{k-h}\prod_{r_b=d_{b-1}+1}^{j_b-1}(-u_{r_b})\right)\left(\prod_{p=d_{k-h}+1}^{n}(-u_p)\right)v_{\{a_1,\dots,a_{h-1},\,i,\,i+j_1,\,i+j_2,\dots,i+j_{k-h}\}}\\
+\sum_{j_1=d_0+1}^{d_1}\dots\sum_{{j_{k-h}=d_{k-h-1}+1}}^{d_{k-h}}\,\,\sum_{{j_{k-h+1}=d_{k-h}+1}}^{n-i}\left(\prod_{b=1}^{k-h+1}\prod_{r_b=d_{b-1}+1}^{j_b-1}(-u_{r_b})\right)\\
\cdot v_{\{a_1,\dots,a_{h-1},\,i+j_1,\,i+j_2,\dots,\,i+j_{k-h},\,i+j_{k-h+1}\}}.
\end{multline*}
This equals the r.h.s.\,\,of \eqref{id1}. This proves identity \eqref{id1}. The proof of identity \eqref{id2} is similar, and it is left to the reader.
\endproof

\proof[Proof of Theorem \ref{satqkz}]Notice that
\begin{empheq}[left =p_1^{e_{11}^{(i)}}p_2^{e_{22}^{(i)}}v_{[a]}{=}\empheqlbrace]{align*}
&p_1v_{[a]},&\text{if $i=a$,}\\
&p_2v_{[a]},&\text{if $i\neq a$,}
\end{empheq}
so that
\begin{align*}
\bigwedge\nolimits^k\left[p_1^{e_{11}^{(i)}}p_2^{e_{22}^{(i)}}\right]v_{[a_1]}\wedge\dots\wedge v_{[a_k]}=&p_1p_2^{k-1}\,v_{[a_1]}\wedge\dots\wedge v_{[a_k]},&\text{if $i\in\{a_1,\dots, a_k\}$,}\\
\bigwedge\nolimits^k\left[p_1^{e_{11}^{(i)}}p_2^{e_{22}^{(i)}}\right]v_{[a_1]}\wedge\dots\wedge v_{[a_k]}=&p_2^k\,v_{[a_1]}\wedge\dots\wedge v_{[a_k]},&\text{if $i\notin\{a_1,\dots,a_k\}$.}
\end{align*}
Hence, we have
\beq
\label{id3}
\theta_{k,n}\circ\bigwedge\nolimits^k\left[p_1^{e_{11}^{(i)}}p_2^{e_{22}^{(i)}}\right]=p_2^{k-1}\,p_1^{e_{11}^{(i)}}p_2^{e_{22}^{(i)}}\circ\theta_{k,n}.
\eeq
Fix $i=1,\dots, n$. Given $v_{[a_1]}\wedge\dots\wedge v_{[a_k]}$, with $1\leq a_1<\dots<a_k\leq n$, let $h$ be such that $a_{h-1}<i\leq a_h$. We have
\[\bigwedge\nolimits^kK_a[\bm z;(-1)^{k-1}p_1,p_2;\kappa]v_{[a_1]}\wedge\dots\wedge v_{[a_k]}=A_3A_2A_1v_{[a_1]}\wedge\dots\wedge v_{[a_k]},
\]where
\[A_1=\bigwedge\nolimits^k\left[((-1)^{k-h})^{e_{11}^{(i)}}R(z_i-z_n)^{(i,n)}\dots R(z_i-z_{i+1})^{(i,i+1)}\right],\qquad A_2=\bigwedge\nolimits^k\left[p_1^{e_{11}^{(i)}}p_2^{e_{22}^{(i)}}\right],
\]
\[A_3=\bigwedge\nolimits^k\left[R(u_{z_i-z_{i-1}+\kappa})^{(i,i-1)}\dots R(z_i-z_1+\kappa)^{(i,1)}((-1)^{h-1})^{e_{11}^{(i)}}\right].
\]We have
\begin{multline*}A_2A_1(v_{[a_1]}\wedge\dots\wedge v_{[a_k]})
=p_2^{k-1}\theta_{k,n}^{-1} p_1^{e_{11}^{(i)}}p_2^{e_{22}^{(i)}}R(z_i-z_n)^{(i,n)}\dots R(z_i-z_{i+1})^{(i,i+1)}\theta_{k,n}(v_{[a_1]}\wedge\dots\wedge v_{[a_k]}).
\end{multline*}
This follows from \eqref{id3} and Lemma \ref{lem1satqkz}. Moreover, an inspection of the proof of formula \eqref{id1} shows that the r.h.s.\,\,is a linear combination of terms of the form $$v_{[a_1]}\wedge\dots\wedge\underset{h\text{-th}}{v_{[i+j_1]}}\wedge v_{[i+j_2]}\wedge\dots\wedge v_{[i+j_{k-h+1}]},\quad 0\leq j_1<j_2<\dots<j_{k-h+1}\leq n-i.$$
Then, we have
\begin{multline*}\theta_{k,n}A_3A_2A_1(v_{[a_1]}\wedge\dots\wedge v_{[a_k]})
=p_2^{k-1}R(u_{z_i-z_{i-1}+\kappa})^{(i,i-1)}\dots R(z_i-z_1+\kappa)^{(i,1)}\\
p_1^{e_{11}^{(i)}}p_2^{e_{22}^{(i)}}R(z_i-z_n)^{(i,n)}\dots R(z_i-z_{i+1})^{(i,i+1)}\theta_{k,n}(v_{[a_1]}\wedge\dots\wedge v_{[a_k]}),
\end{multline*}
since identity \eqref{id2} can be applied. This proves \eqref{idqkzsat}.
\endproof

\begin{example}
Let $n=4$. The space $(\C^2)^{\otm 4}_{1,3}$ has basis
\[v_{[1]}=v_1\otm v_2\otm v_2\otm v_2,\quad v_{[2]}=v_2\otm v_1\otm v_2\otm v_2,\quad v_{[3]}=v_2\otm v_2\otm v_1\otm v_2,\quad v_{[4]}=v_2\otm v_2\otm v_2\otm v_1.
\]The operator $K_1$ has matrix 
\[\left(
\begin{array}{cccc}
 p_1 \left(z_1-z_2\right) \left(z_1-z_3\right) \left(z_1-z_4\right) & p_1 \left(z_3-z_1\right) \left(z_4-z_1\right) & p_1 \left(z_1-z_4\right) & p_1 \\
 p_2 & 0 & 0 & 0 \\
 p_2 \left(z_1-z_2\right) & p_2 & 0 & 0 \\
 p_2 \left(z_2-z_1\right) \left(z_3-z_1\right) & p_2 \left(z_1-z_3\right) & p_2 & 0 \\
\end{array}
\right)
\]
Consider the space $\bigwedge^2(\C^2)^{\otm 4}_{1,3}$ with basis
\[v_{[1]}\wedge v_{[2]},\quad v_{[1]}\wedge v_{[3]},\quad v_{[1]}\wedge v_{[4]},\quad v_{[2]}\wedge v_{[3]},\quad v_{[2]}\wedge v_{[4]},\quad v_{[3]}\wedge v_{[4]}.
\]
In this basis, the operator $p_2^{-1}\bigwedge^2K_1[\bm z;-p_1,p_2;\kappa]$ has matrix
\[\left(
\begin{array}{cccccc}
 p_1 (z_1-z_3) (z_1-z_4) & p_1 (z_1-z_4) & p_1 & 0 & 0 & 0 \\
 0 & p_1 (z_1-z_2) (z_1-z_4) & p_1 (z_1-z_2) & p_1 (z_1-z_4) & p_1 & 0 \\
 0 & 0 & p_1 (z_1-z_2) (z_1-z_3) & 0 & p_1 (z_1-z_3) & p_1 \\
 p_2 & 0 & 0 & 0 & 0 & 0 \\
 p_2 (z_1-z_3) & p_2 & 0 & 0 & 0 & 0 \\
 0 & p_2 (z_1-z_2) & 0 & p_2 & 0 & 0 \\
\end{array}
\right).
\]
This coincides with the matrix of the operator $K_1$ acting on $(\C^2)^{\otm 4}_{2,2}$ with basis
\[v_{1}\otm v_1\otm v_2\otm v_2,\quad v_1\otm v_2\otm v_1\otm v_2,\quad v_1\otm v_2\otm v_2\otm v_1,\quad v_2\otm v_1\otm v_1\otm v_2,\quad v_2\otm v_1\otm v_2\otm v_1,\quad v_2\otm v_2\otm v_1\otm v_1.
\]
Consider the space $\bigwedge^3(\C^2)^{\otm 4}_{1,3}$ with basis
\[v_{[1]}\wedge v_{[2]}\wedge v_{[3]},\quad v_{[1]}\wedge v_{[2]}\wedge v_{[4]},\quad v_{[1]}\wedge v_{[3]}\wedge v_{[4]},\quad v_{[2]}\wedge v_{[3]}\wedge v_{[4]}.
\]In this basis, the operator $p_2^{-2}\bigwedge^3K_1[\bm z;\bm p;\kappa]$ has matrix
\[\left(
\begin{array}{cccc}
 p_1 \left(z_1-z_4\right) & p_1 & 0 & 0 \\
 0 & p_1 \left(z_1-z_3\right) & p_1 & 0 \\
 0 & 0 & p_1 \left(z_1-z_2\right) & p_1 \\
 p_2 & 0 & 0 & 0 \\
\end{array}
\right).
\]
This coincides with the matrix of $K_1$ acting on $(\C^2)^{\otm 4}_{3,1}$ with respect to the basis
\[v_1\otm v_1\otm v_1\otm v_2,\quad v_1\otm v_1\otm v_2\otm v_1,\quad v_1\otm v_2\otm v_1\otm v_1,\quad v_2\otm v_1\otm v_1\otm v_1.
\]
Consider the space $\bigwedge^4(\C^2)^{\otm 4}_{1,3}$ with basis
\[v_{[1]}\wedge v_{[2]}\wedge v_{[3]}\wedge v_{[4]}.
\]The operator $p_2^{-3}\bigwedge^4K_1[\bm z;-p_1,p_2;\kappa]$ is just the rescaling by $p_1$. This coincides with $K_1$ acting on $(\C^2)^{\otm 4}_{4,0}$.
\vskip2mm
Similar results hold for $K_2,K_3,K_4$.\qetr
\end{example}

\subsection{Dynamical operators and Satake identifications}\label{secqkzdyn}
In this section we prove identities \eqref{gauge2} and \eqref{gauge3}.
\begin{thm}\label{satdyn}
The isomorphism $\theta_{k,n}$ intertwines the operators $\widetilde X_1^{(k)}(\bm z;\bm p)$ acting on $\bigwedge^k_\C\left(\C^2\right)^{\otm n}_{1,n-1}$ and the dynamical operator $X_1(\bm z;\bm p)$ acting on $(\C^2)^{\otm n}_{k,n-k}$, that is
\beq
\label{iddynsat1}
\theta_{k,n}\circ \widetilde X_1^{(k)}=X_1\circ \theta_{k,n}.
\eeq
Moreover, we have
\beq\label{iddynsat2}
\theta_{k,n}\circ\widetilde X_2^{(k)}=\left(X_2+(k-1)\sum_{a=1}^nz_a\right)\circ\theta_{k,n}.
\eeq
\end{thm}

\proof
The operators $Q^{a,b}_{i,j}$ act on $(\C^2)^{\otm n}_{1,n-1}$ as follows:\\
\begin{minipage}{0.5\textwidth}
\begin{empheq}[left =Q^{a,b}_{1,2}v_{[p]}{=}\empheqlbrace]{align*}
v_{[p-1]},\quad&\text{if $a=p,\,b=p-1,\,p>1$,}\\
v_{[n]},\quad&\text{if $a=1,\,b=n,\,p=1$,}\\
0,\quad&\text{otherwise},
\end{empheq}
\end{minipage}
\begin{minipage}{0.5\textwidth}
\begin{empheq}[left =Q^{a,b}_{2,1}v_{[p]}{=}\empheqlbrace]{align*}
v_{[p-1]},\quad&\text{if $a=p-1,\,b=p,\,p>1$,}\\
v_{[n]},\quad&\text{if $a=n,\,b=1,\,p=1$,}\\
0,\quad&\text{otherwise}.
\end{empheq}
\end{minipage}\\

So, for any $j=1,\dots,k$, we have
\[X_1(\bm z;(-1)^{k-1}p_1,p_2)v_{[a_j]}=z_{a_j}v_{[a_j]}+(-1)^{k-1}\frac{p_2}{p_1}v_{[n]}\dl_{1,a_j}+v_{[a_j-1]}(1-\dl_{1,a_j}),
\]where $\dl_{i,j}$ is the Kronecker delta symbol. Then, we have:
\begin{itemize}
\item if $a_1=1$,
\[\widetilde X_1^{(k)}v_{[a_1]}\wedge\dots\wedge v_{[a_k]}=\left(\sum_{j=1}^kz_{a_j}\right)v_{[a_1]}\wedge\dots\wedge v_{[a_k]}+(-1)^{k-1}\frac{p_2}{p_1}v_{[n]}\wedge v_{[a_2]}\wedge\dots\wedge v_{[a_k]};
\]
\item if $a_1\neq 1$, 
\begin{multline*}\widetilde X_1^{(k)}v_{[a_1]}\wedge\dots\wedge v_{[a_k]}\\=\left(\sum_{j=1}^kz_{a_j}\right)v_{[a_1]}\wedge\dots\wedge v_{[a_k]}+\sum_{j=1}^kv_{[a_1]}\wedge\dots\wedge v_{[a_{j-1}]}\wedge v_{[a_j-1]}\wedge v_{[a_{j+1}]}\wedge\dots\wedge v_{[a_k]}.
\end{multline*}
\end{itemize}

If we apply $\theta_{k,n}$, and write for short $v_{\{a_1,\dots,a_k\}}$ rather than $v_{\{a_1,\dots,a_k\},\{1,\dots n\}\setminus\{a_1,\dots,a_k\}}$, we obtain:
\begin{itemize}
\item if $a_1=1$,
\[\theta_{k,n}\left(\widetilde X_1^{(k)}v_{[a_1]}\wedge\dots\wedge v_{[a_k]}\right)=\left(\sum_{j=1}^kz_{a_j}\right)v_{\{a_1,\dots,a_k\}}+\frac{p_2}{p_1}v_{\{a_2,\dots,a_k,n\}};
\]
\item if $a_1\neq 1$,
\[\theta_{k,n}\left(\widetilde X_1^{(k)}v_{[a_1]}\wedge\dots\wedge v_{[a_k]}\right)=\left(\sum_{j=1}^kz_{a_j}\right)v_{\{a_1,\dots,a_k\}}+\sum_{j=1}^kv_{\{a_1,\dots,a_{j-1},a_j-1,a_{j+1},\dots,a_k\}}.
\]
\end{itemize}
This equals $X_1(\bm z;\bm p)\theta_{k,n}(v_{[a_1]}\wedge\dots\wedge v_{[a_k]})$. So identity \eqref{iddynsat1} holds.
\vskip2mm
Similarly, for $j=1,\dots,k$ we have
\[X_2(\bm z;(-1)^{k-1}p_1,p_2)v_{[a_j]}=\left(\sum_{a\neq a_j}z_{a}\right)v_{a_j}-(-1)^{k-1}\frac{p_2}{p_1}v_{[n]}\dl_{1,a_j}-v_{[a_j-1]}(1-\dl_{1,a_j}).
\]Then we have:
\begin{itemize}
\item if $a_1=1$,
\[\widetilde X_2^{(k)}v_{[a_1]}\wedge\dots\wedge v_{[a_k]}=\left(\sum_{j=1}^n\sum_{a\neq a_j}z_{a}\right)v_{[a_1]}\wedge\dots\wedge v_{[a_k]}-(-1)^{k-1}\frac{p_2}{p_1}v_{[n]}\wedge v_{[a_2]}\wedge\dots\wedge v_{[a_k]};
\]\item if $a_1\neq 1$, 
\[\widetilde X_2^{(k)}v_{[a_1]}\wedge\dots\wedge v_{[a_k]}=\left(\sum_{j=1}^n\sum_{a\neq a_j}z_{a}\right)v_{[a_1]}\wedge\dots\wedge v_{[a_k]}-\sum_{j=1}^kv_{[a_1]}\wedge\dots\wedge v_{[a_{j-1}]}\wedge v_{[a_j-1]}\wedge v_{[a_{j+1}]}\wedge\dots\wedge v_{[a_k]}.
\]\end{itemize}
In both cases, this equals $\theta_{k,n}^{-1}\circ X_2\circ \theta_{k,n}(v_{[a_1]}\wedge\dots\wedge v_{[a_k]})$.
This proves \eqref{iddynsat2}.
\endproof

\proof[Proof of Theorem \ref{abstractthmsat}]
The theorem follows from Theorem \ref{satqkz} and Theorem \ref{satdyn}.
\endproof

\section{Equivariant cohomology of Grassmannians}\label{SEC3}
General references for this section are \cite{AF,Ful07,GKM98,KT03}.
\subsection{Ring presentation}\label{repcohgr} Consider the complex Grassmannian $G(k,n)$ parametrizing the $k$-dimensional subspaces $F_k$ of $\C^n$.
\vskip2mm
Given a basis of $\C^n$, the group $GL(n,\C)$ acts on $\C^n$. Let ${T}\subseteq GL(n,\C)$ be the torus of diagonal matrices. 
\vskip2mm
Consider the following ${T}$-equivariant bundles on $G(k,n)$:
\begin{itemize}
\item the trivial vector bundle $\underline{\C^n}=\C^n\times G(k,n)$, with equivariant Chern roots $z_1,\dots,z_n$; 
\item the tautological bundle $E_1$, with fiber $F_k$, and with equivariant Chern roots $\gm_{1,1},\dots,\gm_{1,k}$; 
\item the quotient bundle $E_2$, with fiber $\C^n/F_k$, and with equivariant Chern roots $\gm_{2,1},\dots,\gm_{2,n-k}$. 
\end{itemize}
Let $u$ be a formal variable. For any equivariant vector bundle $V$, denote by $c^{T}(V)=\sum_{j\geq 0} c^T_j(V)u^j$ its equivariant total Chern class. We thus have
\[c^{T}(\underline{\C^n})=\prod_{a=1}^n(1+z_au),\qquad c^{T}(E_1)=\prod_{j=1}^k(1+\gm_{1,j}u),\qquad c^{T}(E_2)=\prod_{j=1}^{n-k}(1+\gm_{2,j}u).
\]
Set $\bm z=(z_1,\dots, z_n)$, ${\bm\gm_1}=(\gm_{1,1},\dots,\gm_{1,k})$, $\bm\gm_2=(\gm_{2,1},\dots,\gm_{2,n-k})$. We have
\[c^T_j(\underline{\C^n})=e_j(\bm z),\qquad c^T_j(E_1)=e_j(\bm\gm_1),\qquad c^T_j(E_2)=e_j(\bm\gm_2),
\]where $e_j$ denotes the $j$-th elementary symmetric polynomial.
\vskip2mm
Set $\bm\gm:=\left(\gm_{1,1},\dots,\gm_{1,k},\gm_{2,1},\dots,\gm_{2,n-k}\right)$, and denote by $\C[\bm\gm]^{S_k\times S_{n-k}}$ the space of polynomials in $\bm\gm$ symmetric in $\gm_{1,1},\dots,\gm_{1,k}$ and $\gm_{2,1},\dots,\gm_{2,n-k}$, separately.
\vskip2mm
The bundles $E_1,E_2$ fit in the short exact sequence
\beq\label{Euler}
0\to E_1\to\underline{\C^n}\to E_2\to 0.
\eeq
From this, the equivariant cohomology algebra $H_T^\bullet(G(k,n),\C)$ admits the following presentation
\beq
\label{presentation}
H_T^\bullet(G(k,n),\C)\cong \C[\bm\gm]^{S_k\times S_{n-k}}\otimes\C[\bm z]\Big/\Big\langle c^{T}(E_1)c^{T}(E_2)=c^{T}(\underline{\C^n})\Big\rangle.
\eeq
Here we mean that the ideal is generated by the relations obtained by equating all the powers $u^i$, with $i=0,\dots,n$, in the identity $c^{T}(E_1)c^{T}(E_2)=c^{T}(\underline{\C^n})$. That is 
\[\sum_{j=0}^he_j(\bm\gm_1)e_{h-j}(\bm \gm_2)=e_h(\bm z),\quad h=1,\dots,n.
\]

\subsection{Partitions, Grassmannian permutations, and Young diagrams} \label{secpgrpy}
Recall the definition of the set $\mc I_{k,n}$ of all double partitions $I=(I_1,I_2)$ of $\{1,\dots,n\}$ with ${\rm card\,} I_1=k,$ and ${\rm card\,} I_2=n-k$.
We denote by $I^{\rm min}_k$ the partition
\[I^{\rm min}_k=\left(\{1,\dots,k\},\quad\left\{k+1,\dots,n\right\}\right).
\]Permutations in $S_n$ act on $\mc I_{k,n}$, via $\si(I)=\left(\si(I_1),\si(I_2)\right)$. Such an action is transitive but not free: it descends to a free transitive action of $S_n/(S_k\times S_{n-k})$ on $\mc I_{k,n}$.
It is well known that minimal length representatives for cosets in $S_n/(S_k\times S_{n-k})$ are the so-called {\it Grassmannian permutations with descent at $k$}: these are the partitions $\si\in S_n$ such that
\[\si(1)<\si(2)<\dots<\si(k),\quad \si(k+1)<\si(k+2)<\dots<\si(n).
\]We denote by $GrS_{k,n}$ the set of all such permutations.
Notice that $\si\in GrS_{k,n}$ is uniquely determined by the values $(\si(1),\dots,\si(k))$. \vskip2mm
Following \cite{TV23}, given $I\in\mc I_{k,n}$, we denote by $\si_I$ the unique permutation in $GrS_{k,n}$ such that $\si_I(I^{\rm min}_k)=I$. Conversely, given $\si\in GrS_{k,n}$, we can associate to it an element $I^\si\in\mc I_{k,n}=(I^\si_1,I^\si_2)$ by setting $I_1^\si:=\{\si(1),\dots,\si(k)\}$.
\begin{lem}
The maps $I\mapsto \si_I$ and $\si\mapsto I^\si$ define a bijection.\qed
\end{lem}
\vskip2mm
Let us construct an explicit bijection between $GrS_{k,n}$ and the set $P_{k,n}$ of partitions $\bm\la\subseteq k\times (n-k)$.
Given a Grassmannian permutation $\si\in GrS_{k,n}$, we associate to it a partition $\bm\la^\si\subseteq k\times (n-k)$ defined by
\[\la_j^\si=\si(k-j+1)-k+j-1,\quad j=1,\dots,k.
\]Conversely, given a partition $\bm\la\subseteq k\times (n-k)$, we can associate to it the Grassmannian permutation $\si^{\bm\la}\in GrS_{k,n}$ defined by
\[
\si^{\bm\la}(j)=\la_{k-j+1}+j,\quad j=1,\dots, k,\qquad \si^{\bm\la}(k+j)=-\la'_j+j+k,\quad j=1,\dots,n-k,
\]where $\bm\la'$ denotes the dual partition of $\bm\la$. In other words, to construct $\si^{\bm \la}$ consider the Young diagram of $\bm \la\subseteq k\times (n-k)$, and starting form the SW corner, walk to the NE corner. Define the tuple $(i_1,\dots, i_k)$ as follows: 
\begin{itemize}
\item if the $i_h$-th step is in upper direction (vertical one), then $i_h$ is the $h$-th element of the tuple, 
\item if $i_h$ is in right direction (horizontal), then $i_h$ is not in the tuple.
\end{itemize}
The permutation $\si^{\bm \la}\in GrS_{k,n}$ is the one defined by $\si^{\bm \la}(h)=i_{h}$ for $h=1,\dots,k$.
\begin{lem}
The maps $\si\mapsto \bm\la^\si$ and $\bm\la\mapsto\si^{\bm\la}$ define a bijection.\qed
\end{lem}

\subsection{Fixed points, and localization} The ${T}$-action on $G(k,n)$ admits $\binom{n}{k}$ fixe points ${\rm pt}_I$, labelled by $I\in \mc I_{k,n}$. If $(v_1,\dots, v_n)$ denotes the standard basis of $\C^n$, with $v_i=(0,\dots,0,1_i,0,\dots,0)$ for $i=1,\dots, n$, the fixed point ${\rm pt}_I$ is the $k$-dimensional subspace ${\rm span}(v_{i_1},\dots, v_{i_k})$, with $I_1=\{i_1,\dots, i_k\}$.
\vskip2mm
If we denote by $G(k,n)^{T}$ the set of fixed points, the inclusion $\iota\colon G(k,n)^{T}\to G(k,n)$ is ${T}$-equivariant. Then, we have an induced pullback map $\iota^*\colon H^\bullet_T(G(k,n),\C)\to H^\bullet_T(G(k,n)^{T},\C)\cong \bigoplus_{I\in\mc I_{k,n}}H^\bullet_T({\rm pt}_I,\C)$. This map is injective \cite[Th.\,1.2.2]{GKM98}, so that any class $\al$  is determined by its image $\iota^*\al$. Moreover, for each $I\in \mc I_{k,n}$, the inclusion $\iota_I\colon \{{\rm pt}_I\}\to G(k,n)$ is ${T}$-equivariant, and it induces a pullback map. The image $\iota_I^*\al\in H^\bullet_T({\rm pt}_I,\C)$ is called the {\it localization}, or {\it restriction}, of the class $\al$ at the point ${\rm pt}_I$. It is standard to denote by $\al|_{{\rm pt}_I}$ the localization $\iota_I^*\al$.
\vskip2mm
In terms of the polynomial presentation \eqref{presentation}, the restriction to a fixed point ${\rm pt}_I$ is given by the evaluation map
\[
f(\bm\gm;\bm z)\mapsto f(\bm z_{\si_I};\bm z),\quad \bm z_{\si_I}=(z_{\si_I(1)},\dots,z_{\si_I(n)}),\quad I\in\mc I_{k,n}.
\]
The fundamental classes $[{\rm pt}_I]$ of fixed points generate the localization of the ring $H^\bullet_T(G(k,n),\C)$ at $H^\bullet_T({\rm pt},\C)\cong \C[\bm z]$: for each cohomology class $\al$, we have
\[\al=\sum_{I\in\mc I_{k,n}}\al|_{{\rm pt}_I}\frac{[{\rm pt}_I]}{[{\rm pt}_I]|_{{\rm pt}_I}}.
\]See Theorem \ref{thlocal} and Remark \ref{remfixedpt} below.

\subsection{Extension of scalars and idempotent basis}\label{secext} Let $\bm x=(x_1,\dots,x_n)$ be a $n$-tuple of formal variables, and $\C[\bm x]^{S_k\times S_{n-k}}$ the ring of polynomials symmetric in $(x_1,\dots,x_k)$ and $(x_{k+1},\dots,x_n)$ separately.
Consider the polynomials
\beq\label{defVR}
V_{k,n}(\bm x;\bm z)=\prod_{j=1}^k\prod_{a=k+1}^n(x_j-z_a),\quad R_{k,n}(\bm z):=V_{k,n}(\bm z;\bm z)=\prod_{j=1}^k\prod_{a=k+1}^n(z_j-z_a).
\eeq
Given $\si\in S_n$, set $\bm z_\si=(z_{\si(1)},\dots, z_{\si(n)})$. For each $I\in\mc I_{k,n}$ define
\beq
\Dl_I(\bm x;\bm z):=\frac{V_{k,n}(\bm x;\bm z_{\si_I})}{R_{k,n}(\bm z_{\si_I})} \in\C[\bm x]^{S_k\times S_{n-k}}\otm\C(\bm z).
\eeq
\begin{lem}\cite[Th.\,2.1]{CL96}\label{lemidem}
For any $f\in \C[x_1,\dots,x_k]^{S_k}$, with partial degrees $\deg_{x_i}f\leq n-k$, we have
\beq\label{idemid0}
\pushQED{\qed} 
f(\bm x)=\sum_{I\in\mc I_{k,n}}f(\bm z_{\si_I})\Dl_I(\bm x;\bm z).\qedhere
\popQED
\eeq
\end{lem}

Given $\kappa\in\C^*$, let ${L_\kappa}$ be the complement in $\C^n$ of the union of hyperplanes 
\[z_i-z_j=\kappa m,\quad i,j=1,\dots,n,\quad m\in\Z.
\]Denote by $\mc O_{{L_\kappa}}$ the ring of holomorphic functions on ${L_\kappa}$. This ring is a module over the ring $H^\bullet_{{T}}({\rm pt},\C)\cong \C[\bm z]$. Set
\beq
\label{hext}
H^{(\kappa)}_T(G(k,n)):=H^\bullet_T(G(k,n),\C){\otm_{H^\bullet_T({\rm pt},\C)}}\mc O_{{L_\kappa}}.
\eeq
The polynomials $\Dl_I(\bm \gm;\bm z)$ define classes in $H^{(\kappa)}_{{T}}(G(k,n))$, denoted by the same symbols.

\begin{thm}\label{thlocal}
The polynomials $\Dl_I$, with $I\in \mc I_{k,n}$, define a basis of idempotent elements of $H^{(\kappa)}_{{T}}(G(k,n))$, that is $\Dl_I\Dl_J=\Dl_I\dl_{IJ}$. For each $f\in H^{(\kappa)}_{{T}}(G(k,n))$, we have
\beq\label{idemid0}
f(\bm \gm;\bm z)=\sum_{I\in\mc I_{k,n}}f(\bm z_{\si_I};\bm z)\Dl_I(\bm\gm;\bm z).
\eeq
\end{thm}
\proof
This follows from Lemma \ref{lemidem}. Notice indeed that any cohomological class $f(\bm \gm;\bm z)$ can be represented by a symmetric polynomial in $\bm \gm_1$ only, with partial degrees at most $n-k$.
\endproof

\begin{rem}\label{remfixedpt}
The classes $\Dl_I(\bm\gm;\bm z)\in H^\bullet_T(G(k,n),\C)$, with $I\in\mc I_{k,n}$, are the classes representing the fixed points ${\rm pt}_I$.
\end{rem}

\subsection{Equivariant Poincar\'e pairing}\label{poimet} The map $\pi\colon G(k,n)\to \{\rm pt\}$, induces a push-forward map $\pi_*\colon H^\bullet_T(G(k,n),\C)\to H^\bullet_T({\rm pt}, \C)\cong \C[\bm z]$. The map $\pi_*$ defines the {\it equivariant integration} of cohomological classes. In terms of the polynomial presentation \eqref{presentation}, it is given by the following Atiyah--Bott integration formula \cite{AB84}.
\vskip2mm
Given any permutation $\si\in S_n$, set $\bm z_\si:=(z_{\si(1)},\dots, z_{\si(n)})$. For each equivariant cohomology class $\al(\bm\gm;\bm z)$ in $H^\bullet_T(G(k,n),\C)$, we have\footnote{Notice that the metric considered in \cite{CV21,TV23} differs from the metric $\eta$ by a sign $(-1)^{k(n-k)}$. This sign is dictated by the identity $[{\rm pt}_I]_{{\rm pt}_I}=c_{\rm top}^{T}(TG(k,n))|_{{\rm pt}_I}=c^T_{\rm top}(E_1^*\otimes E_2)|_{{\rm pt}_I}=(-1)^{k(n-k)}R_{k,n}(\bm z_{\si_I})$.}
\beq
\pi_*\al=\int^{\rm eq}_{G(k,n)}\al=(-1)^{k(n-k)}\sum_{I\in\mc I_{k,n}}\frac{\al(\bm z_{\si_I};\bm z)}{R_{k,n}(\bm z_{\si_I})},\quad R_{k,n}(\bm z):=V_{k,n}(\bm z;\bm z)=\prod_{j=1}^k\prod_{a=k+1}^n(z_j-z_a).
\eeq
The {\it equivariant Poincar\'e metric} on $G(k,n)$ is the $H^\bullet_T({\rm pt},\C)$-bilinear map
\[\eta\colon H^\bullet_T(G(k,n),\C)\times H^\bullet_T(G(k,n),\C)\to H^\bullet_T({\rm pt},\C)\cong \C[\bm z]
\]defined by
\[\eta(v_1,v_2)=\int^{\rm eq}_{G(k,n)}v_1\cdot v_2=(-1)^{k(n-k)}\sum_{I\in\mc I_{k,n}}\frac{v_1(\bm z_{\si_I},\bm z)v_2(\bm z_{\si_I},\bm z)}{R_{k,n}(\bm z_{\si_I})}.
\]

\subsection{Schubert basis} 
In this section we introduce a basis of $H_T^\bullet(G(k,n),\C)$ as a $\C[\bm z]$-module, called the {\it equivariant Schubert basis}.
It consists of a basis $([\Om_{\bm\la}])_{\bm\la}$, parametrized by partitions $\bm\la$ contained in a $k\times(n-k)$ Young diagram, that is $$\bm\la\colon \la_1\geq \la_2\geq\dots\geq\la_k\geq 0,\quad \la_1\leq n-k.$$ Its construction depends on the choice of a complete flag in $\C^n$. Let $(v_1,\dots, v_n)$ be the standard basis of $\C^n$, with $v_i=(0,\dots,0,1_i,0,\dots,0)$ for $i=1,\dots, n$. Set $L_i:={\rm span}(v_1,\dots, v_i)$, for $i=1,\dots, n$, and consider the standard flag
$L_\bullet\colon\,\, 0=L_0\subset L_1\subset L_2\subset\dots\subset L_n=\C^n.$ 

Given $\bm\la$, the Schubert class $[\Om_{\bm\la}]\in H_T^\bullet(G(k,n),\C)$ is the cohomology class representing the Schubert variety $\Om_{\bm\la}$ defined by
\[\Om_{\bm\la}:=\left\{F\in G(k,n)\colon \dim\left(F\cap L_{n-k+i-\la_i}\right)\geq i,\quad i=1,\dots, k\right\}.
\]It is easily seen that the varieties $\Om_{\bm\la}$ are ${T}$-invariant, so that they can be represented by equivariant cohomological classes.
\vskip2mm
We can similarly introduce the {\it opposite} equivariant Schubert basis $([\widetilde\Om_{\bm\la}])_{\bm\la}$. For $i=1,\dots, n$, set $\widetilde L_i:={\rm span}(v_{n-i+1}, v_{n-i+2},\dots, v_n)$ for $ i=1,\dots, n$. Consider the opposite complete flag $\widetilde L_\bullet$ of $\C^n$, $\widetilde L_\bullet\colon\,\, 0=\widetilde L_0\subset \widetilde L_1\subset \widetilde L_2\subset\dots\subset \widetilde L_n=\C^n$.

The class $[\widetilde\Om_{\bm\la}]$ is the cohomology class representing the opposite Schubert variety
\[\widetilde \Om_{\bm\la}:=\left\{F\in G(k,n)\colon \dim\left(F\cap \widetilde L_{n-k+i-\la_i}\right)\geq i,\quad i=1,\dots, k\right\}.
\]
Given a partition $\bm\la\subseteq k\times (n-k)$, denote by $\bm\la^\vee$ the {\it complementary} partition, that is the partition $\bm\la^\vee\subseteq k\times (n-k)$ such that
\[\la_j+\la^\vee_{k+1-j}=n-k,\quad j=1,\dots, k.
\]

\begin{thm}
The classes $[\Om_{\bm\la}]$ and $[\Tilde\Om_{\bm\la^\vee}]$ are Poincar\'e dual:
\[\pushQED{\qed} 
\eta([\Om_{\bm\la}], [\Tilde\Om_{\bm\mu}])=\dl_{\bm\mu,\bm\la^{\vee}}.\qedhere
\popQED\]
\end{thm}For proofs and more details, see e.g.\,\,\cite{AF,Ful07}.
\subsection{Kempf--Laksov--Lascoux formula and Schur polynomials} Given two ${T}$-equivariant bundles $E,F$, the Chern classes of the $K$-theoretical difference $E-F$ are defined via the formula
\begin{align*}c^{T}(E-F)=\frac{c^{T}(E)}{c^{T}(F)}&=\frac{1+c_1^{T}(E)t+c_2^{T}(E)t^2+\dots}{1+c_1^{T}(F)t+c_2^{T}(F)t^2+\dots}\\
&=1+\left(c_1^{T}(E)-c_1^{T}(F)\right)t+\left(c_2^{T}(E)-c_2^{T}(F)-c_1^{T}(E)c_1^{T}(F)+c_1^{T}(F)^2\right)t^2+\dots
\end{align*}
Denote by $\underline{L_j}$ and $\underline{\widetilde L_j}$, with $j=1,\dots, n$, the trivial bundles on $G(k,n)$ with fibers $L_j$ and $\widetilde L_j$, respectively. Set
\begin{align*}
c_j(i)&=c_j^{T}(E_2-\underline{L_{n-k+i-\la_i}}),& i=1,\dots, k,\quad j\geq 0,\\
\tilde c_j(i)&=c_j^{T}(E_2-\underline{\widetilde L_{n-k+i-\la_i}}),& i=1,\dots, k,\quad j\geq 0.
\end{align*}
The following result expresses the Schubert classes $[\Om_{\bm\la}]$ and $[\widetilde\Om_{\bm\la}]$ in terms of the presentation \eqref{presentation}, that is in terms of the classes $c^T_i(E_j)$, with $j=1,2$. Such a formula was obtained by G.\,Kempf and D.\,Laksov in \cite{KL74}, and independently by A.\,Lascoux in \cite{Las74}.

\begin{thm}[Kempf--Laksov--Lascoux formula]
We have
\[[\Om_{\bm\la}]=\det\left(c_{\la_i+j-i}(i)\right)_{1\leq i,j\leq \ell(\la)}=\left|\begin{array}{cccc}
c_{\la_1}(1)&c_{\la_1+1}(1)&c_{\la_1+2}(1)&\dots\\
c_{\la_2-1}(2)&c_{\la_2}(2)&c_{\la_2+1}(2)&\dots\\
\vdots&\ddots&\ddots&\\
\\
&&c_{\la_{\ell(\la)}-1}(\ell(\la))&c_{\la_{\ell(\la)}}(\ell(\la))
\end{array}\right|,
\]
and 
\[[\widetilde\Om_{\bm\la}]=\det\left(\tilde c_{\la_i+j-i}(i)\right)_{1\leq i,j\leq \ell(\la)}=\left|\begin{array}{cccc}
\tilde c_{\la_1}(1)&\tilde c_{\la_1+1}(1)&\tilde c_{\la_1+2}(1)&\dots\\
\tilde c_{\la_2-1}(2)&\tilde c_{\la_2}(2)&\tilde c_{\la_2+1}(2)&\dots\\
\vdots&\ddots&\ddots&\\
\\
&&\tilde c_{\la_{\ell(\la)}-1}(\ell(\la))&\tilde c_{\la_{\ell(\la)}}(\ell(\la))
\end{array}\right|.
\]
\qed
\end{thm}

\begin{example}\label{ex1}
For $\bm\la=(1)$, the Kempf--Laksov formula is 
\begin{align*}
[\Om_{(1)}]=c_1^{T}(E_2-\underline{L_{n-k}})=c_1^{T}(E_2)-c_1^{T}(\underline{L_{n-k}})=c_1^{T}(E_2)-\sum_{i=1}^{n-k}z_i=-c_1^{T}(E_1)+\sum_{i=n-k+1}^nz_i.
\end{align*}
Similarly, we have
\[[\widetilde \Om_{(1)}]=c_1^{T}(E_2-\underline{\widetilde L_{n-k}})=c_1^{T}(E_2)-c_1^{T}(\underline{\widetilde L_{n-k}})=c_1^{T}(E_2)-\sum_{i=k+1}^{n}z_i=-c_1^{T}(E_1)+\sum_{i=1}^kz_i.
\]\qetr
\end{example}
\begin{example}\label{ex2} Consider the case $k=1$, that is the case of $\Pb^{n-1}=G(1,n)$. If $\bm\la=(a)$, with $a\leq n-1$, we have $[\Om_{{\bm\la}}]=c_a^{T}(E_2-\underline{L_{n-a}})$. To evaluate this Characteristic class, first notice that, in virtue of the Euler sequence \eqref{Euler}, we can rewrite the $K$-theoretical difference as follows:
\[E_2-\underline{L_{n-a}}=\underline{\C^n}-E_1-\underline{L_{n-a}}=\underline{\C^n/L_{n-a}}-E_1.
\]We can now invoke a simplifying formula:
\[c_a^{T}(E-L)=c_a^{T}(E\otimes L^*),\quad E\text{ vector bundle, $L$ line bundle.} 
\]Indeed, $E_1=\mc O(-1)$ is a line bundle, and we have
\[[\Om_{(a)}]=c_a^{T}(\underline{\C^n/L_{n-a}}-E_1)=c_a^{T}\left(\underline{\C^n/L_{n-a}}\otimes E_1^*\right)=\prod_{i=n-a+1}^n\left(-c_1^{T}(E_1)+z_i\right).
\]
Similarly, we have
\[[\widetilde\Om_{(a)}]=c_a^{T}\left(\underline{\C^n/\widetilde L_{n-a}}-E_1\right)=c_a^{T}\left(\underline{\C^n/\widetilde L_{n-a}}\otimes E_1^*\right)=\prod_{i=1}^a\left(-c_1^{T}(E_1)+z_i\right).
\]\qetr
\end{example}

\begin{rem}
The basis $\left((-1)^{n-1}[\Om_{(n-1)}],(-1)^{n-2}[\Om_{(n-2)}],\dots,-[\Om_{(1)}],[\Om_{(0)}]\right)$ coincides with the $\bm g$-basis of $H^\bullet_T(\Pb^{n-1},\C)$ considered in \cite{CV21}.
\end{rem}

We can express the Schubert classes as polynomials in $(\bm\gm_1,\bm z)$ only. For that purpose, let us introduce a set of symmetric polynomials called {\it factorial Schur polynomials}. Given $\bm x=(x_1,\dots,x_k)$ and $\bm y=(y_1,y_2,y_3,\dots)$ formal variables, for any partition $\bm \la\subseteq k\times(n-k)$ we define the  factorial Schur polynomial $s_{\bm\la}(\bm x|\bm y)$ by the determinantal formula
\beq\label{facschurpol}
s_{\bm\la}(\bm x|\bm y):=\frac{\det\left((x_i|\bm y)^{\la_j+k-j}\right)_{1\leq i,j\leq k}}{\det\left((x_i|\bm y)^{k-j}\right)_{1\leq i,j\leq k}},\quad (x_i|\bm y)^p:=\prod_{a=1}^p(x_i-y_a),\quad p\geq 0.
\eeq
\begin{rem}
We have $\det\left((x_i|\bm y)^{k-j}\right)_{1\leq i,j\leq k}=\prod_{1\leq i<j\leq k}(x_i-x_j)$.
\end{rem}
Given any permutation $\si\in S_n$, set $\bm z_\si:=(z_{\si(1)},\dots, z_{\si(n)})$. In particular, if $\si_0\in S_n$ denotes the longest permutation, we set $\bm z_{\si_0}=(z_n,z_{n-1},\dots,z_1)$. The following result is well known, see for example \cite{Mih08}.
\begin{thm}\label{schubertclassesschurpoly}
We have
\begin{align*}
[\Om_\la]&=s_\la(-\bm\gm_1|-\bm z_{\si_0})=(-1)^{|\la|}s_\la(\bm\gm_1|\bm z_{\si_0}),\\
[\widetilde\Om_\la]&=s_\la(-\bm\gm_1|-\bm z)=(-1)^{|\la|}s_\la(\bm\gm_1|\bm z).&
\end{align*}where $|\la|:=\sum_i\la_i$\qed
\end{thm}

\subsection{Double Schubert polynomials and stable envelope maps} Following \cite{TV23}, we introduce a further family of  bases of $H^\bullet_T(G(k,n),\C)$, the {\it stable envelope bases}. In order to define them, let us first recall the notion of {\it double Schubert polynomials} \cite{LS82a,Las12}. 
\vskip2mm
Let $\bm x=(x_1,\dots,x_n)$ and $\bm y=(y_1,y_2,\dots)$ be formal variables. Define the operators $\Dl^{\bm x}_1,\dots, \Dl_{n-1}^{\bm x}$, acting on functions of $\bm x$, by
\[\Dl_a^{\bm x}f(\bm x)=\frac{f(\bm x)-f(x_1,\dots,x_{a+1},x_a,\dots, x_n)}{x_a-x_{a+1}}.
\]They satisfy the nil-Coxeter relations
\beq\label{nilcox}
(\Dl_a^{\bm x})^2=0,\quad \Dl_a^{\bm x}\Dl_b^{\bm x}=\Dl_b^{\bm x}\Dl_a^{\bm x},\quad |a-b|>1,\quad \Dl_a^{\bm x}\Dl_{a+1}^{\bm x}\Dl_a^{\bm x}=\Dl_{a+1}^{\bm x}\Dl_{a}^{\bm x}\Dl_{a+1}^{\bm x},
\eeq
for any $a,b$.
\vskip2mm
Given $\si\in S_n$, define the operator $\Dl_\si^{\bm x}$ as follows. For $a=1,\dots,n-1$, let $\si=s_a$ is the transposition of $a$ and $a+1$. If $\si=s_{a_1}s_{a_2}\dots s_{a_k}$ is the reduced presentation, set $\Dl_{\si}^{\bm x}=\Dl_{a_1}^{\bm x}\Dl_{a_2}^{\bm x}\dots\Dl_{a_k}^{\bm x}$. In particular, $\Dl_{s_a}^{\bm x}=\Dl_a^{\bm x}$. The operators $\Dl_\si^{\bm x}$ are well defined because of the relations \eqref{nilcox}.
\vskip2mm
The $A$-type Schubert polynomials $\frak S_\si$, with $\si\in S_n$, are defined as follows:
\begin{enumerate}
\item for the longest permutation $\si_0\in S_n$, that is $\si_0(i):=n+1-i$ for any $i=1,\dots,n$, set
\[\frak S_{\si_0}(\bm x;\bm y):=\prod_{i+j\leq n}(x_i-y_j);
\]
\item for any $\si\in S_n$, set $\frak S_\si(\bm x;\bm y)=\Dl^{\bm x}_{\si^{-1}\si_0}\frak S_{\si_0}(\bm x;\bm y)$.
\end{enumerate}
\vskip2mm
We are now ready to define the {\it stable envelope maps}. 
Given $I\in\mc I_{k,n}$, define the polynomials
\[ {\rm Stab}_I(\bm\gm;\bm z):=\frak S_{\si_{\si_0(I)}}(\bm\gm;\bm z_{\si_0}),\qquad {\rm Stab}^{\rm op}_I(\bm\gm;\bm z):=\frak S_{\si_{I}}(\bm\gm;\bm z),\quad \text{elements of }\C[\bm\gm]^{\rm sym}\otimes\C[\bm z].
\]
The polynomials ${\rm Stab}_I(\bm\gm;\bm z)$ and ${\rm Stab}^{\rm op}_I(\bm\gm;\bm z)$ define two  classes in $H_T^\bullet(G(k,n),\C)$, that we denote by ${\rm Stab}_I$ and ${\rm Stab}^{\rm op}_I$. 
The following result clarifies the relations between these classes and equivariant Schubert classes: 
\vskip2mm
\begin{prop}
We have
\[{\rm Stab}_I=s_{\bm\la^\vee}(\bm\gm_1;\bm z_{\si_0})=(-1)^{|\bm\la^\vee|}[\Om_{\bm\la^\vee}],\quad {\rm Stab}^{\rm op}_I=s_{\bm\la}(\gm_1;\bm z)=(-1)^{|\bm\la|}[\Tilde\Om_{\bm\la}],
\]where $\bm\la=\bm\la^{\si_I}$. The classes ${\rm Stab}_I$ and ${\rm Stab}^{\rm op}_I$ 
satisfy
\beq\label{stabort}
\eta\left({\rm Stab}_I,{\rm Stab}^{\rm op}_J\right)=(-1)^{k(n-k)}\dl_{I,J}.
\eeq
\end{prop}
\proof These identities immediately follow from Theorem \ref{schubertclassesschurpoly}, and the well known fact that double Schubert polynomials attached with Grassmannian permutations reduce to factorial Schur polynomials, see for example \cite[Ch.\,10]{AF} \cite[Sec.\,1]{CLL02}. More precisely, we have the following lemma.
\begin{lem}
For $\si\in GrS_{k,n}$, we have $\frak S_\si(\bm x;\bm y)=s_{\bm\la}(\bm x|\bm y)$, where $\bm\la=\bm\la^\si$.\qed
\end{lem}
This directly implies the identity for ${\rm Stab}^{\rm op}_I$. Moreover, notice that if $\bm\la=\bm\la^{\si_I}$, then $\bm\la^\vee=\bm\la^{\si_{\si_0(I)}}$, and the identity for ${\rm Stab}_I$ follows. Equation \eqref{stabort} follows from orthogonality properties of double Schubert polynomials, see \cite{TV23}.
\endproof
\begin{rem}
In \cite{CV21,TV23}, the Poincar\'e metric has a different normalization: it differs by the current one by a factor $(-1)^{k(n-k)}$, see Section \ref{poimet}. In the metric of \cite{TV23}, the classes ${\rm Stab}_I$ and ${\rm Stab}^{\rm op}_I$ are Poincar\'e dual.
\end{rem}
Recall the weight decomposition $(\C^2)^{\otm n}=\bigoplus_{j=0}^n(\C^2)^{\otm n}_{j,n-j}$ from Section \ref{sec1.1}. Define the {\it stable envelope map} ${\rm Stab}_{k,n}$ by
\beq
\label{stabkn}
{\rm Stab}_{k,n}\colon (\C^2)^{\otm n}_{k,n-k}\otm \C[\bm z]\to H^\bullet_T(G(k,n),\C),\quad v_I\mapsto {\rm Stab}_I,\quad I\in\mc I_{k,n}.
\eeq

\begin{lem}
The map ${\rm Stab}_{k,n}$ is an isomorphism of free $\C[\bm z]$-modules.
\end{lem}
\proof
The statement follows from \eqref{stabort}. For any $f\in\C[\bm\gm]^{S_k\times S_{n-k}}\otm \C[\bm z]$, one has
\[
\pushQED{\qed}
({\rm Stab}_{k,n})^{-1}\colon [f]\mapsto\sum_{J\in\mc I_{k,n}}\eta\left(f(\bm\gm;\bm z),{\rm Stab}^{\rm op}_J(\bm\gm;\bm z)\right)v_J.\qedhere
\popQED
\]

\section{Equivariant derived category, $K$-theory, and exceptional collections for Grassmannians}\label{SEC4}
General references for this section are \cite{GM03,GK04}\cite[Sec.\,2]{CV21}\cite[Ch.\,5]{CG97}\cite[Sec.\,15.1]{Ful98}.

\subsection{Derived category} We denote by 
\begin{itemize}
\item ${\on{Rep}_\C}({T})$ the category of finite dimensional complex representations of the torus ${T}=(\C^*)^n$,
\item by $R({T}):=K_0({\on{Rep}_\C}({T}))$ the ring of finite dimensional complex representations of ${T}$,
\item by $R({T})_\C$ the complexification $R({T})\otm_\Z\C$.
\end{itemize}
In particular, if we set $\bm Z=(Z_1,\dots, Z_n)$ and $\C[\bm Z^{\pm 1}]=\C[Z_1^{\pm 1},\dots, Z_n^{\pm 1}]$, we have $R({T})_\C=\C[\bm Z^{\pm 1}]$. Each variable $Z_i$ corresponds to a factor of ${T}$, and $\bm Z$ are $K$-theoretic Chern roots (see Section \ref{Kth}) of the representation of ${T}$ defined by its diagonal action on 
$\C^n$.
\vskip2mm
We consider the derived category ${\mathsf D}^b(G(k,n))$ of coherent sheaves on $G(k,n)$, together with its equivariant analog ${\mathsf D}^b_T(G(k,n))$, that is the derived category of ${T}$-equivariant coherent sheaves on $G(k,n)$. All standard functors (such as left derived tensor product, and left and right derived functors induced by morphisms of varieties) are well defined on ${\mathsf D}^b_T(G(k,n))$, see \cite[Ch.\,5]{CG97}.
\vskip2mm
We denote their Grothendieck groups by $K_0(G(k,n))$ and $K_0^{T}(G(k,n))$, respectively. Their complexifications will be denote by $K_0(G(k,n))_\C$ and $K_0^{T}(G(k,n))_\C$, respectively.
\vskip2mm
The structural morphism $\pi\colon G(k,n)\to{\rm Spec}\,\C$ endows $K_0(G(k,n))_\C$ and $K^T_0(G(k,n))_\C$ with a $\C$-algebra and a $R({T})_\C$-algebra structures, respectively. Moreover, we have two push-forward morphisms 
\[\pi_*\colon K_0(G(k,n))\to K_0({\rm Spec}\,\C)\cong\Z,\qquad \pi_*^{T}\colon K_0^{T}(G(k,n))\to K_0^{T}({\rm Spec}\,\C)\cong R({T})=\Z[\bm Z^{\pm 1}],
\]and functors
\[R\pi_*\colon {\mathsf D}^b(G(k,n))\to {\mathsf D}^b(\C),\qquad R\pi_*^{T}\colon {\mathsf D}^b_T(G(k,n))\to {\mathsf D}^b({\on{Rep}_\C}({T})),
\]fitting in the commutative diagram
\[\xymatrix{
{\mathsf D}^b_T(G(k,n))\ar[r]^{\,\, R\pi_*^{T}}\ar[d]_{\frak F_{G(k,n)}}&{\mathsf D}^b({\on{Rep}_\C}({T}))\ar[d]^{\frak F_{\rm pt}}\\
{\mathsf D}^b(G(k,n))\ar[r]^{\,\,\,R\pi_*}&{\mathsf D}^b(\C)
}
\]where $\frak F_{G(k,n)},\frak F_{\rm pt}$ denote the forgetful functors. A ${T}$-{\it equivariant structure} on $V\in{\rm Ob}\left({\mathsf D}^b(G(k,n))\right)$ is an object $V'\in{\rm Ob}\left({\mathsf D}^b_T(G(k,n))\right)$ such that $\frak F_{G(k,n)}(V')=V$.

\subsection{Equivariant $K$-theory, and Grothendieck polynomials}\label{Kth} Given a ${T}$-equivariant complex vector bundle $V$ on $G(k,n)$, we denote by $[V]$ its class in the Grothendieck group $K_0^{T}(G(k,n))$. Following F.\,Hirzebruch \cite{Hir95}, we associate with $V$ a $K$-theoretic analog of the equivariant Chern polynomial, called {\it equivariant} {\it Hirzebruch's} $\lambda$-{\it class}. This is the polynomial $\la^{T}(V)\in K_0^{T}(G(k,n))[u]$, in a formal variable $u$, defined by
\[\la^{T}(V)=\sum_{p\geq 0}[\wedge^p V]u^p.
\]Virtual line bundles $\xi_1,\dots,\xi_{{\rm rk} V}$ such that $\la^{T}(V)=\prod_{i=1}^{{\rm rk}V}(1+\xi_iu)$ will be called $K$-{\it theoretic Chern roots} of $V$. 
The $\la$-class is multiplicative: given a short exact sequence of ${T}$-equivariant bundles 
\[0\to V_1\to V_2\to V_3\to 0,
\]we have $\la^{T}(V_2)=\la^{T}(V_1)\la^{T}(V_3)$. See \cite{Hir95}.  
\vskip2mm
Introduce $K$-theoretic Chern roots $\bm \Gm_1=(\Gm_{1,1},\dots,\Gm_{1,k})$ of the tautological bundle $E_1$, $\bm \Gm_2=(\Gm_{2,1},\dots,\Gm_{2,n-k})$ of the quotient bundle $E_2$. Set $\bm \Gm=(\bm\Gm_1,\bm\Gm_2)$, and denote by $\C[\bm\Gm^{\pm 1}]^{S_k\times S_{n-k}}$ the ring of Laurent polynomials symmetric in $\bm\Gm_1$ and $\bm\Gm_2$, separately. From the exact sequence \eqref{Euler}, we obtain the following presentation of the ${T}$-equivariant $K$-theory algebra of $G(k,n)$:
\beq\label{kthpres}
K^T_0(G(k,n))_\C\cong \C[\bm\Gm^{\pm 1}]^{S_k\times S_{n-k}}\otimes \C[\bm Z^{\pm 1}]\Big/\Big\langle\la^{T}(E_1)\la^{T}(E_2)=\la^{T}(\underline{\C^n})\Big\rangle.
\eeq
Here we mean that the ideal is generated by the relations obtained by equating all the powers $u^i$, with $i=0,\dots,n$, in the identity $\la^{T}(E_1)\la^{T}(E_2)=\la^{T}(\underline{\C^n})$. That is 
\[\sum_{j=0}^he_j(\bm\Gm_1)e_{h-j}(\bm \Gm_2)=e_h(\bm Z),\quad h=1,\dots,n.
\]
Under the presentation \eqref{kthpres}, the duality involution 
\beq\label{dinv1} (-)^*\colon K^T_0(G(k,n))_\C\to K^T_0(G(k,n))_\C,\,\,[V]\mapsto [V^*],
\eeq is given by
\beq\label{dinv2} f(\bm\Gm;\bm Z)^*= f(\bm\Gm^{-1};\bm Z^{-1}).
\eeq
The pushforward morphism $\pi_*^{T}\colon K^T_0(G(k,n))\to K^T_0({\rm Spec}\,\C)\cong R({T})=\Z[\bm Z^{\pm 1}]$ is given by
\[\pi_*^{T}([V])=\sum(-1)^i[H^i(G(k,n),V)],
\]where $[H^i(G(k,n),V)]$ denotes the $R({T})$-class of the cohomology space $H^i(G(k,n),V)$ seen as a representation of ${T}$. The morphism $\pi_*^{T}$ defines the ${T}$-equivariant analog of the Euler-Poincar\'e-Grothendieck characteristic of (the isomoprhism class) of an object $V$, and it will also be denoted by $\chi^{T}(V)$.
\vskip2mm
In terms of the presentation \eqref{kthpres}, we obtain the following expression for the ${T}$-equivariant Euler-Poincar\'e-Grothendieck characteristic:
\begin{align}\nonumber\chi^{T}(f(\bm \Gm;\bm Z))&=\frac{(-1)^k}{k!}\underset{\bm \Gm_1\neq 0,\infty}{\rm Res}\left(\frac{f(\bm\Gm;\bm Z)\prod_{i\neq j}(1-\Gm_{1,i}/\Gm_{1,j})}{\prod_{i=1}^k\prod_{j=1}^n(1-\Gm_{1,i}/Z_j)}\frac{{\rm d}\Gm_{1,1}\wedge\dots\wedge{\rm d}\Gm_{1,k}}{\Gm_{1,1}\dots\Gm_{1,k}}\right)\\
\nonumber&=\frac{(-1)^k}{k!}\sum_{\substack{i_1,\dots,i_k=1\\ i_a\neq i_b}}^n\underset{\Gm_{1,1}=Z_{i_1}}{\rm Res}\dots\underset{\Gm_{1,k}=Z_{i_k}}{\rm Res}\left(\frac{f(\bm\Gm;\bm Z)\prod_{i\neq j}(1-\Gm_{1,i}/\Gm_{1,j})}{\prod_{i=1}^k\prod_{j=1}^n(1-\Gm_{1,i}/Z_j)}\frac{{\rm d}\Gm_{1,1}\wedge\dots\wedge{\rm d}\Gm_{1,k}}{\Gm_{1,1}\dots\Gm_{1,k}}\right)\\
\label{defTR}&=\sum_{I\in\mc I_{k,n}}\frac{f(\bm Z_{\si_I};\bm Z)}{\Tilde R_{k,n}(\bm Z_{\si_I})},\quad \Tilde R_{k,n}(\bm Z)=\prod_{i=1}^k\prod_{j=k+1}^n\left(1-\frac{Z_i}{Z_j}\right),\quad Z_{\si_I}=(Z_{\si_I(1)},\dots, Z_{\si_I(n)}).
\end{align}
The formula is a $K$-theoretical analog of Atiyah--Bott localization formula \cite{AB84}, and it has appeared in the literature much earlier, see \cite{AB64,BFQ79}\cite[Cor.\,6.12]{Gro77}\cite[\S4.7]{Nie74}. The computation above is a direct application of the ``Holomorphic Lefschetz Theorem'' of \cite[(4.6)]{AS68}.

By using the duality involution \eqref{dinv1}, \eqref{dinv2}, we can define the non-symmetric {\it Euler--Poincar\'e--Grothendieck pairing} (or also {\it Mukai pairing}):
\beq\label{epgm}
\chi^{T}(E,F):=\chi^{T}(E^*\otimes F).
\eeq  This pairing naturally extend to the complexified algebra $K_0^{T}(G(k,n))_{\mathbb C}$, and it is $R({T})_{\mathbb C}$-sesquilinear w.r.t.\,\,the duality involution naturally defined on $R({T})_{\mathbb C}$:
\beq\label{dualrap}(-)^*\colon R({T})_{\mathbb C}\to R({T})_{\mathbb C},\quad [V]\mapsto [V^*].
\eeq
That is,
 $\chi^{T}(\rho_1E_1,\rho_2 E_2)=\rho_1^*\rho_2\,\chi^{T}(E_1,E_2)$ for $E_1,E_2\in K_0^{T}(G(k,n))_\C$ and $\rho_1,\rho_2\in R({T})_\C$.
\vskip2mm
Denote by $\mc O_{\bm\la}$ (resp.\,$\Tilde{\mc O}_{\bm\la}$) the structural sheaf of the Schubert variety $\Om_{\bm\la}$ (resp.\,$\Tilde{\Om}_{\bm\la}$), with $\bm\la\subseteq k\times (n-k)$. Since $\mc O_{\bm\la}$ define ${T}$-equivariant sheaves on $G(k,n)$, they define $K$-theoretic classes in $K_0^{T}(G(k,n))$. 
\begin{thm}
The Grothendieck group $K_0^{T}(G(k,n))_\C$ is a free $\C[\bm Z^{\pm 1}]$-module of rank $\binom{n}{k}$. The Schubert structure sheaves $[\mc O_{\bm\la}]$ , with $\bm\la\subseteq k\times (n-k)$, define a $\C[\bm Z^{\pm 1}]$-basis of $K_0^{T}(G(k,n))_\C$.\qed
\end{thm}
For more details and proofs, see for example \cite{LS82b,KK90,AGM11}.
\vskip2mm
In terms of the presentation \eqref{kthpres}, the Schubert structures sheaves $[\mc O_{\bm\la}]$ can be represented by suitable polynomials, called\footnote{Grothendieck polynomials were introduced by A.\,Lascoux and M.P.\,Sch\"utzenberger to represent Schubert structures sheaves in (equivariant) $K$-theory of partial flag varieties \cite{LS82b}. They are inhomogeneous polynomials, generalizing Schubert polynomials (recovered by taking the lowest degree monomials), and they are indexed by permutations. In this paper, we consider only the case of Grassmannian permutations $\si$. We will directly label them by the corresponding partition $\bm\la=\bm\la^\si$. The double Grothendieck polynomials labelled by Grassmannian permutations coincide with the factorial Grothendieck polynomials as defined by P.J.\,McNamara in \cite{McN06} using set-valued tableaux. The determinantal formula is due to T.\,Ikeda and H.\,Naruse \cite{IN13}.} {\it double Grothendieck polynomials} \cite{LS82b,Las12}.  

Consider formal variables $\bm x=(x_1,\dots, x_k)$ and $\bm y=(y_1,\dots, y_n)$. For any $i=1,\dots, n-1$, define the Demazure divided difference operator $D_i^{\bm y}$, acting on functions $f(\bm y)$, by
\[D_i^{\bm y}f(\bm y)=\frac{y_if(\bm y)-y_{i+1}f(y_1,\dots,y_{i+1},y_i,\dots,y_n)}{y_i-y_{i+1}}
\]
For the bigger partition $\bm\la^{\si_0}=(n-k)^k$, the Grothendieck polynomial $\frak G_{\bm\la^{\si_0}}(\bm x;\bm y)$ is defined as
\beq
\frak G_{\bm\la^{\si_0}}(\bm x;\bm y):=\prod_{i=1}^k\prod_{j=1}^{n-k}\left(1-\frac{x_i}{y_j}\right).
\eeq
If $\bm\la,\bm\mu\subseteq k\times (n-k)$ are such that the set $I^{\si^{\bm\mu}}$ can be obtained from the set $I^{\si^{\bm\la}}$ by replacing $i+1$ with $i$, then we set
\[\frak G_{\bm\mu}(\bm x;\bm y)=D_i^{\bm y}\frak G_{\bm \la}(\bm x;\bm y).
\]Notice that $\frak G_{\bm \la}(\bm x;\bm y)$ are symmetric in $\bm x$ but not in $\bm y$. Alternatively, the double Grothendieck polynomial $\frak G_{\bm\la}(\bm x;\bm y)$ can be defined by a determinantal formula: 
given $\bm\la\subseteq k\times (n-k)$, we have
\beq
\frak G_{\bm\la}(\bm x;\bm y)=\frac{\det\left((x_i|\bm y)^{\la_j+k-j}x_i^{j-1}\right)_{1\leq i<j\leq k}}{\prod_{1\leq i<j\leq k}(x_j-x_i)},\quad (x_i|\bm y)^p:=\prod_{j=1}^p\left(1-\frac{x_i}{y_j}\right).
\eeq

\begin{prop}
We have 
\[\pushQED{\qed}
[\mc O_{\bm\la}]=\frak G_{\bm\la}(\bm\Gm_1;\bm Z)\quad\text{ and }\quad[\Tilde{\mc O_{\bm\la}}]:=\prod_{i=1}^k\Gm_{1,i}\prod_{j\in\bm\la}Z_j^{-1}\frak G_{\bm\la^\vee}(\bm\Gm;\bm Z_{\si_0}).\qedhere
\popQED
\]
\end{prop}
See \cite{LS82b,AGM11,Oet21}.

\subsection{Equivariant characteristic classes} 
In this sections we are going to define several characteristic classes of ${T}$-equivariant vector bundles on $G(k,n)$. These characteristic classes will take values in the extension $H^{(\kappa)}_T(G(k,n))$, defined in Section \ref{secext}, for a given $\kappa\in\C^*$.
\vskip2mm
Let $V$ be a rank $r$, ${T}$-equivariant vector bundle on $G(k,n)$, with equivariant Chern roots $\xi_1,\dots,\xi_r$. 
The {\it $\kappa$-graded equivariant Chern character} ${\rm Ch}^{(\kappa)}_T(V)$ 
is defined as the sum
\beq\label{kchch}
{\rm Ch}^{(\kappa)}_T(V):=\sum_{i=1}^r\exp\left(-\frac{2\pi\sqrt{-1}}{\kappa}\xi_i\right)\in H^{(\kappa)}_T(G(k,n)).
\eeq
For $\kappa=-2\pi\sqrt{-1}$ we obtain the standard equivariant Chern character. For $\kappa=-1$, we obtain the graded Chern character used in \cite{CDG1,CV21}. In the sequel, we will use the following notation:
\beq\label{acute1}
\acute{\bm z}:=\left(\acute{z}_1,\dots,\acute{z}_n\right),\quad \acute{z}_i=\exp\left(-\frac{2\pi\sqrt{-1}}{\kappa}z_i\right),\quad i=1,\dots,n,
\eeq
\beq\label{acute2}
\acute{\bm\gm}=\left(\acute{\bm\gm}_1;\acute{\bm\gm}_2\right)=\left(\acute{\gm}_{1,1},\dots,\acute{\gm}_{1,k};\acute{\gm}_{2,1},\dots,\acute{\gm}_{2,n-k}\right),\quad \acute{\gm}_{a,b}:=\exp\left(-\frac{2\pi\sqrt{-1}}{\kappa}\gm_{a,b}\right).
\eeq In terms of the presentations \eqref{presentation} and \eqref{kthpres}, the $\kappa$-graded equivariant Chern character acts as follows
\[f(\bm\Gm;\bm Z)\mapsto {\rm Ch}^{(\kappa)}_T(f)=f(\acute{\bm\gm};\acute{\bm z}).
\]

Let us identify the equivariant $K$-theory and equivariant cohomology of a point via the isomorphism
\beq\label{identkhpt}
i\colon K_0^{T}({\rm pt})_\C\to H^{(\kappa)}_T({\rm pt})\cong \mc O_{L_\kappa},\quad Z_i\mapsto \acute{z}_i,
\quad i=1,\dots,n.
\eeq
The map $i$ allows us to equip the space $H^{(\kappa)}_T(G(k,n))$ with a $K_0^{T}({\rm pt})_\C$-algebra structure. 
\begin{lem} The map ${\rm Ch}^{(\kappa)}_T\colon K^T_0(G(k,n))_\C\to H^{(\kappa)}_T(G(k,n))$ is $K_0^{T}({\rm pt})_\C$-linear with respect the structure induced by $i$. That is, for $V\in K_0^{T}(G(k,n))_\C$ and $Q\in K_0^{T}({\rm pt})_\C\cong \C[\bm Z^{\pm1}]$, we have
\beq
{\rm Ch}^{(\kappa)}_T(Q(\bm Z)V)=Q(\acute{\bm z}){\rm Ch}^{(\kappa)}_T(V). 
\eeq
\end{lem}
\proof
If $\xi_1,\dots,\xi_r$ are the equivariant Chern roots of $V$, then $\xi_1+\sum_{j=1}^n\al_jz_j,\dots,\xi_r+\sum_{j=1}^n\al_jz_j$ are the equivariant Chern roots of $Q(\bm Z)V$ with $Q(\bm Z)=Z_1^{\al_1}\dots Z_n^{\al_n}$. By additivity, the result follows.
\endproof

\vskip2mm
Let $t$ be a formal variable, and $F\in\C[\![t]\!]$ of the form $F(t)=1+\sum_{k>0}F_kt^k$. We define the {\it equivariant} $\Hat F$-{\it class} of $V$ as the characteristic class 
\[\Hat F(V):=\prod_{i=1}^rF(\xi_r),\quad \xi_1,\dots,\xi_r\text{ equivariant Chern roots of }V.
\]In this paper we will consider only vector bundles $V$ and series $F(t)$ such that the corresponding class $\Hat F(V)$ is an element of $H^{(\kappa)}_T(G(k,n))$. If $V=TG(k,n)$, the $\Hat F$-class of $V$ is called the $\Hat F$-{\it class} of $G(k,n)$. We will denote it by $\Hat F_{G(k,n)}$. Since $TG(k,n)\cong E_1^*\otimes E_2$, this class is given by
\[\Hat F_{G(k,n)}=\prod_{i=1}^k\prod_{j=1}^{n-k}F(-\gm_{1,i}+\gm_{2,j}).
\]This is the only class whose restriction at the fixed point ${\rm pt}_I$, with $I\in\mc I_{k,n}$, equals $\prod_{i\in I_1}\prod_{j\in I_2}F(z_j-z_{i})$. If $F(t)$ defines a meromorphic function on $\C$ with poles in $\kappa\Z$ only, then $\Hat F_{G(k,n)}$ is an element of $H^{(\kappa)}_T(G(k,n))$.
\vskip2mm
Relevant examples used in this paper will be the characteristic classes obtained from the series
\begin{align}
\label{ser1}
F(t)&=-\frac{2\pi\sqrt{-1}}{\kappa}\frac{t}{1-\exp\left(\frac{2\pi\sqrt{-1}}{\kappa}t\right)}=1-\frac{\pi\sqrt{-1}}{\kappa}t-2\sum_{i=1}^\infty\frac{\zeta(2i)}{\kappa^{2i}}t^{2i},\\
\label{ser2}
F(t)&=\Gm\left(1-\frac{t}{\kappa}\right)=\exp\left(\gm_{\rm EM}\frac{t}{\kappa}+\sum_{n=2}^\infty\frac{\zeta(n)}{n}\frac{t^n}{\kappa^n}\right),\\
\label{ser3}
F(t)&=\exp\left(-\frac{\log\kappa}{\kappa}t\right)\Gm\left(1-\frac{t}{\kappa}\right),
\end{align}where $\Gamma(s)$ is the Euler Gamma function, $\zeta(s)$ denotes the Riemann zeta function, and $\gm_{\rm EM}=0.577216\dots$ is the Euler--Mascheroni constant.

The characteristic class obtained from \eqref{ser1} will be called {\it $\kappa$-deformed equivariant Todd class}, and denoted by ${\rm Td}_{G(k,n)}^{(\kappa)}$.
The characteristic class obtained from \eqref{ser2} will be called {\it $\kappa$-deformed equivariant Gamma class}, and denoted by $\Hat{\Gm}_{G(k,n)}^{(\kappa)}$. Consequently, the characteristic class obtained from \eqref{ser3} equals $\kappa^{-\frac{1}{\kappa}c_1^{T}(G(k,n))}\,\Hat{\Gm}_{G(k,n)}^{(\kappa)}$.

\begin{rem}For $\kappa=-2\pi\sqrt{-1}$, we obtain the standard Todd class. In \cite{CDG1,CV21} the Gamma classes with $\kappa=\pm1$ appeared.
\end{rem}

The following is a generalization of Hirzebruch--Riemann--Roch Theorem.
\begin{thm}
For any $V\in K^T_0(G(k,n))_\C$, we have
\beq\label{HRR}
\left(-\frac{\kappa}{2\pi\sqrt{-1}}\right)^{\dim G(k,n)}\int_{G(k,n)}^{\rm eq}{\rm Ch}_T^{(\kappa)}(V){\rm Td}^{(\kappa)}_{G(k,n)}\underset{i}{=}\chi^{T}(V),
\eeq
where ${=}_i$ means ``modulo the identification $i$'' defined in \eqref{identkhpt}.
\end{thm}
\proof
If $f_o(\bm\gm;\bm z)={\rm Td}^{(\kappa)}_{G(k,n)}$, we have $f_o(\bm z;\bm z)=(\frac{2\pi\sqrt{-1}}{\kappa})^{k(n-k)}R_{k,n}(\bm z)/\Tilde R_{k,n}(\acute{\bm z})$, where $R_{k,n}$ and $\Tilde R_{k,n}$ are defined in \eqref{defVR} and \eqref{defTR}, respectively. Then, for an arbitrary $f(\bm\Gm;\bm Z)$ in $K^T_0(G(k,n))_\C$, we have
\begin{multline*}\int_{G(k,n)}^{\rm eq}{\rm Ch}_T^{(\kappa)}(f){\rm Td}^{(\kappa)}_{G(k,n)}=(-1)^{k(n-k)}\sum_{I\in\mc I_{k,n}}\frac{f(\acute{\bm z}_{\si_I};\acute{\bm z}_{\si_I}) f_o(\bm z_{\si_I};\bm z_{\si_I})}{R_{k,n}(\bm z_{\si_I})}\\=\left(-\frac{2\pi\sqrt{-1}}{\kappa}\right)^{k(n-k)}\sum_{I\in\mc I_{k,n}}\frac{f(\acute{\bm z}_{\si_I};\acute{\bm z}_{\si_I})}{\Tilde R_{k,n}(\acute{\bm z}_{\si_I})}=\left(-\frac{2\pi\sqrt{-1}}{\kappa}\right)^{k(n-k)} i\left(\chi^{T}(f)\right).
\end{multline*}
\endproof

\begin{rem}
By setting $\kappa=-2\pi\ic$, and taking the non-equivariant limit, one recovers the classical Hirzebruch--Riemann--Roch Theorem.
\end{rem}

\subsection{The \textcyr{B}-morphism} Let $\kappa$ be a point of the universal cover $\Tilde{\C^*}$, so that $\log\kappa$ is a well-defined function. We define the \textcyr{B}-class of $V\in K^T_0(G(k,n))_\C$ to be the cohomology class
\beq
\text{\textcyr{B}}(V;\kappa):= \Hat{\Gamma}^{(\kappa)}_{G(k,n)}\,\exp\left[-\frac{\log\kappa}{\kappa}c_1^{T}(G(k,n))\right]{\rm Ch}^{(\kappa)}_T(V).
\eeq
This defines a deformation of the Chern character, that we will call \textcyr{B}-morphism
\beq
\text{\textcyr{B}}(-;\kappa)\colon K^T_0(G(k,n))_\C\to H^{(\kappa)}_T(G(k,n)).
\eeq
By setting $\kappa=e^{\pi\sqrt{-1}}$, the class above specializes to the \textcyr{B}-class of \cite{CV21}.

\begin{thm}\label{HRR2}
For any $V_1,V_2\in K^T_0(G(k,n))_\C$, we have
\beq
\left(-\frac{\kappa}{2\pi\sqrt{-1}}\right)^{\dim G(k,n)}\int^{\rm eq}_{G(k,n)}\textnormal{\textcyr{B}}(V_1;e^{-\pi\sqrt{-1}}\kappa)\textnormal{\textcyr{B}}(V_2;\kappa)\underset{i}{=}\chi^{T}(V_1,V_2).
\eeq
\end{thm}
\proof
From the identity $$\Gm(1+t)\Gm(1-t)=\frac{2\pi\ic t}{1-e^{-2\pi\ic t}}e^{-\pi\ic t},$$ we deduce that
\[\textnormal{\textcyr{B}}\left(V_1;e^{-\pi\sqrt{-1}}\kappa\right)\textnormal{\textcyr{B}}(V_2;\kappa)=
{\rm Td}^{(\kappa)}_{G(k,n)}{\rm Ch}_T^{(\kappa)}(V_1^*\otimes V_2),
\]
and the result follows from \eqref{HRR}.
\endproof

\subsection{Exceptional collections and their mutations} Given $F_1,F_2\in{\rm Ob}\left({\mathsf D}^b_T(G(k,n))\right)$, we define
\[\Hom^\bullet_T(F_1,F_2):=R\pi_*^{T}(F_1^*\otm F_2)\in{\rm Ob}\left({\mathsf D}^b({\on{Rep}_\C}({T}))\right),
\]where $F_1^*:=R\HOM(F_1,\mathcal O_X)$ is the ordinary dual sheaf of $F_1$.
\begin{defn}\label{exccol}
An object $E\in{\rm Ob}({\mathsf D}^b_T(G(k,n)))$ is called an \emph{exceptional object} if and only if
\[\Hom^\bullet_T(E,E)\cong \mathbb C_T,
\]
where $\mathbb C_T$ denotes the object of ${\mathsf D}^b({\on{Rep}_\C}({T}))$ given by the trivial complex one dimensional representation of ${T}$, concentrated in degree zero.

An ordered collection $(E_1,\dots, E_m)$ is said to be an \emph{exceptional collection} if and only if
\begin{itemize}
\item all objects $E_i$'s are exceptional objects,
\item and $\Hom^\bullet_T(E_j, E_i)=0$ for $j>i$.
\end{itemize}
\end{defn}

Given $E\in {\rm Ob}({\mathsf D}^b_T(G(k,n)))$ and $V^\bullet\in{\rm Ob}({\mathsf D}^b({\on{Rep}_\C}({T})))$, the tensor product $E\otimes V^\bullet$ is defined as the object of ${\mathsf D}^b_T(G(k,n))$ defined as $\bigoplus_{i}E[-i]\otimes V^i.$
This extends the obvious operation of tensor product between objects of $Coh_T(G(k,n))$ and ${\on{Rep}_\C}({T})$. 

If $\mathcal A_1,\dots,\mathcal A_m$ are subcategories of ${\mathsf D}^b_T(G(k,n))$, we denote by $\left\langle \mathcal A_1,\dots,\mathcal A_m\right\rangle$ the smallest full triangulated subcategory of ${\mathsf D}^b_T(G(k,n))$ containing $\mathcal A_1,\dots, \mathcal A_m$.

\begin{defn}
Let $\frak E:=(E_1,\dots, E_{\binom{n}{k}})$ be an exceptional collection in ${\mathsf D}^b_T(G(k,n))$. We  say that $\frak E$ is ${T}$-\emph{full} if
\beq\label{Gfull}{\mathsf D}^b_T(G(k,n))=\left\langle E_1\otimes{\mathsf D}^b({\on{Rep}_\C}({T})),\dots,  E_n\otimes{\mathsf D}^b({\on{Rep}_\C}({T}))\right\rangle.
\eeq
\end{defn}

\begin{thm}\cite{Ela09,Pol11}\label{Tfullcoll}
Let $(E_1,\dots, E_{\binom{n}{k}})$ be a full exceptional collection in ${\mathsf D}^b(G(k,n))$. Each object $E_i$, with $i=1,\dots,{\binom{n}{k}}$, admits a ${T}$-equivariant structure $\mc E_i$, unique up to tensoring by a character of ${T}$. The resulting collection $(\mc E_1,\dots,\mc E_{\binom{n}{k}})$ is a ${T}$-full exceptional collection in ${\mathsf D}^b_T(G(k,n))$.\qed
\end{thm}

\begin{defn}[Mutations of objects]
Let $E\in{\rm Ob}({\mathsf D}^b_T(G(k,n)))$ be an exceptional object. For any $F\in{\rm Ob}({\mathsf D}^b_T(G(k,n)))$ we  define two new objects
\[\mathbb L_EF,\quad \mathbb R_EF,
\]called the \emph{left} and \emph{right mutations} of $F$ with respect to $E$. These two objects are defined (up to unique isomoprhism) through the distinguished triangles
\beq\label{left}
\xymatrix{\mathbb L_EF[-1]\ar[r]&\Hom^\bullet_G(E,F)\otimes E\ar[r]^{\quad\quad\quad j^*}&F\ar[r]&\mathbb L_EF,}
\eeq
\beq\label{right}
\xymatrix{\mathbb R_EF\ar[r]&F\ar[r]^{j_*\quad\quad\quad}&\Hom^\bullet_G(F,E)^*\otimes E\ar[r]&\mathbb R_EF[1],}
\eeq
where $j^*,j_*$ denote the natural evaluation and coevaluation morphisms.
\end{defn}

\begin{lem}
Let $E\in{\rm Ob}({\mathsf D}^b_T(G(k,n)))$ be an exceptional object. We have
\[\Hom^\bullet_G(E,\mathbb L_EF)=0,\quad \Hom_G^\bullet(\mathbb R_EF,E)=0,
\]for all objects $F\in{\rm Ob}({\mathsf D}^b_T(G(k,n)))$.
\end{lem}

\proof
Apply the functor $\Hom^\bullet_G(E,-)$ (resp. $\Hom^\bullet_G(-,E)$) to the distinguished triangle \eqref{left} (resp. \eqref{right}), and use the exceptionality of $E$.
\endproof

\begin{defn}
Let $\frak E:=(E_1,\dots,E_{\binom{n}{k}})$ be an exceptional collection in ${\mathsf D}^b_T(G(k,n))$. For any integer $i$, with $0< i<{\binom{n}{k}}$, we define two new collections
\begin{align*}\mathbb L_i(\frak E):&=(E_1,\dots, \mathbb L_{E_i}E_{i+1}, E_i,\dots, E_{\binom{n}{k}}),\\
\mathbb R_i(\frak E):&=(E_1,\dots, E_{i+1}, \mathbb R_{E_{i+1}}E_i, \dots, E_{\binom{n}{k}}).
\end{align*}
\end{defn}

\begin{prop}
For any $i$, with $0<i<{\binom{n}{k}}$, the
 collections $\mathbb L_i(\frak E), \mathbb R_i(\frak E)$ are exceptional.
  Moreover, the mutation operators $\mathbb L_i,\mathbb R_i$ satisfy the following identities:
\beq\label{br1}
\mathbb L_i\mathbb R_i=\mathbb R_i\mathbb L_i={\rm Id},
\eeq
\beq\label{br2}
\mathbb R_i\mathbb R_j=\mathbb R_j\mathbb R_i,\quad \text{if}\ \
|i-j|>1,\quad \mathbb R_{i+1}\mathbb R_i\mathbb R_{i+1}=\mathbb R_i\mathbb R_{i+1}\mathbb R_i.
\eeq
\end{prop}

\proof
The same as in the non-equivariant case, see \cite{GK04}, \cite[Section 3.3]{CDG1}.
\endproof

Denote by $\tau_1,\dots, \tau_{{\binom{n}{k}}-1}$ the generators of the braid group $\mathcal B_{\binom{n}{k}}$, satisfying the relations 
\[\tau_i\tau_{i+1}\tau_i=\tau_{i+1}\tau_i\tau_{i+1},\quad \tau_i\tau_j=\tau_j\tau_i,\quad
\text{if}\ \ |i-j|>1.
\]
Following \cite{CV21}, we define the left action of $\mathcal B_{\binom{n}{k}}$ on the set of exceptional collections of length ${\binom{n}{k}}$ by identifying the action of  $\tau _i$ with $\mathbb R_{{\binom{n}{k}}-i}$, for $i=1,\dots,{\binom{n}{k}}-1$.
 See identities \eqref{br1}-\eqref{br2}.

\subsection{Exceptional bases}In this Section we focus on the $K$-theoretical counterpart of the notion of exceptional collections introduced in Definition \ref{exccol} and of the action of the braid group on them.

\begin{defn}\label{excbas}An element $e\in K_0^{T}(G(k,n))$ is \emph{exceptional} if 
 $$\chi^{T}(e,e)=\mathbb C_T.$$ A basis $\varepsilon:=(e_1,\dots,e_{\binom{n}{k}})$ of $K_0^{T}(G(k,n))$ as an
  $R(G)$-module, is  \emph{exceptional} if 
\beq
\chi^{T}(e_i,e_i)=\mathbb C_T,\quad \chi^{T}(e_j,e_i)=0,\quad \text{for } j>i.
\eeq 
\end{defn}

The following result is a $K$-theoretical analogue of Theorem \ref{Tfullcoll}.

\begin{thm}[{\cite[ Lemma 2.1]{Pol11}}]
\label{pol2}
Let $(E_1,\dots, E_{\binom{n}{k}})$ be a full exceptional collection in ${\mathsf D}^b(G(k,n))$. If $\mc E_i$ is a ${T}$-equivariant structure on $E_i$, 
then the classes $([\mc E_1],\dots,[\mc E_{\binom{n}{k}}])$ form an exceptional basis of $K_0^{T}(G(k,n))$ as an
$R({T})$-module.\qed
\end{thm}

\begin{prop}
Let $E\in{\rm Ob}({\mathsf D}^b_T(G(k,n)))$ be an exceptional object. For any $F\in {\rm Ob}({\mathsf D}^b_T(G(k,n)))$ we have
\beq
[\mathbb L_EF]=[F]-\chi^{T}(E,F)\cdot [E],\quad [\mathbb R_EF]=[F]-\chi^{T}(F,E)^*\cdot [E].
\eeq
\end{prop}

\proof See \cite[Sec.\,2.5, Prop.\,2.18]{CV21}.
\endproof

\begin{defn}
Let $e\in K_0^{T}(G(k,n))$ be an exceptional element. Given $f\in K_0^{T}(G(k,n))$, we define its \emph{left} and \emph{right mutations} wrt $e$ as the elements
\beq
\mathbb L_{e}f:=f-\chi^{T}(e,f)\cdot e,\quad \mathbb R_{e}f:=f-\chi^{T}(f,e)^*\cdot e.
\eeq
\end{defn}

\begin{lem}\label{lem}
Let $e\in K_0^{T}(G(k,n))$ be an exceptional element. We have
\[\chi^{T}(e,\mathbb L_ef)=0,\quad \chi^{T}(\mathbb R_ef,e)=0,
\]for any $f\in K_0^{T}(G(k,n))$.
\qed
\end{lem}

\begin{defn}
Let $K^G_0(G(k,n))$ be a free $R(G)$-module of finite rank and  $\varepsilon:=(e_1,\dots,e_{\binom{n}{k}})$ an exceptional basis of
 $K^G_0(G(k,n))$. For any $0<i<n$ define the two new exceptional bases 
\beq
\mathbb L_i\varepsilon:=(e_1,\dots,e_{i-1},\mathbb L_{e_i}e_{i+1},e_i,e_{i+2},\dots, e_{\binom{n}{k}}),
\eeq
\beq
\mathbb R_i\varepsilon:=(e_1,\dots,e_{i-1}, e_{i+1}, \mathbb R_{e_{i+1}}e_i,e_{i+2},\dots, e_{\binom{n}{k}}).
\eeq

This construction defines the action of the braid group $\mathcal B_{\binom{n}{k}}$ on the set of exceptional bases of $K^G_0(G(k,n))$,
in which the action of a generator $\tau_i$,  $i=1,\dots, {\binom{n}{k}}-1$, is identified with the action of the mutation
 $\mathbb R_{{\binom{n}{k}}-i}$.
\end{defn}

\begin{rem}\label{remsignedbraid}
If \( e \in K^T_0(G(k,n)) \) is an exceptional object, then \( -e \) is also exceptional. Consequently, the group \( \mathbb{Z}^{\binom{n}{k}} \) acts on the set of exceptional bases. This action, together with the braid group action, can be combined to provide an action of the \textit{signed braid group}. This is the semidirect product \( \mathbb{Z}^{\binom{n}{k}} \rtimes \mathcal{B}_{\binom{n}{k}} \), where the action \( \mathcal{B}_{\binom{n}{k}} \to \mathrm{Aut}(\mathbb{Z}^{\binom{n}{k}}) \) is given by the permutations of the entries of the tuples \( (\pm1, \dots, \pm1) \).
\end{rem}

\subsection{Schur functors and Kapranov collection} Let $E\to G(k,n)$ be a ${T}$-equivariant vector bundle. Given a partition $\bm\la\subseteq k\times (n-k)$, the {\it Schur vector bundle} $\Si_{\bm \la}(E)\to G(k,n)$ is the ${T}$-equivariant vector bundle satisfying the following property: if ${\rm pt}_I\in G(k,n)$ is a ${T}$-fixed point, with $I\in \mc I_{k,n}$, then the fiber $\Si_{\bm \la}(E)_{{\rm pt}_I}$ is the ${T}$-module with character the Schur polynomial $s_{\bm\la}\in K^T_0({\rm pt}_I)\cong \Z[\bm Z_I^{\pm 1}]$, $\bm Z_I=(Z_i)_{i\in I_1}$. Main references are \cite{Kap84,Kap88,Wey03}.

\begin{example}
If $\bm\la=(1)^k$, then $\Si_{\bm\la}(E)=\wedge^kE$. If $\bm\la=(k)$, then $\Si_{\bm \la}(E)={\rm Sym}^k E$.
\end{example}
\vskip2mm
Recall the short exact sequence \eqref{Euler} on $G(k,n)$. The next result is due to M.M.\,Kapranov, who extended Beilinson's method for constructing full exceptional collections on $\mathbb P^n$, see \cite{Bei78}.

\begin{thm}\cite{Kap84,Kap88}
The collection $\left(\Si_{\bm\la}(E_1^*)\right)_{\bm\la\subseteq k\times (n-k)}$, in any linear order refining the partial inclusion order on the Young diagrams, is a full exceptional collection in ${\mathsf D}^b(G(k,n))$.\qed
\end{thm}

From Theorem \ref{Tfullcoll}, we deduce the following result.
\begin{cor}\label{kapfullcoll}
The collection $\left(\Si_{\bm\la}(E_1^*)\right)_{\bm\la\subseteq k\times (n-k)}$ is a ${T}$-full exceptional collection in ${\mathsf D}^b_T(G(k,n))$.\qed
\end{cor}

\begin{prop}\cite{Kap84,Kap88}
We have the following isomorphisms of ${T}$-modules:
\begin{enumerate}
\item For any $i\geq 0$, $H^i(G(k,n),\Si_{\bm \la}(E_1))=0$.
\item 
\begin{empheq}[left ={H^i(G(k,n),\Si_{\bm \la}(E_1^*))}{=}\empheqlbrace]{align*}
&\Si_{\bm\la}(\underline{\C^n})^*\quad\text{if $i=0$,}\\
&0,\quad\text{if $i>0$.}
\end{empheq}
\item \begin{empheq}[left ={H^i(G(k,n),\Si_{\bm \la}(E_2))}{=}\empheqlbrace]{align*}
&\Si_{\bm\la}(\underline{\C^n})\quad\text{if $i=0$,}\\
&0,\quad\text{if $i>0$.}
\end{empheq}
\end{enumerate}\qed
\end{prop}

\begin{cor}\label{corchikapr}
Given $\bm\la_1,\bm\la_2\subseteq k\times (n-k)$. The complex $\Hom^\bullet_T(\Si_{\bm\la_1}(E_1^*),\Si_{\bm\la_2}(E_1^*))$ is concentrated in degree 0, and it equals 
\beq
H^0(G(k,n),\Si_{\la_1}(E_1)\otimes \Si_{\la_2}(E_1^*))\cong \bigoplus_{\bm\gm} \Si_{\bm\gm}(\underline{\C^n})^*,
\eeq
where $\bm\gm\subseteq k\times (n-k)$ runs over all possible summands in the decomposition of $\Si_{\la_1}(E_1)\otimes \Si_{\la_2}(E_1^*)$ into irreducibles.\qed
\end{cor}
\begin{rem}
The decomposition into irreducibles of of $\Si_{\la_1}(E)\otimes \Si_{\la_2}(E)$, where $E$ is a ${T}$-equivariant vector bundle, is dictated by the Littlewood--Richardson rule, see \cite[pg.\,121]{Ful97}.
\end{rem}

\section{Equivariant quantum cohomology, qDE and qKZ equations for Grassmannians}\label{sec5}
\subsection{Equivariant Gromov--Witten invariants} Let $d\in H_2(G(k,n),\Z)$, and $g,m\geq 0$. Denote by $\overline{\mc M}_{g,m}(G(k,n),d)$ the moduli stack of genus $g$ stable maps with target $G(k,n)$, degree $d$ and $m$ marked points. If $d$ is the {\it Mori cone} ${\rm Eff}(G(k,n))$ (this is the semigroup in $H_2(G(k,n),\Z)$ generated by those classes representable by holomoprhic curves in $G(k,n)$), and in the regime $2g+m>2$, the stack $\overline{\mc M}_{g,m}(G(k,n),d)$ is non-empty. 

The ${T}$-action on $G(k,n)$ induces an action on $\overline{\mc M}_{g,m}(G(k,n),d)$. Given $m$ cohomological classes $\al_1,\dots,\al_m\in H^\bullet_T(G(k,n),\C)$, and integers $d_1,\dots, d_m\in\mathbb Z_{\geq 0}$, we define the \emph{genus $g$, 
   degree $d$, ${T}$-equivariant descendant Gromov-Witten invariants} 
   of $G(k,n)$ to be the polynomials
  \beq
\label{eqgw}
\langle\tau_{d_1}(\al_1),\dots,\tau_{d_m}(\al_m)\rangle^{G(k,n), {T}}_{g,m,d}:=\left(\int^{\rm eq}_{[\overline{\mathcal M}_{g,m}(G(k,n),d)]^{\rm vir}_{{T}}}
\prod_{j=1}^m\psi^{d_j}_{j}{\rm ev}^*_j(\al_j)\right)\in 
H^\bullet_{{T}}({\rm pt},\mathbb C),
\eeq 
where 
\begin{itemize}
\item the map ${\rm ev}_j\colon \overline{\mathcal M}_{g,m}(G(k,n),d)\to G(k,n)$ is the evaluation at the $j$-th marked point, which is ${T}$-equivariant,
\item the classes $\psi_j\in A^{{T}}_{1}(\overline{\mathcal M}_{g,m}(G(k,n),d))$ denote any equivariant lift of the first Chern classes of the universal cotangent line bundles $\mathcal L_j$ on $\overline{\mathcal M}_{g,m}(G(k,n),d)$,
\item the equivariant integral is over the equivariant virtual fundamental class $[\overline{\mathcal M}_{g,m}(G(k,n),d)]^{\rm vir}_{{T}}$, element of the Chow group $A^T_{D_{\rm vir}}(\overline{\mathcal M}_{g,m}(G(k,n),d),\Q)$ with
\begin{center}$D_{\rm vir}=
\int_dc_1(TG(k,n))+(\dim_\C G(k,n)-3)(1-g)+m.$\end{center}
\end{itemize}
When all $d_i$'s equal zero, we call the polynomial \eqref{eqgw} a {\it primary} Gromov--Witten invariant.
\begin{rem}
The existence of equivariant virtual fundamental class $[\overline{\mathcal M}_{g,m}(G(k,n),d)]^{\rm vir}_{{T}}$ is ensured by the \emph{properness} of $\overline{\mathcal M}_{g,m}(G(k,n),d)$. From this property, it also follows that equivariant Gromov--Witten invariants are polynomials in $\bm z$. By localization, one can reduce the computation of equivariant Gromov--Witten invariants to the intersection theory of the fixed torus substack $\overline{\mc M}_{g,m}(G(k,n),d)^{T}$ of $\overline{\mc M}_{g,m}(G(k,n),d)$. See \cite[Section 3]{liush} and references therein.
\end{rem}

\subsection{Equivariant Gromov--Witten potential} Fix a basis $(\eps_1,\dots, \eps_{\binom{n}{k}})$ of $H^\bullet_T(G(k,n),\C)$ as $H^\bullet_T({\rm pt},\C)$-module. Let $\bm t=(t^1,\dots,t^{\binom{n}{k}})$ be the dual coordinates.
Introduce a formal variable $q$, and let $\La=\C[\![q^d;\,d\in{\rm Eff}(G(k,n))]\!]$ be the formal semigroup ring of the Mori cone ${\rm Eff}={\rm Eff}(G(k,n))$, 
with well-defined\footnote{The Mori cone consists of lattice points in a strongly convex cone. This implies that, given $d\in{\rm Eff}(G(k,n))$, the equation $d=d_1+d_2$ has a finite number of solutions $(d_1,d_2)$. So the product \eqref{moriprod} is well-defined.} product 
\beq\label{moriprod}\sum_{d}a_dq^d\cdot \sum_db_dq^d=\sum_{d}\Big(\sum_{d=d_1+d_2}a_{d_1}b_{d_2}\Big)q^d.
\eeq

The {\it equivariant Gromov--Witten potential} $F^{G(k,n),{T}}_0\in H^\bullet_T({\rm pt},\La)[\![\bm t]\!]$ is the generating function of genus 0 primary equivariant Gromov--Witten invariants, defined by
\[F^{G(k,n),{T}}_0(\bm t)=\sum_{m=0}^\infty\sum_{d\in{\rm Eff}}\sum_{i_1,\dots,i_m=1}^{\binom{n}{k}}\frac{t^{i_1}\dots t^{i_m}}{m!}\langle \eps_{i_1},\dots,\eps_{i_m}\rangle^{G(k,n),{T}}_{0,m,d}q^d.
\]
For each $a,b=1,\dots,\binom{n}{k}$, set $\eta_{ab}=\eta(\eps_a,\eps_b)$ entries of the Gram matrix $\eta=\left(\eta_{ab}\right)_{ab}$. Denote by $\eta^{ab}$ the entries of the inverse Gram matrix $\eta^{-1}$.
\begin{thm}[{\cite[Theorem 3.1]{giv1}}]
The function $F^{G(k,n),{T}}_0(\bm t)$ satisfies the WDVV equations:
\[
\pushQED{\qed}\sum_{c,e}\frac{\der^3 F^{G(k,n),{T}}_0}{\der t^a\der t^b\der t^c}\eta^{ce}\frac{\der^3 F^{G(k,n),{T}}_0}{\der t^e\der t^\ell\der t^k}=
\sum_{c,e}\frac{\der^3 F^{G(k,n),{T}}_0}{\der t^k\der t^b\der t^c}\eta^{ce}\frac{\der^3 F^{G(k,n),{T}}_0}{\der t^e\der t^\ell\der t^a}.\qedhere
\popQED\]
\end{thm}

\subsection{Equivariant quantum cohomology}
The {\it big equivariant quantum product} $*_{\bm t,q}$ is the product on $H^\bullet_T(G(k,n),\La[\![\bm t]\!])$ defined by
\beq\label{qprod}
\eps_a*_{\bm t,q}\eps_b=\sum_{c,e}\frac{\der^3 F^{G(k,n),{T}}_0(\bm t)}{\der t^a\der t^b\der t^c}\eta^{ce}\eps_e.
\eeq
The space $H^\bullet_T(G(k,n),\La[\![\bm t]\!])$ equipped with the product $*_{\bm t,q}$ and the equivariant Poincar\'e metric $\eta$ (extended by $\La[\![\bm t]\!]$-bilinearity) is a Frobenius algebra. This consists of a commutative, associative algebra with unit (the cohomology class 1) whose product is compatible with the metric:
\beq
\eta(a*_{\bm t,q}b,c)=\eta(a,b*_{\bm t,q}c),\quad a,b,c\in H^\bullet_T(G(k,n),\La[\![\bm t]\!]).
\eeq
This algebra structure on $H^\bullet_T(G(k,n),\La[\![\bm t]\!])$ is called the {\it big equivariant quantum cohomology} of $G(k,n)$. It provides an example of a {\it formal Frobenius manifold} in the sense of Yu.I.\,Manin \cite[Ch.\,III]{manin}. The quantum product \eqref{qprod} is a deformation of both the quantum and the equivariant product in cohomology of Grassmannians.
\vskip2mm
If we specialize $\bm t=0$, we obtain the {\it small equivariant quantum product} $*_{q}=*_{0,q}$. This defines a Frobenius algebra structure on $H^\bullet_T(G(k,n),\La)$. It turns out that the structural constants of $*_{q}$ are polynomials in $q$, see e.g.\,\,\cite{givkim,Mih06,Mih08}. This allows us to define a family of Frobenius algebra structures on $(H^\bullet_T(G(k,n),\C)$ parametrized by a parameter $q\in\C^*$:
\begin{align}\label{seqp}
\eps_i*_q\eps_j=&\sum_{\substack{d\in{\rm Eff}}}\sum_{a,b}\langle\eps_i,\eps_j,\eps_a\rangle_{0,3,d}^{G(k,n),{T}}\eta^{ab}\eps_b q^d\\
\nonumber
=&\eps_i\eps_j+\sum_{\substack{d\in{\rm Eff}\setminus\{0\}}}\sum_{a,b}\langle\eps_i,\eps_j,\eps_a\rangle_{0,3,d}^{G(k,n),{T}}\eta^{ab}\eps_b q^d,\qquad q^d=q^{-\int_d c_1(E_1)}=q^{\int_dc_1(E_2)}.
\end{align}
The first term $\eps_i\eps_j$ is the (classical) equivariant cup product, and $c_1(E_1),c_1(E_2)\in H^2(G(k,n),\Z)$ denotes the non-equivariant first Chern class of the bundles $E_1,E_2$, repsectively.
\subsection{Quantum multiplication and stable envelope maps} Given $a\in H^\bullet_T(G(k,n),\C)$, consider the operator $a*_{q}$ of (small) equivariant multiplication by $a$, as in \eqref{seqp}. This is a $\C[\bm z]$-module endomorphism of $H^\bullet_T(G(k,n),\C)$.
\vskip2mm
The dynamical operators $X_1(\bm z;\bm p)$ and $X_2(\bm z;\bm p)$, defined in \eqref{dynop1} and \eqref{dynop2}, are linear in $\bm z$. Their action on $(\C^2)^{\otm n}_{k,n-k}$ extends $\C[\bm z]$-linearly on $(\C^2)^{\otm n}_{k,n-k}\otm\C[\bm z]$.

\begin{thm}\cite{TV23}\label{stabint}
If $q=p_2p_1^{-1}$, the isomorphism ${\rm Stab}_{k,n}$ intertwines the dynamical operators $X_1(\bm z;\bm p),X_2(\bm z;\bm p)$ acting on $(\C^2)^{\otm n}_{k,n-k}\otm\C[\bm z]$ and the operators of quantum multiplication $c_1^{T}(E_1)*_q, $ $c_1^{T}(E_2)*_q$ acting on $H^\bullet_T(G(k,n),\C)$:
\beq\label{semqm}
{\rm Stab}_{k,n}\circ X_i(\bm z;\bm p)=c_1^{T}(E_i)*_q\circ\,{\rm Stab}_{k,n},\quad i=1,2.
\eeq
\end{thm}
\proof The small equivariant quantum multiplication by divisor classes in $H^\bullet_T(G(k,n),\C)$ (and for more general homogeneous spaces) is described in \cite[Th.\,6.4]{Mih07}. 
For $i=1$, identity \eqref{semqm} follows by a term by term comparison of formula \eqref{dynop1} with formula (6.1) in \cite{Mih07}. For $i=2$, the result then follows since $c_1^{T}(E_1)+c_1^{T}(E_2)=z_1+\dots+z_n$ in $H^\bullet_T(G(k,n),\C)$, and $X_1(\bm z;\bm p)+X_2(\bm z;\bm p)=z_1+\dots+z_n$ in $(\C^2)^{\otm n}_{k,n-k}\otm\C[\bm z]$.
\endproof

\subsection{qDE and qKZ equations}\label{secqdeqkz}Given $\bm z_o\in\C^n$, consider the algebra 
\[H^\bullet_T(G(k,n),\C)_{\bm z_o}=\C[\bm \gm]^{S_k\times S_{n-k}}\Big/\Big\langle\sum_{h=0}^\ell e_h(\bm\gm_1)e_{\ell-h}(\bm\gm_2)=e_\ell(\bm z_o),\,\ell=1,\dots,n\Big\rangle.
\]The evaluation map $\C[\bm \gm]^{S_k\times S_{n-k}}\otm \C[\bm z]\to \C[\bm \gm]^{S_k\times S_{n-k}}$, $f(\bm\gm;\bm z)\mapsto f(\bm\gm;\bm z_o)$, identifies the algebra  $H^\bullet_T(G(k,n),\C)_{\bm z_o}$ with the quotient $H^\bullet_T(G(k,n),\C)/\langle\bm z=\bm z_o\rangle$.
\vskip2mm
Given $I\in\mc I_{k,n}$, let ${\rm Stab}_{I,\bm z_o}$, ${\rm Stab}^{\rm op}_{I,\bm z_o}$ be the classes of ${\rm Stab}_{I}(\bm\gm;\bm z_o)$, ${\rm Stab}^{\rm op}_{I}(\bm\gm;\bm z_o)$ in the quotient $H^\bullet_T(G(k,n),\C)_{\bm z_o}$. Consider the maps
\[{\rm Stab}_{k,n;\bm z_o}\colon (\C^2)^{\otm n}_{k,n-k}\to H^\bullet_T(G(k,n),\C)_{\bm z_o},\quad v_I\mapsto {\rm Stab}_{I,\bm z_o}.
\]
\begin{lem}
For every $\bm z_o\in\C^n$, the map ${\rm Stab}_{k,n;\bm z_o}$ is an isomoprhism of $\C$-vector spaces. That is the classes ${\rm Stab}_{I,\bm z_o}$, with $I\in \mc I_{k,n}$, define a basis of $H^\bullet_T(G(k,n),\C)_{\bm z_o}$.\qed
\end{lem}

\begin{lem}\label{stabintzo}
If $q=p_2p_1^{-1}$, for every $\bm z_o\in\C^n$, the isomorphism ${\rm Stab}_{k,n;\bm z_o}$ intertwines the dynamical operators $X_1(\bm z_o;\bm p),X_2(\bm z_o;\bm p)$ acting on $(\C^2)^{\otm n}_{k,n-k}$ and the operators of quantum multiplication $c_1^{T}(E_1)_{\bm z_o}*_{q,\bm z_o}, $ $c_1^{T}(E_2)_{\bm z_o}*_{q,\bm z_o}$ acting on $H^\bullet_T(G(k,n),\C)_{\bm z_o}$:
\beq\label{semqm}
{\rm Stab}_{k,n;\bm z_o}\circ X_i(\bm z_o;\bm p)=c_1^{T}(E_i)_{\bm z_o}*_{q,\bm z_o}\circ\,{\rm Stab}_{k,n;\bm z_o},\quad i=1,2.
\eeq\qed
\end{lem}

Consider the space $\C^n\times \C^2$ with coordinates $(\bm z;\bm p)$. We have two trivial vector bundles 
\begin{itemize}
\item $H_{k,n}\to\C^n\times \C^2$ with fiber over a point $(\bm z_o;\bm p_o)$ given by $H^\bullet_T(G(k,n),\C)_{\bm z_o}$;
\item $U_{k,n}\to\C^n\times \C^2$ with fiber $(\C^2)^{\otm n}_{k,n-k}$.
\end{itemize}

\begin{lem}
The map ${\rm Stab}_{k,n}^\diamond\colon U_{k,n}\to H_{k,n}$, $(z_o,\bm p_o,v)\mapsto (\bm z_o;\bm p_o;{\rm Stab}_{k,n;\bm z_o}v)$ is an isomorphism of vector bundles. \qed
\end{lem}

The ({\it small}) {\it equivariant quantum differential equations} for sections of $H_{k,n}$ is the system of differential equations
\beq\label{qde.0}
\kappa p_i\frac{\der f}{\der p_i}=c_1^{T}(E_i)*_{q,\bm z}f,\quad q=\frac{p_2}{p_1},\quad i=1,2.
\eeq

By Lemma \ref{stabintzo}, the isomorphism ${\rm Stab}_{k,n}^\diamond$ identifies equations \eqref{qde.0} with the dynamical differential equations \eqref{dyn.0} for sections of $U_{k,n}$.

Furthermore, the isomorphism ${\rm Stab}_{k,n}^\diamond$ and the qKZ equations \eqref{qKZ.0} for sections of $U_{k,n}$ define the {\it qKZ difference equations in cohomology}
\beq\label{qkz.H}
f(z_1,\dots,z_a+\kappa,\dots,z_n;\bm p;\kappa)=K_a^H(\bm z;\bm p;\kappa)f(\bm z;\bm p;\kappa),\quad a=1,\dots,n,
\eeq
where for each $a=1,\dots,n$, and fixed $\bm z,\bm p$ the operator $K_a^H(\bm z;\bm p;\kappa)$ is a map of fibers
\[K_a^H(\bm z;\bm p;\kappa)\colon H^\bullet_T(G(k,n),\C)_{\bm z}\to H^\bullet_T(G(k,n),\C)_{(z_1,\dots,z_a+\kappa,\dots,z_n)},
\]
\[K_a^H(\bm z;\bm p;\kappa)={\rm Stab}_{k,n;(z_1,\dots,z_a+\kappa,\dots,z_n)}\circ K_a(\bm z;\bm p;\kappa)\circ {\rm Stab}_{k,n;\bm z}^{-1},
\]and the operator $K_1(\bm z;\bm p;\kappa)$ is given by \eqref{qkzop}.
The following result directly follows from Theorem \ref{cjs}.
\begin{thm}
The joint system of quantum differential equations \eqref{qde.0} and qKZ difference equations \eqref{qkz.H} is a compatible joint system of differential and difference equations for sections of $H_{k,n}$.\qed
\end{thm}

In \cite{CV21}, the authors studied the joint system of qDE and qKZ equations associated with projective spaces $\Pb^{n-1}$. To compare formulae of the current paper with their counterparts in \cite{CV21}, the following substitution are required:
\[p_1=q^{-1},\qquad p_2=1,\qquad \kappa=-1.
\]

\section{Solutions of the equivariant qDE and qKZ difference equations}\label{SEC6}
\subsection{Topological--enumerative morphism}\label{sectopenmor} Let $(\eps_1,\dots,\eps_{\binom{n}{k}})$ be a fixed basis of $H^\bullet_T(G(k,n),\C)$ as module over $H^\bullet_T({\rm pt},\C)\cong \C[\bm z]$. Set $\eta_{ab}=\eta(\eps_a,\eps_b)$, with $a,b=1,\dots,\binom{n}{k}$, $\eta=(\eta_{ab})_{a,b}$, and $\eta^{-1}=(\eta^{ab})_{ab}$. Introduce the dual basis $(\eps^1,\dots,\eps^{\binom{n}{k}})$, defined by $\eps^a=\sum_b\eta^{ab}\eps_b$.
\vskip2mm
Define $\mc S\in\End_{\C[\bm z]}\left(H^\bullet_T(G(k,n),\C)\right)[\![\bm p]\!][\![\kappa^{-1}]\!]$ by
\beq
\mc S(\bm z; \bm p;\kappa)a=\left[\widetilde{S}(\bm z;\bm p;\kappa)\circ \left(p_1^{c_1^{T}(E_1)}p_2^{c_1^{T}(E_2)}\right)^\frac{1}{\kappa}\right]a,\quad a\in H^\bullet_T(G(k,n),\C),
\eeq
where 
\begin{multline}
\label{tildes}
\widetilde{S}(\bm z; \bm p;\kappa)a=a+\sum_{\bt\in {\rm Eff}\setminus\{0\}}\sum_\mu p_1^{\int_\bt c_1(E_1)}p_2^{\int_\bt c_1(E_2)}\Big \langle\eps_\mu,\frac{a}{\kappa-\psi}\Big\rangle^{G(k,n), {T}}_{0,2,\bt}\eps^\mu\\
=a+\sum_{\bt\in {\rm Eff}\setminus\{0\}}\sum_\mu\sum_{j=0}^\infty\kappa^{-j-1}p_1^{\int_\bt c_1(E_1)}p_2^{\int_\bt c_1(E_2)}\langle\eps_\mu,\,\tau_ja\rangle^{G(k,n), {T}}_{0,2,\bt}\eps^\mu,\quad a\in H^\bullet_T(G(k,n),\C).
\end{multline}
\begin{rem}
The series $\widetilde S(\bm z; \bm p;\kappa)$, and consequently $S(\bm z; \bm p;\kappa)$, is independent of the chosen fixed basis $(\eps_1,\dots,\eps_{\binom{n}{k}})$. The $\bm z$-dependence of $\widetilde S(\bm z; \bm p;\kappa)$ and $S(\bm z; \bm p;\kappa)$ is hidden in the Gromov--Witten invariants.
\end{rem}
\begin{thm}\label{thmtopsol}
For each $a\in H^\bullet_T(G(k,n),\C)$, the $H^\bullet_T(G(k,n),\C)$-valued function $\mc S(\bm z;\bm p;\kappa)a$ is a solution of the qDE \eqref{qde.0}.
\end{thm}
\proof Let us prove the validity of \eqref{qde.0} for $i=1$ (the proof for the $i=2$ is identical). Consider the morphism $R\colon H^\bullet_T(G(k,n),\C)\to H^\bullet_T(G(k,n),\C)$, $a\mapsto c_1^{T}(E)a$ (equivariant cup product).
We have
\[\kappa p_1\frac{\der}{\der p_1}\mc S(\bm z;\bm p;\kappa)a=\left[\kappa p_1\frac{\der\widetilde S}{\der p_1}(\bm z;\bm p;\kappa)+\widetilde S(\bm z;\bm p;\kappa)\circ R)\right]\circ\left(p_1^{c_1^{T}(E_1)}p_2^{c_1^{T}(E_2)}\right)^\frac{1}{\kappa}a.
\]We need to prove
\[\kappa p_1\frac{\der\widetilde S}{\der p_1}(\bm z;\bm p;\kappa)+\widetilde S(\bm z;\bm p;\kappa)\circ R=[c_1^{T}(E_1)*_{q}]\circ \widetilde S(\bm z;\bm p;\kappa).
\]For any $a\in H^\bullet_T(G(k,n),\C)$, we have
\begin{multline*}
\left[\kappa p_1\frac{\der\widetilde S}{\der p_1}(\bm z;\bm p;\kappa)+\widetilde S(\bm z;\bm p;\kappa)\circ R\right]a\\
=\sum_{\bt\in {\rm Eff}\setminus\{0\}}\sum_\mu\sum_{j=0}^\infty\kappa^{-j}\left(\int_\bt c_1(E_1)\right)p_1^{\int_\bt c_1(E_1)}p_2^{\int_\bt c_1(E_2)}\langle\eps_\mu,\,\tau_ja\rangle^{G(k,n), {T}}_{0,2,\bt}\eps^\mu\\
+c_1^{T}(E_1)a+\sum_{\bt\in {\rm Eff}\setminus\{0\}}\sum_\mu\sum_{j=0}^\infty\kappa^{-j-1}p_1^{\int_\bt c_1(E_1)}p_2^{\int_\bt c_1(E_2)}\langle\eps_\mu,\,\tau_jc_1^{T}(E_1)a\rangle^{G(k,n), {T}}_{0,2,\bt}\eps^\mu.
\end{multline*}
By the divisor property of equivariant Gromov--Witten invariants, we have: 
\[\left(\int_\bt c_1(E_1)\right)\langle\eps_\mu,\,a\rangle^{G(k,n), {T}}_{0,2,\bt}=\langle c_1^{T}(E_1),\eps_\mu,\,a\rangle^{G(k,n), {T}}_{0,3,\bt},
\]
\[\left(\int_\bt c_1(E_1)\right)\langle\eps_\mu,\,\tau_ja\rangle^{G(k,n), {T}}_{0,2,\bt}+\langle\eps_\mu,\,\tau_{j-1}c_1^{T}(E_1)a\rangle^{G(k,n), {T}}_{0,2,\bt}=\langle c_1^{T}(E_1),\eps_\mu,\,\tau_ja\rangle^{G(k,n), {T}}_{0,3,\bt},\quad j\geq 1.
\]So, we have
\begin{multline}\label{topsol1}
\left[\kappa p_1\frac{\der\widetilde S}{\der p_1}(\bm z;\bm p;\kappa)+\widetilde S(\bm z;\bm p;\kappa)\circ R\right]a\\
=c_1^{T}(E_1)a+\sum_{\bt\in {\rm Eff}\setminus\{0\}}\sum_\mu\sum_{j=0}^\infty\kappa^{-j}p_1^{\int_\bt c_1(E_1)}p_2^{\int_\bt c_1(E_2)}\langle c_1^{T}(E_1),\eps_\mu,\,\tau_ja\rangle^{G(k,n), {T}}_{0,3,\bt}\eps^\mu.
\end{multline}
On the other hand, by \eqref{seqp}, we have
\begin{multline}\label{topsol2}
c_1^{T}(E_1)*_{q}S(\bm z;\bm p;\kappa)a=\sum_{\bt_1\in{\rm Eff}}\sum_\mu\langle c_1^{T}(E_1),a,\eps_\mu\rangle^{G(k,n), {T}}_{0,3,\bt_1}\eps^\mu\\+\sum_{\bt\in{\rm Eff}\setminus\{0\}}\sum_{\bt_1\in{\rm Eff}}\sum_{j=0}^\infty\sum_{\mu,\la}\kappa^{-j-1}p_1^{\int_{\bt+\bt_1}c_1(E_1)}p_2^{\int_{\bt+\bt_1}c_1(E_2)}\langle\eps_\mu,\tau_j a\rangle_{0,2,\bt}\langle c_1^{T}(E_1),\eps^\mu,\eps_\la\rangle^{G(k,n), {T}}_{0,3,\bt_1}\eps^\la.
\end{multline}
By the splitting property of equivariant Gromov--Witten invariants, we have:
\[\langle\tau_{j+1}a,c_1^{T}(E_1),\eps_\la\rangle_{0,3,\bt}=\sum_{\mu}\sum_{\bt=\bt_1+\bt_2}\langle\eps_\mu,\tau_j a\rangle_{0,2,\bt_1}\langle c_1^{T}(E_1),\eps^\mu,\eps_\la\rangle^{G(k,n), {T}}_{0,3,\bt_2},\quad j\geq 0.
\]Then \eqref{topsol1} and \eqref{topsol2} coincide. 
\endproof

\subsection{Levelt fundamental system of solutions}\label{seclev} We say that a function $f(\bm p)$ is {\it entire in $\bm p_\div$} if $f(\bm p)=g(p_2/p_1)$ for an entire function $g(s)$, and $f(\oslash):=g(0)$.
\vskip1,5mm
The dynamical operators $X_1(\bm z;\bm p)$ and $X_2(\bm z;\bm p)$, defined in \eqref{dynop1} and \eqref{dynop2}, are entire in $\bm p_\div$, and 
\beq
X_1(\bm z;\oslash)=\sum_{a=1}^nz_a e_{11}^{(a)}+\sum_{1\leq b<a\leq n}Q^{a,b}_{1,2},\qquad
X_2(\bm z;\oslash)=\sum_{a=1}^nz_a e_{22}^{(a)}-\sum_{1\leq a<b\leq n}Q^{a,b}_{2,1}.
\eeq

For $I\in\mc I_{k,n}$, with $I=(I_1,I_2)$, set $E_I(\bm z)=\sum_{a\in I_1}z_a$. For $I,J\in\mc I_{k,n}$, denote by $D_{I,J}$ the set of points $\bm z$ such that $E_I(\bm z)-E_J(\bm z)\in\Z_{\geq 0}$ and $E_I(\bm z)\neq E_J(\bm z)$. Set $D_{k,n}:=\bigcup_{I,J\in\mc I_{k,n}}D_{I,J}$.

\begin{thm}\label{thmlev}\cite[Sec.\,5.5]{TV23} \begin{enumerate}
\item For any $\bm z$ such that $z_a-z_b\notin \kappa\Z_{\neq 0}$ for all $1\leq a<b\leq n$, there exists an $\on{End}\left((\C^2)^{\otm n}_{k,n-k}\right)$-valued function $\Psi^\bullet(\bm z;\bm p;\kappa)$ entire in $\bm p_\div$ such that $\Psi^\bullet(\bm z;\oslash;\kappa)$ is the identity operator, and the function
\beq\label{leveq1}
\Hat\Psi(\bm z;\bm p;\kappa)=\Psi^\bullet(\bm z;\bm p;\kappa)p_1^{X_1(\bm z;\oslash)/\kappa}p_2^{X_2(\bm z;\oslash)/\kappa}
\eeq is a solution of the dynamical differential equations \eqref{dyn.0}. For given $\bm z$, the function $\Psi^\bullet(\bm z;\bm p;\kappa)$ with the specified properties is unique if and only if $\bm z\notin D_{k,n}$. Furthermore
\beq
\det \Hat\Psi(\bm z;\bm p;\kappa)=p_1^{d_1\sum_{a=1}^nz_a/\kappa}p_2^{d_2\sum_{a=1}^nz_a/\kappa},\quad d_1=\frac{(n-1)!k}{k!(n-k)!},\quad d_2=\frac{(n-1)!(n-k)}{k!(n-k)!}. 
\eeq
\item Define the function $\Psi^\bullet(\bm z;\bm p;\kappa)$ for generic $\bm z$ as in item (1). Then $\Psi^\bullet(\bm z;\bm p;\kappa)$ is holomorphic in $\bm z$ if $z_a-z_b\notin\kappa\Z_{\neq 0}$ for all $1\leq a<b\leq n$. The singularities of $\Psi^\bullet(\bm z;\bm p;\kappa)$ at the hyperplanes $z_a-z_b\in\kappa\Z_{\neq 0}$ are simple poles.\qed
\end{enumerate}
\end{thm}

We call $\Hat\Psi(\bm z;\bm p;\kappa)$ the {\it Levelt fundamental system of solutions} of the dynamical equations \eqref{dyn.0} on $(\C^2)^{\otm n}_{k,n-k}$.
\vskip2mm
Consider the vector bundle $E_{k,n}\to \C^n\times \C^2$ with fiber over a point $(\bm z_o,\bm p_o)$ given by the $\C$-vector space $\on{End}\left(H^\bullet_T(G(k,n),\C)_{\bm z_o}\right)$. Consider the section $\on{Stab}_{k,n}\Hat\Psi$ defined by
\beq
\on{Stab}_{k,n}\Hat\Psi(\bm z;\bm p;\kappa)=\on{Stab}_{k,n,\bm z}\circ\Hat\Psi(\bm z;\bm p;\kappa)\circ \left(\on{Stab}_{k,n,\bm z}\right)^{-1}.
\eeq
By Theorem \ref{thmlev} and Lemma \ref{stabintzo}, for any section $f$ of $H_{k,n}$ not depending on $\bm p$, the section $\on{Stab}_{k,n}\Hat\Psi f$ of $H_{k,n}$ with values $\on{Stab}_{k,n}\Hat\Psi(\bm z;\bm p;\kappa) f(\bm z)$ is a solution of the quantum differential equations \eqref{qde.0}. The section $\on{Stab}_{k,n}\Hat\Psi$ of $E_{k,n}$ is called the {\it Levelt fundamental system of solutions} of the equivariant quantum differential equations \eqref{qde.0}.

\begin{thm}\label{sollev} We have $\on{Stab}_{k,n}\Psi^\bullet(\bm z;\bm p;\kappa)=\Tilde{\mc S}(\bm z;\bm p;\kappa)$. Hence, 
the topological enumerative morphism $\mc S(\bm z;\bm p;\kappa)$ equals the Levelt fundamental system of solutions $\on{Stab}_{k,n}\Hat\Psi(\bm z;\bm p;\kappa)$.
\end{thm}
\proof
By formulae \eqref{semqm} and \eqref{leveq1}, we have 
\[\on{Stab}_{k,n}\Hat\Psi(\bm z;\bm p;\kappa)=\left(\on{Stab}_{k,n}\Psi^\bullet(\bm z;\bm p;\kappa)\right)p_1^{c_1^{T}(E_1)/\kappa}p_2^{c_1^{T}(E_2)/\kappa},
\]where the operator $p_1^{c_1^{T}(E_1)/\kappa}p_2^{c_1^{T}(E_2)/\kappa}$ acts on the fiber $H^\bullet_T(G(k,n),\C)_{\bm z}$ as classical equivariant multiplication $a\mapsto p_1^{c_1^{T}(E_1)/\kappa}p_2^{c_1^{T}(E_2)/\kappa}a$. Notice that $\Tilde S(\bm z;\bm p;\kappa)$ defined by \eqref{tildes} is entire in $\bm p_\div$, since $\int_\bt c_1(E_1)=-\int_\bt c_1(E_2)$ (non-equivariant first Chern classes). Moreover, $\Tilde S(\bm z;\oslash;\kappa)$ equals the identity operator. The result then follows by the uniqueness statement asserted in Theorem \ref{thmlev}.
\endproof

\begin{rem}
Theorem \ref{sollev} proves Conjecture 6.32 formulated in \cite{TV23} for the Grassmannian $G(k,n)$. For more general partial flag varieties $\mc F_{\bm \la}$ the same result holds true, and it can be proved similarly.
\end{rem}

\subsection{Multidimensional hypergeometric solutions} Introduce the variables $\bm t=(t_1,\dots, t_k)$. For each $I\in\mc I_{k,n}$, with $I=(I_1,I_2)$ and $I_1=\{i_1<\dots<i_k\}$, set 
\[U_I(\bm t;\bm z):=\prod_{a=1}^k\left(\prod_{c=1}^{i_a-1}(t_a-z_c)\prod_{b=a+1}^k\frac{1}{t_b-t_a}\right).
\]For $I\in\mc I_{k,n}$, define the {\it weight function} 
\beq W_I(\bm t;\bm z):={\rm Sym}_{\bm t}U_I(\bm t;\bm z)=\sum_{\si\in S_k}U_I(t_{\si(1)},\dots, t_{\si(k)};\bm z).
\eeq
The {\it master function} $\Phi_{k,n}$ is defined by
\beq
\Phi_{k,n}(\bm t;\bm z;\bm p;\kappa)=(\kappa^{-k}p_2)^{\sum_{a=1}^nz_a/\kappa}\left(\kappa^n\frac{p_1}{p_2}\right)^{\sum_{a=1}^kt_a/\kappa}
\prod_{a=1}^k\left(\prod_{\substack{b=1\\b\neq a}}^k\frac{1}{\Gm\left(\frac{t_a-t_b}{\kappa}\right)}\prod_{c=1}^n\Gm\left(\frac{t_a-z_c}{\kappa}\right)\right).
\eeq
For $J\in\mc I_{k,n}$, with $J=(J_1,J_2)$ and $J_1=\{j_1<\dots<j_k\}$, set $\Si_J=(z_{j_1},\dots,z_{j_k})$. For $J\in \mc I_{k,n}$, and any polynomial $f(\bm t;\bm z)$, define the {\it Jackson integral}
\[\mc M_J\left(\Phi_{k,n}f\right)(\bm z;\bm p;\kappa):=\sum_{\bm\ell\in\Z^k}\underset{\bm t=\Si_J-\kappa\bm\ell}{\rm Res}\left(\Phi_{k,n}(\bm t;\bm z;\bm p;\kappa)f(\bm t;\bm z)\right).
\]Here the multi-dimensional residue $\underset{\bm t=\bm s}{\rm Res}\,f(\bm t)$, $\bm s\in\C^k$, of a function $f(\bm t)$ stands for the iterated residue
\[\underset{\bm t=\bm s}{\rm Res}\,f(\bm t)=\underset{t_1=s_1}{\rm Res}\,\underset{t_2=s_2}{\rm Res}\dots\underset{t_k=s_k}{\rm Res}\,f(\bm t).
\]Define the $(\C^2)^{\otm n}_{k,n-k}$-valued function
\[\Psi_J(\bm z;\bm p;\kappa)=\kappa^{-k(n-k)-k}\sum_{I\in\mc I_{k,n}}\mc M_J\left(\Phi_{k,n}W_I\right)(\bm z;\bm p;\kappa)v_I,\quad J\in\mc I_{k,n}.
\]
\begin{thm}\cite[Th.\,5.11]{TV23}
The function $\Psi_J(\bm z;\bm p;\kappa)$ is a solution of the joint system of qKZ difference equations \eqref{qKZ.0} and dynamical differential equations \eqref{dyn.0}.\qed
\end{thm}

The solutions $\Psi_J$, with $J\in\mc I_{k,n}$, are called {\it multi-dimensional hypergeometric solutions} of the joint system \eqref{qKZ.0} and \eqref{dyn.0}.
\vskip2mm
Since $\Psi_J(\bm z;\bm p;\kappa)$ is a $(\C^2)^{\otm n}_{k,n-k}$-valued function, for each fixed $\kappa\in\C^*$ it defines a section of the trivial bundle $U_{k,n}$. Denote by ${\rm Stab}_{k,n}\Psi_J$ the corresponding section of $H_{k,n}$, with values
\[[{\rm Stab}_{k,n}\Psi_J](\bm z;\bm p;\kappa)={\rm Stab}_{k,n;\bm z}[\Psi_J(\bm z;\bm p;\kappa)].
\]
\begin{cor}
The function ${\rm Stab}_{k,n}\Psi_J$, with $J\in\mc I_{k,n}$, is a section of $H_{k,n}$ solving the joint system of qDE \eqref{qde.0} and qKZ difference equations \eqref{qkz.H}.
\end{cor}
\proof
The result follows from Theorem \ref{stabint} and the definition of the qKZ difference equations \eqref{qkz.H}.
\endproof

\subsection{Identification of solutions with $K$-theoretical classes}\label{secideksol}For a Laurent polynomial $P\in\C[\bm\Gm^{\pm 1}]^{S_k\times S_{n-k}}\otm\C[\bm Z^{\pm 1}]$, define the polynomial $\acute{P}(\bm\gm;\bm z;\kappa)$ by
\beq
\acute{P}(\bm\gm;\bm z;\kappa):=P(\acute{\bm\gm};\acute{\bm{z}}),
\eeq see \eqref{acute1} and \eqref{acute2}, and set
\beq
\Psi_P(\bm z;\bm p;\kappa):=\sum_{J\in \mc I_{k,n}}\acute{P}(\bm z_{\si_J};\bm z;\kappa)\Psi_J(\bm z;\bm p;\kappa).
\eeq
Recall the ring presentation for the $K$-theory $K^T_0(G(k,n))_\C$ given by \eqref{kthpres}. Each polynomial $P\in\C[\bm\Gm^{\pm 1}]^{S_k\times S_{n-k}}\otm\C[\bm Z^{\pm 1}]$ defines a $K$-theoretical class $[P]$.

\begin{thm}\cite[Sec.\,5.4]{TV23} The function $\Psi_P(\bm z;\bm p;\kappa)$ is a solution of the joint system of qKZ difference \eqref{qKZ.0} and dynamical differential equations \eqref{dyn.0}. Furthermore, for $P\in\C[\bm\Gm^{\pm 1}]^{S_k\times S_{n-k}}\otm\C[\bm Z^{\pm 1}]$, the function $\Psi_P(\bm z;\bm p;\kappa)$ is entire in $\bm z$ and holomorphic in $\bm p$ provided that a branch of $\log p_i$ is fixed for each $i=1,2$. Moreover, the function $\Psi_P(\bm z;\bm p;\kappa)$ only depends on the class $[P]$.  \qed
\end{thm}

Denote by $\mathscr S_{H_{k,n}}$ the space spanned over $\C$ by the solutions $\on{Stab}_{k,n}\Psi_P$ of the joint system of equations \eqref{qde.0} and \eqref{qkz.H}. Each element of $\mathscr S_{H_{k,n}}$ is a section of $H_{k,n}$ holomorphic in $\bm p$, when a branch of $\log p_i$ is fixed for $i=1,2$, and entire in $\bm z$. The space $\mathscr S_{H_{k,n}}$ is a $\C[\acute{\bm z}^{\pm 1}]$-module, with $f(\acute{\bm z})$ acting as multiplication by $\acute{f}(\bm z;\kappa)$.
\vskip1,5mm
We have a well-defined map
\beq\label{defmukn}
\mu_{k,n}\colon K^T_0(G(k,n))_\C\to \mathscr S_{H_{k,n}},\quad [P]\mapsto\on{Stab}_{k,n}\Psi_P.
\eeq
\begin{thm}\cite[Sec.\,5.4]{TV23}
The map $\mu_{k,n}$ is an isomorphism of $\C[\acute{\bm z}^{\pm 1}]$-modules.\qed
\end{thm}

\begin{example}\label{ex:intforpn}
Consider the case of the projective space $G(1,n)=\Pb^{n-1}$. Let $\gm=\gm_{1,1}^\Pb$. The equivariant cohomology ring has presentation
\[H^\bullet_T(\Pb^{n-1},\C)\cong\C[\gm]\otm \C[\bm z]\Big/\Big\langle\prod_{i=1}^n(\gm -z_i)\Big\rangle.
\]Denote by $[a]$ the element $(\{a\},\{1,\dots,a-1,a+1,\dots,n\})\in\mc I_{k,n}$. The polynomials $\on{Stab}_{[a]}(\gm;\bm z)$ and $\on{Stab}^{\on{op}}_{[a]}(\gm;\bm z)$ are
\[\on{Stab}_{[a]}(\gm;\bm z)=\prod_{c=a+1}^n(\gm-z_c),\quad \on{Stab}^{\on{op}}_{[a]}(\gm;\bm z)=\prod_{c=1}^{a-1}(\gm-z_c),\quad a=1,\dots,n.
\]The weight functions $W_{[a]}(t;\bm z)$ are
\[W_{[a]}(t;\bm z)=\prod_{c=1}^{a-1}(t-z_c),\quad a=1,\dots,n.
\]Consider the polynomial 
\[\bar{W}(t;\gm;\bm z):=\frac{1}{t-\gm}\left(\prod_{a=1}^n(t-z_a)-\prod_{a=1}^n(\gm-z_a)\right).
\]We have the equality of polynomials 
\beq\label{vectorvaluedweight}
\bar{W}(t;\gm;\bm z)=\sum_{a=1}^nW_{[a]}(t;\bm z)\on{Stab}_{[a]}(\gm;\bm z).
\eeq\begin{rem}The image $[\bar{W}(t;\gm;\bm z)]$ in $\C[t]\otm H^\bullet_T(\Pb^{n-1},\C)$ also satisfies
\[\left[\bar{W}(t;\gm;\bm z)\right]=\left[\Hat W(t;\bm \gm)\right],\quad \Hat W(t;\bm\gm):=\prod_{a=1}^{n-1}(t-\gm_{2,a}).
\]\end{rem}
Consider the trivial vector bundle $H_{1,n}\to\C^n\times\C^2$, with fiber over a point $(\bm z^o,\bm p^o)$ given by 
\[H^\bullet_T(\Pb^{n-1},\C)_{\bm z^o}=\C[\gm]\Big/\Big\langle\prod_{i=1}^n(\gm -z^o_i)\Big\rangle.
\]Given a polynomial $f(\gm;\bm z)$, denote by $[f(\gm;\bm z^o)]_{\bm z^o}$ its class in $H^\bullet_T(\Pb^{n-1},\C)_{\bm z^o}$. Taking into account formula \eqref{vectorvaluedweight}, we have an integral formula for the section $\on{Stab}_{1,n}\Psi_P(\bm z;\bm p;\kappa)$,
\beq\label{intpn-1}
\on{Stab}_{1,n}\Psi_P(\bm z;\bm p;\kappa)=\frac{\kappa^{-n}}{2\pi\sqrt{-1}}\int_C\acute{P}(t;\bm z;\kappa)\Phi_{1,n}(t;\bm z;\bm p;\kappa)[\bar W(t;\gm;\bm z)]_{\bm z}dt,
\eeq
where
\beq
\Phi_{1,n}(t;\bm z;\bm p;\kappa)=(p_2/\kappa)^{\sum_{a=1}^nz_a/\kappa}(\kappa^np_1/p_2)^{t/\kappa}\prod_{a=1}^n\Gm((t-z_a)/\kappa),
\eeq and the contour $C$ encircles the poles of $\Gm((t-z_a)/\kappa)$ counterclockwise. For instance, we can take
\[C=\{\kappa(A-s^2+s\sqrt{-1})\colon s\in\R\},
\]with $A$ a sufficiently large number.
\end{example}

\subsection{Solutions via equivariant Satake correspondence} The results of Section \ref{sec2} give us a further way to construct solutions of the joint system of qDE \eqref{qde.0} and qKZ difference equations \eqref{qkz.H} for $G(k,n)$, starting from solutions of the corresponding joint system for $\Pb^{n-1}$.
\vskip2mm
For short, we set $\si_{\bm\la}=[\Om_{\bm\la}]$, for any partition $\bm\la$. In what follows we will identify the cohomology of $G(k,n)$ with the $k$-th exterior power of the cohomology of $\Pb^{n-1}$. Schubert classes of $\Pb^{n-1}$ will be denoted by $\si_{\bm\la}^{\Pb}$, whereas Schubert classes of $G(k,n)$ will be denoted by $\si_{\bm\la}^G$. The tautological bundle on $G(k,n)$ (resp.\,\,$\Pb^{n-1}$) will be denoted by $E_1^G$ (resp.\,\,$E_1^{\Pb}$), and similarly for the quotient bundles $E_2^G$ (resp.\,\,$E_2^{\Pb}$).
\vskip1,5mm
Consider the $\C[\bm z]$-module morphism $\thi_{k,n}\colon\bigwedge\nolimits^k_{\C[\bm z]}H^\bullet_T(\Pb^{n-1},\C)\to H^\bullet_T(G(k,n),\C)$ defined as the composition
\beq\label{sathi}
\thi_{k,n}={\rm Stab}_{k,n}\circ \theta_{k,n}\circ\left(\bigwedge\nolimits^k{\rm Stab}_{1,n}\right)^{-1},
\eeq
where $\theta_{k,n}$ is defined in \eqref{thk}, and ${\rm Stab}_{k,n}$ in \eqref{stabkn}. 
\begin{prop}
The map $\thi_{k,n}$ is an isomorphism of $\C[\bm z]$-modules. It acts on Schubert classes as follows:
\[\thi_{k,n}(\si_{\la_k}^\Pb\wedge \si_{\la_{k-1}+1}^\Pb\wedge\dots\wedge\si_{\la_1+k-1}^\Pb)=\si_{\bm\la}^G,\quad \bm\la\subseteq k\times (n-k).
\]
\end{prop}
\proof
The map $\thi_{k,n}$ is an isomorphism because $\theta_{k,n}$ and ${\rm Stab}_{k,n}$ are isomorphisms for any $k,n$. We have
\begin{multline*}\left(\bigwedge\nolimits^k{\rm Stab}_{1,n}\right)^{-1}[\si_{\la_k}^\Pb\wedge \si_{\la_{k-1}+1}^\Pb\wedge\dots\wedge\si_{\la_1+k-1}^\Pb]=(-1)^{|\bm\la|+\frac{k(k-1)}{2}}v_{[n-\la_k]}\wedge\dots\wedge v_{[n-\la_1-k+1]}\\
=(-1)^{|\bm\la|}v_{[n-\la_1-k+1]}\wedge\dots\wedge v_{[n-\la_k]}.
\end{multline*}
If $I=(\{n-\la_1-k+1,\dots,n-\la_k\},\{1,\dots,n\}\setminus\{n-\la_1-k+1,\dots,n-\la_k\})$, then $\si_I(j)=n-\la_j-k+j$ for $j=1,\dots,k$, so that $\bm\la^{\si_I}_j=n-k-\la_{k-j+1}$ for $j=1,\dots,k$. Hence, we have $(\bm\la^{\si_I})^\vee=\bm\la$, and $({\rm Stab}_{k,n}\circ\theta_{k,n})v_I=(-1)^{|\bm\la|}\si_{\bm\la}^G$. The claim follows.
\endproof

Denote $\gm=\gm_{1,1}^\Pb$. Given a tuple of indeterminates $\bm x=(x_1,\dots,x_n)$, set $\frak D_k(\bm x)$ to be the Vandermonde determinant in the first $k$ variables, that is
\beq\label{vandet}
\frak D_k(\bm x)=\det\left(\begin{array}{cccc}
x_1^{k-1}&\dots&x_1&1\\
x_2^{k-1}&\dots&x_2&1\\
\vdots&&&\vdots\\
x_k^{k-1}&\dots&x_k&1
\end{array}\right)=\prod_{1\leq i< j\leq k}(x_i-x_j).
\eeq
In terms of the polynomial presentations \eqref{presentation}, the morphism $\thi_{k,n}$ takes the following form
\begin{multline}\label{polthik}
\thi_{k,n}\colon \bigwedge\nolimits^k_{\C[\bm z]}\left(\C[\gm,\bm z]\Big/\Big\langle\prod_{i=1}^n(\gm-z_i)\Big\rangle\right)\to \C[\bm\gm^G]^{S_{k}\times S_{n-k}}[\bm z]\Big/\Big\langle   e_j(\bm z)-\sum_{h=0}^j e_h(\bm\gm_1^G)e_{j-h}(\bm\gm_2^G)\Big\rangle_{j=1}^n,\\
f_1(\gm;\bm z)\wedge\dots\wedge f_k(\gm;\bm z)\mapsto \frac{\det
\left(\begin{array}{ccc}f_1(\gm_{1,1}^G;\bm z)&\dots&f_k(\gm_{1,1}^G;\bm z)\\
\vdots&&\vdots\\
f_1(\gm_{1,k}^G;\bm z)&\dots&f_k(\gm_{1,k}^G;\bm z)
\end{array}\right)}{\frak D_k(\bm\gm^G)}.
\end{multline}
By extensions of scalars, the isomorphism $\thi_{k,n}$ induces $\thi_{k,n}\colon \bigwedge\nolimits^k_{\mc O_{L_\kappa}} H^{(\kappa)}_T(\Pb^{n-1})\to H^{(\kappa)}_T(G(k,n))$ defined by \eqref{polthik}.

Given $a\in\{1,\dots,n\}$, the element $[a]\in\mc I_{1,n}$ is defined by $[a]=(\{a\},\{1,\dots,\hat a,\dots,n\})$. Let $a_1,\dots,a_k\in\{1,\dots,n\}$ pairwise distinct, we have the idempotents $\Dl_{[a_i]}\in H^{(\kappa)}_T(\Pb^{n-1})$, see Section \ref{secext}.

\begin{prop}\label{propidemp}
We have 
\[\thi_{k,n}\left(\Dl_{[a_1]}\wedge\dots\wedge\Dl_{[a_k]}\right)=\frak D_k(\bm z_{\si_I})^{-1}\Dl_{I},\quad I=(\{a_1,\dots,a_k\},\{1,\dots,n\}\setminus\{a_1,\dots,a_k\}).
\]
\end{prop}
\proof
By definition, we have $\thi_{k,n}\left(\Dl_{[a_1]}\wedge\dots\wedge\Dl_{[a_k]}\right)=\det(\Dl_{[a_i]}(\gm_{1,j}^G))_{i,j=1}^k\Big/\frak D_k(\bm\gm^G)$. Given $J\in\mc I_{k,n}$, the restriction $\thi_{k,n}\left(\Dl_{[a_1]}\wedge\dots\wedge\Dl_{[a_k]}\right)|_{{\rm pt}_J}$ equals
\[\pushQED{\qed}
\det\left(\Dl_{[a_i]}(z_{\si_J(j)})\right)_{i,j=1}^k\Big/\frak D_k(\bm z_{\si_J})=\frak D_k(\bm z_{\si_I})^{-1}\dl_{IJ}.\qedhere
\popQED
\]
The Poincar\'e metric $\eta^\Pb$ on $H^\bullet_T(\Pb^{n-1},\C)$ induces a metric $\eta^{\wedge\Pb}$ on $\bigwedge\nolimits^k_{\C[\bm z]} H^\bullet_T(\Pb^{n-1},\C)$, defined by
\[\eta^{\wedge\Pb}(v_1\wedge\dots\wedge v_k,w_1\wedge\dots\wedge w_k)=\det\left(\eta^\Pb(v_i,w_j)\right)_{i,j=1}^k.
\]
\begin{prop}\label{isometry1}
The isomorphism $\thi_{k,n}\colon \left(\bigwedge\nolimits^k_{\C[\bm z]} H^\bullet_T(\Pb^{n-1},\C),\,(-1)^{\binom{k}{2}}\eta^{\wedge\Pb}\right)\to\left(H^\bullet_T(G(k,n),\C),\eta^G\right)$ is an isometry.
\end{prop}
\proof Let us compute the Gram matrix of $\eta^G$ in the idempotent basis $(\Dl_I)_{I\in I_{k,n}}$, and the Gram matrix of $\eta^{\wedge\Pb}$ in the basis induced by the idempotent basis $(\Dl_{[i]})_{i=1}^n$. In both cases, the only nonzero entries are the diagonal ones.
For $I\in\mc I_{k,n}$, we have
\[\eta^G(\Dl_I,\Dl_I)=\eta^G(\Dl_I,1)=\int_{G(k,n)}^{\rm eq}\Dl_I=(-1)^{k(n-k)}\prod_{i\in I_1}\prod_{j\in I_2}\frac{1}{z_i-z_j}.
\]
On the other hand, if $I_1=\{i_1<\dots<i_k\}$, we have
\[\det\left(\eta^\Pb(\Dl_{[i_a]},\Dl_{[i_b]})\right)_{a,b=1}^k=\prod_{a=1}^k(-1)^{n-1}\prod_{j\neq i_a}\frac{1}{z_{i_a}-z_j}=(-1)^{k(n-1)}\prod_{a=1}^k\prod_{j\neq i_a}\frac{1}{z_{i_a}-z_j}.
\]The last product can be manipulated as follows:
\begin{multline*}\prod_{a=1}^k\prod_{j\neq i_a}\frac{1}{z_{i_a}-z_j}=\prod_{a=1}^k\prod_{j\in I_1,j<i_a}\frac{1}{z_{i_a}-z_j}\prod_{j\in I_1,j>i_a}\frac{1}{z_{i_a}-z_j}\prod_{j\not\in I_1}\frac{1}{z_{i_a}-z_j}\\=(-1)^{\binom{k}{2}}\frak D_k(\bm z_{\si_I})^{-2}\prod_{i\in I_1}\prod_{j\in I_2}\frac{1}{z_i-z_j}.
\end{multline*}
The claim follows, by Proposition \ref{propidemp}.
\endproof

\begin{example}\label{ex:p2g23-1}
Let $n=3$, and $k=1,2$. Let us compare the case of $\Pb^2=G(1,3)$ and $G(2,3)$. 
\vskip2mm
For $\Pb^2$, denote $\gm=\gm_{1,1}$. We have partitions $\bm\la\subseteq 1\times 2 $, that is 
\[\bm\la=(0),\quad \bm\la=(1),\quad \bm\la=(2),
\]with corresponding Grassmannian permutations
\[\si=\on{id},\quad \si=(213),\quad \si=(312),
\]
and corresponding partitions
\[I=[1],\quad I=[2],\quad I=[3].
\]
The Schubert basis of $H^\bullet_T(\Pb^2,\C)$ is
\begin{align*}&\si_{(0)}^\Pb=(-1)^{|0|}s_{(0)}(\gm|\bm z_{\si_0})=1,\\
&\si_{(1)}^\Pb=(-1)^{|(1)|}s_{(1)}(\gm|\bm z_{\si_0})=-(\gm-z_3),\\
&\si_{(2)}^\Pb=(-1)^{|(2)|}s_{(2)}(\gm|\bm z_{\si_0})=(\gm-z_2)(\gm-z_3).
\end{align*}
The stable envelop map $\on{Stab}_{1,3}\colon (\C^2)^{\otm 3}_{1,2}\otm \C[\bm z]\to H^\bullet_T(\Pb^2,\C)$ is given by
\begin{align*}
&v_{[1]}=v_1\otm v_2\otm v_2\mapsto \on{Stab}_{[1]}=\si_{(2)}^\Pb,\\
&v_{[2]}=v_2\otm v_1\otm v_2\mapsto \on{Stab}_{[2]}=-\si_{(1)}^\Pb,\\
&v_{[3]}=v_2\otm v_2\otm v_1\mapsto \on{Stab}_{[3]}=\si_{(0)}^\Pb.
\end{align*}
In \cite{CV21} the basis $(\on{Stab}_{[1]},\on{Stab}_{[2]},\on{Stab}_{[3]})$ is called $\bm g$-basis.\newline
The idempotent basis of $H^{(\kappa)}_T(\Pb^2)$ is given by
\[
\Dl_{[j]}=\prod_{i\neq j}\frac{\gm-z_i}{z_j-z_i},\quad j=1,2,3.
\]The Poincar\'e metric $\eta^\Pb$ is determined by the identities
\begin{align*}
\eta^\Pb(\Dl_{[i]},\Dl_{[j]})=\eta^\Pb(\Dl_{[i]}\Dl_{[j]},1)=\dl_{ij}\int_{\Pb^2}^{\rm eq}\Dl_{[i]}=\dl_{ij}\prod_{a\neq i}\frac{1}{z_i-z_a},\quad i,j=1,2,3.
\end{align*}
\vskip2mm
For $G(2,3)$, we have partitions $\bm\la\subseteq 2\times 1$, that is
\[\bm\la=(0),\quad \bm\la=(1),\quad \bm\la=(1,1),
\]with corresponding Grassmannian permutations
\[\si=\on{id},\quad \si=(132),\quad \si=(231),
\]and corresponding partitions
\[I=(\{1,2\},\{3\}),\qquad I=(\{1,3\},\{2\}),\qquad I=(\{2,3\},\{1\}).
\]The Schubert basis of $H^\bullet_T(G(2,3),\C)$ is
\begin{align*}
&\si_{(0)}^G=(-1)^{|(0)|}s_{(0)}(\bm\gm_1|\bm z_{\si_0})=1,\\
&\si_{(1)}^G=(-1)^{|(1)|}s_{(1)}(\bm\gm_1|\bm z_{\si_0})=-\frac{\left|\begin{array}{cc}
(\gm_{1,1}|\bm z_{\si_0})^2&(\gm_{1,1}|\bm z_{\si_0})^0\\
(\gm_{1,2}|\bm z_{\si_0})^2&(\gm_{1,2}|\bm z_{\si_0})^0
\end{array}\right|}{(\gm_{1,1}-\gm_{1,2})}\\&=-(\gm_{1,1}+\gm_{1,2}-z_2-z_3),\\
&\si_{(1,1)}^G=(-1)^{|(1,1)|}s_{(1,1)}(\bm\gm_1|\bm z_{\si_0})=\frac{\left|\begin{array}{cc}
(\gm_{1,1}|\bm z_{\si_0})^2&(\gm_{1,1}|\bm z_{\si_0})^1\\
(\gm_{1,2}|\bm z_{\si_0})^2&(\gm_{1,2}|\bm z_{\si_0})^1
\end{array}\right|}{(\gm_{1,1}-\gm_{1,2})}\\&=(\gm_{1,1}-z_3)(\gm_{1,2}-z_3).
\end{align*}The formulas for the Schubert basis can also be recovered from the Satake identification:
\[\si_{(0)}^G=\thi_{2,3}(\si^\Pb_0\wedge\si_1^\Pb)=\frac{\left|\begin{array}{cc}
1&-(\gm_{1,1}-z_3)\\
1&-(\gm_{1,2}-z_3)
\end{array}\right|}{\gm_{1,1}-\gm_{1,2}}=1,
\]
\[\si_{(1)}^G=\thi_{2,3}(\si_0^\Pb\wedge\si_2^{\Pb})=\frac{\left|\begin{array}{cc}
1&(\gm_{1,1}-z_2)(\gm_{1,1}-z_3)\\
1&(\gm_{1,2}-z_2)(\gm_{1,2}-z_3)
\end{array}\right|}{\gm_{1,1}-\gm_{1,2}}=-(\gm_{1,1}+\gm_{1,2}-z_2-z_3),
\]
\[\si_{(1,1)}^G=\thi_{2,3}(\si_1^\Pb\wedge\si_2^\Pb)=\frac{\left|\begin{array}{cc}
-(\gm_{1,1}-z_3)&(\gm_{1,1}-z_2)(\gm_{1,1}-z_3)\\
-(\gm_{1,2}-z_3)&(\gm_{1,2}-z_2)(\gm_{1,2}-z_3)
\end{array}\right|}{\gm_{1,1}-\gm_{1,2}}=(\gm_{1,1}-z_3)(\gm_{1,2}-z_3).
\]
The stable envelop map $\on{Stab}_{2,3}\colon (\C^2)^{\otm 3}_{2,1}\otm \C[\bm z]\to H^\bullet_T(G(2,3),\C)$ is given by
\begin{align*}
&v_{\{1,2\},\{3\}}=v_1\otm v_1\otm v_2\mapsto \on{Stab}_{\{1,2\},\{3\}}=\si_{(1,1)}^G,\\
&v_{\{1,3\},\{2\}}=v_1\otm v_2\otm v_1\mapsto \on{Stab}_{{\{1,3\},\{2\}}}=-\si_{(1)}^G,\\
&v_{\{2,3\},\{1\}}=v_2\otm v_1\otm v_1\mapsto \on{Stab}_{\{2,3\},\{1\}}=\si_{(0)}^G.
\end{align*}
The idempotent basis of $H^{(\kappa)}_T(G(2,3))$ is 
\[\Dl_{\{1,2\},\{3\}}=\frac{\prod_{j=1}^2(\gm_{1,j}-z_3)}{(z_1-z_3)(z_2-z_3)},\quad \Dl_{\{1,3\},\{2\}}=\frac{\prod_{j=1}^2(\gm_{1,j}-z_2)}{(z_1-z_2)(z_3-z_2)},\quad \Dl_{\{2,3\},\{1\}}=\frac{\prod_{j=1}^2(\gm_{1,j}-z_1)}{(z_2-z_1)(z_3-z_1)}.
\]We have
\[\thi_{2,3}(\Dl_{[1]}\wedge\Dl_{[2]})=\frac{(\gm_{1,1}-z_3)(\gm_{1,2}-z_3)}{(z_1-z_2)(z_1-z_3)(z_2-z_3)}=\frac{\Dl_{\{1,2\},\{3\}}}{(z_1-z_2)},
\]
\[\thi_{2,3}(\Dl_{[1]}\wedge\Dl_{[3]})=\frac{(\gm_{1,1}-z_2)(\gm_{1,2}-z_2)}{(z_1-z_3)(z_1-z_2)(z_3-z_2)}=\frac{\Dl_{\{1,3\},\{2\}}}{(z_1-z_3)},
\]
\[\thi_{2,3}(\Dl_{[2]}\wedge\Dl_{[3]})\frac{(\gm_{1,1}-z_1)(\gm_{1,2}-z_1)}{(z_2-z_3)(z_2-z_1)(z_3-z_1)}=\frac{\Dl_{\{2,3\},\{1\}}}{(z_2-z_3)}.
\]
\qetr
\end{example}

The following result is due to several people \cite{Gat05,GS09,GS10,Lak08,Lak09,LT09,GM11}. Moreover, it also follows by the application of more general theories on the geometry of affine Grassmannians, perverse sheaves on loop groups, and Langlands duality \cite{Lus83,BD96,Gin90,Gin95,MV00,MV07}. See also \cite[Sec.\,7.5]{CDG1} and \cite{AN20}. Our results in Section \ref{sec2} and Section \ref{sec5} provide a simple direct proof.
\begin{thm}\label{gSc}
The operators of quantum multiplication $c_1^{T}(E_i^G)*_{q,\bm z}$, with $i=1,2$, on $H^\bullet_T(G(k,n),\C)$ are related to the quantum multiplication operators $c_1^{T}(E_i^\Pb)*_{(-1)^{k-1}q,\bm z}$, with $i=1,2$, 
on $H^\bullet_T(\Pb^{n-1},\C)$, by the following formulas:
\[c_1^{T}(E_1^G)*_{q,\bm z}\,\si_{\bm \la}^G=\thi_{k,n}\left(\sum_{j=1}^k \si_{\la_k}^\Pb\wedge\dots\wedge c_1^{T}(E_1^\Pb)*
_{(-1)^{k-1}q,\bm z}\si_{\la_j+k-j}^\Pb\wedge\dots\wedge\si_{\la_1+k-1}^\Pb\right),\quad  
\]
\begin{multline*}
c_1^{T}(E_2^G)*_{q,\bm z}\,\si_{\bm \la}^G=\thi_{k,n}\left(\sum_{j=1}^k \si_{\la_k}^\Pb\wedge\dots\wedge c_1^{T}(E_2^\Pb)*
_{(-1)^{k-1}q,\bm z}\si_{\la_j+k-j}^\Pb\wedge\dots\wedge\si_{\la_1+k-1}^\Pb\right)\\
+(1-k)\left(\sum_{a=1}^nz_a\right)\si_{\bm\la}^G,
\end{multline*}for any $\bm\la\subseteq k\times (n-k)$.
\end{thm}
\proof
The identities immediately follow from Theorems \ref{satdyn} and \ref{stabint}.
\endproof
\begin{example}
Let $k=2$ and $n=4$. The Schubert basis of $H^\bullet_T(G(2,4),\C)$ is 
\[\si_0^G,\quad \si_{1}^G,\quad \si_{2}^G,\quad \si_{1,1}^G,\quad \si_{2,1}^G,\quad \si_{2,2}^G.
\]
We have $c_1^{T}(E_1^G)=-\si^G_1+z_3+z_4$, by Example \ref{ex1}. We deduce the multiplication table:
\begin{align*}
c_1^{T}(E_1)*_{q,\bm z}\si_0^G&=-\si^G_1+z_3+z_4,\\
c_1^{T}(E_1)*_{q,\bm z}\si_1^G&=-\si_2^G-\si_{1,1}^G+(z_2+z_4)\si^G_1,\\
c_1^{T}(E_1)*_{q,\bm z}\si_2^G&=-\si_{2,1}^G+(z_1+z_4)\si_2^G,\\
c_1^{T}(E_1)*_{q,\bm z}\si_{1,1}^G&=-\si_{2,1}^G+(z_2+z_3)\si_{1,1}^G,\\
c_1^{T}(E_1)*_{q,\bm z}\si_{2,1}^G&=-\si_{2,2}^G+(z_1+z_3)\si_{2,1}^G-q,\\
c_1^{T}(E_1)*_{q,\bm z}\si_{2,2}^G&=(z_1+z_2)\si_{2,2}^G-q\si_1^G.
\end{align*}
The table above can be computed from the multiplication table in \cite[Sec.\,8.2, pag.\,27]{Mih06}. (Notice that the equivariant parameters $T_i$ of {\it loc.\,cit.\,\,}are the {\it opposite} of our parameters $z_i$, that is $T_i=-z_i$ for $i=1,\dots, n$. See also \cite[Sec.\,5.1, pag.\,2299]{Mih08}.)
\vskip2mm
The above formulas can actually be recovered from the quantum multiplication $c_1^{T}(\mc O(-1))*_{q,\bm z}$ on $H^\bullet_T(\Pb^3,\C)$.
The Schubert basis for $H^\bullet_T(\Pb^3,\C)$ is (for short, we omit the label $\Pb$)
\[\si_0,\quad \si_{1},\quad \si_{2},\quad \si_{3}.
\]For short, set $x:=c_1^{T}(\mc O(-1))=-\si_1+z_4$. In the Schubert basis, the operator $x*_{q,\bm z}\in\End\left(H^\bullet_T(\Pb^3,\C)\right)$ is represented by the matrix
\[x*_{q,\bm z}=\begin{pmatrix}
z_4&0&0&-q\\
-1&z_3&0&0\\
0&-1&z_2&0\\
0&0&-1&z_1
\end{pmatrix}.
\]
We identify $H^\bullet_T(G(2,4),\C)$ with $\bigwedge\nolimits^2 H^\bullet_T(\Pb^3,\C)$, via the isomorphism $\thi_{2,4}^{-1}$:
\begin{align*}
\si^G_{0}&\mapsto \si_0\wedge \si_1,&
\si^G_{1}&\mapsto \si_0\wedge\si_2,&
\si^G_{2}&\mapsto \si_0\wedge\si_3,\\
\si^G_{1,1}&\mapsto \si_1\wedge\si_2,&
\si^G_{2,1}&\mapsto \si_1\wedge\si_3,&
\si^G_{2,2}&\mapsto \si_2\wedge\si_3.
\end{align*}
We have
\begin{align*}
\thi_{2,4}^{-1}(c_1^{T}(E_1)*_{q,\bm z} \si_0^G)&=& (x*_{-q,\bm z}\si_0)\wedge\si_1 +\si_0\wedge(x*_{-q,\bm z}\si_1)&\\
&=&(-\si_1+z_4)\wedge\si_1+\si_0\wedge(z_3\si_1-\si_2)&\\
&=&-\si_0\wedge\si_2+(z_3+z_4)\si_0\wedge\si_1&\\
&=&\thi_{2,4}^{-1}(-\si_1^G+z_3+z_4),&\\%
&&&\\
\thi_{2,4}^{-1}(c_1^{T}(E_1)*_{q,\bm z} \si_1^G)&=&(x*_{-q,\bm z}\si_0)\wedge \si_2 +\si_0\wedge(x*_{-q,\bm z}\si_2)&\\
&=&(-\si_1+z_4)\wedge\si_2+\si_0\wedge(z_2\si_2-\si_3)&\\
&=&-\si_0\wedge\si_3-\si_1\wedge\si_2+(z_2+z_4)\si_0\wedge\si_2&\\
&=&\thi_{2,4}^{-1}(-\si_2^G-\si_{1,1}^G+(z_2+z_4)\si_1^G),&\\%
&&&\\
\thi_{2,4}^{-1}(c_1^{T}(E_1)*_{q,\bm z} \si_2^G)&=&(x*_{-q,\bm z}\si_0)\wedge\si_3 +\si_0\wedge(x*_{-q,\bm z}\si_3)&\\
&=&(-\si_1+z_4)\wedge\si_3+\si_0\wedge(q\si_0+z_1\si_3)&\\
&=&-\si_1\wedge\si_3+(z_1+z_4)\si_0\wedge\si_3&\\
&=&\thi_{2,4}^{-1}(-\si_{2,1}^G+(z_1+z_4)\si_2^G),&\\%
&&&\\
\thi_{2,4}^{-1}(c_1^{T}(E_1)*_{q,\bm z} \si_{1,1}^G)&=&(x*_{-q,\bm z}\si_1)\wedge\si_2+\si_1\wedge(x*_{-q,\bm z}\si_2)&\\
&=&(z_3\si_1-\si_2)\wedge\si_2+\si_1\wedge(z_2\si_2-\si_3)&\\
&=&-\si_1\wedge \si_3+(z_2+z_3)\si_1\wedge\si_2&\\
&=&\thi_{2,4}^{-1}(-\si_{2,1}^G+(z_2+z_3)\si_{1,1}^G),&\\%
&&&\\
\thi_{2,4}^{-1}(c_1^{T}(E_1)*_{q,\bm z} \si_{2,1}^G)&=&(x*_{-q,\bm z}\si_1)\wedge\si_3+\si_1\wedge(x*_{-q,\bm z}\si_3)&\\
&=&(z_3\si_1-\si_2)\wedge\si_3+\si_1\wedge(q\si_0+z_1\si_3)&\\
&=&-\si_2\wedge\si_3+(z_1+z_3)\si_1\wedge\si_3+q\si_1\wedge\si_0&\\
&=&\thi_{2,4}^{-1}(-\si_{2,2}^G+(z_1+z_3)\si_{2,1}^G-q\si_0^G),&\\%
&&&\\
\thi_{2,4}^{-1}(c_1^{T}(E_1)*_{q,\bm z} \si_{2,2}^G)&=&(x*_{-q,\bm z}\si_2)\wedge\si_3+\si_2\wedge(x*_{-q,\bm z}\si_3)&\\
&=&(z_2\si_2-\si_3)\wedge\si_3+\si_2\wedge(q\si_0+z_1\si_3)&\\
&=&(z_1+z_2)\si_2\wedge\si_3+q\si_2\wedge\si_0&\\
&=&\thi_{2,4}^{-1}((z_1+z_2)\si_{2,2}^G-q\si_{1}^G)&.%
\end{align*}
\qetr
\end{example}

The isomorphisms $\theta_{k,n},\thi_{k,n}$ induce two isomorphisms of bundles 
\[\theta_{k,n}^\diamond\colon \bigwedge\nolimits^k U_{1,n}\to U_{k,n},\qquad \thi_{k,n}^\diamond\colon \bigwedge\nolimits^k H_{1,n}\to H_{k,n},
\]fitting in the following commutative diagram
\[
\xymatrix{
\bigwedge\nolimits^k U_{1,n}\ar[rr]^{\bigwedge^k{\rm Stab}^\diamond_{1,n}}\ar[d]_{\theta_{k,n}^\diamond}&&\bigwedge\nolimits^k H_{1,n}\ar[d]^{\thi_{k,n}^\diamond}\\
U_{k,n}\ar[rr]_{{\rm Stab}^\diamond_{k,n}}&&H_{k,n}}
\]
For short, denote by $\mathsf{X}_i^\Pb(\bm z;\bm p)$ and $\mathsf Z_a^\Pb(\bm z;\bm p;\kappa)$, with $i=1,2$ and $a=1,\dots,n$, the sections of $\End H_{1,n}$ defined by
\begin{align*}
&\mathsf{X}_i^\Pb(\bm z;\bm p)f(\bm z;\bm p):=c_1^T(E_i^\Pb)*_{q,\bm z}f(\bm z;\bm p),& q=\frac{p_2}{p_1},\quad i=1,2,\\
&\mathsf Z_a^\Pb(\bm z;\bm p;\kappa)f(\bm z;\bm p)=K^{H,\Pb}_a(\bm z;\bm p;\kappa)^{-1}\exp\left(\kappa\frac{\der}{\der z_a}\right)f(\bm z;\bm p),& a=1,\dots,n.
\end{align*}
Set $\tilde{\bm{p}} = ((-1)^{k-1} p_1, p_2)$, and let $\Hat{\mathsf X}_i^{(k)}, \Hat{\mathsf Z}_a^{(k)}$, with $i=1,2$ and $a=1,\dots,n$, be the sections of $\End\bigwedge^kH_{1,n}$ acting on decomposable sections of $\bigwedge^kH_{1,n}$ by
\begin{align}
&\Hat{\mathsf{X}}^{(k)}_{i}(\bm z;\bm p)\left[\bigwedge_{j=1}^m f_j(\bm{z}; \bm{p})\right] := \sum_{j=1}^m f_1(\bm{z}; \bm{p}) \wedge \dots \wedge \mathsf X^\Pb_i(\bm z;\tilde{\bm p}) f_j(\bm{z}; \bm{p}) \wedge \dots \wedge f_m(\bm{z}; \bm{p}), \\
&\Hat{\mathsf{Z}}^{(k)}_{a}(\bm z;\bm p;\kappa)\left[\bigwedge_{j=1}^m f_j(\bm{z}; \bm{p})\right] := \bigwedge_{j=1}^m \mathsf Z^\Pb_a(\bm z;\tilde{\bm p};\kappa) f_j(\bm{z}; \bm{p}).
\end{align}
The joint system of differential-difference equations 
\begin{align}
\label{eqh1nk.1}
&\kappa p_i\frac{\der}{\der p_i}f(\bm z;\bm p)=\Hat{\mathsf{X}}^{(k)}_{i}(\bm z;\bm p)f(\bm z;\bm p),& i=1,2,\\
\label{eqh1nk.2}
&\Hat{\mathsf{Z}}^{(k)}_{a}(\bm z;\bm p;\kappa)f(\bm z;\bm p)=f(\bm z;\bm p),&a=1,\dots,n,
\end{align}
where $f(\bm z;\bm p)$ is a section of $\bigwedge^kH_{1,n}$, is compatible.
\begin{thm}\label{mainthm1}
Under the isomoprhism $\thi_{k,n}^\diamond$, the qDE and qKZ system of equations \eqref{qde.0}, \eqref{qkz.H} for $G(k,n)$ is gauge equivalent to the joint system of equations \eqref{eqh1nk.1}, \eqref{eqh1nk.2}. More precisely, if $f(\bm z;\bm p)$ is a solution of \eqref{eqh1nk.1}, \eqref{eqh1nk.2}, then the function
\[G(\bm z;\bm p;\kappa)f(\bm z;\bm p),\quad \text{with }G(\bm z;\bm p;\kappa)=p_2^{(1-k)\sum_{a=1}^n z_a/\kappa},
\]is a solution of the qDE and qKZ equations \eqref{qde.0} and \eqref{qkz.H}.
\end{thm}
\proof
It is a direct consequence of Theorem \ref{abstractthmsat}.
\endproof
\begin{thm}\label{satsolpg}
Let $f_1,\dots, f_k$ be sections of $H_{1,n}$, solving the joint system \eqref{qde.0} and \eqref{qkz.H} for the projective space $\Pb^{n-1}$. The section of $H_{k,n}$ defined by
\beq\label{finalsolprod}
G(\bm z;\bm p;\kappa)\cdot \thi_{k,n}^\diamond\left[f_1(\bm z;(-1)^{k-1}p_1,p_2)\wedge\dots\wedge f_k(\bm z;(-1)^{k-1}p_1,p_2)\right],\quad G(\bm z;\bm p;\kappa)=p_2^{(1-k)\sum_{a=1}^nz_a/\kappa},
\eeq
is a solution of the joint system of equations \eqref{qde.0} and \eqref{qkz.H} for the Grassmannian $G(k,n)$.
\end{thm}
\proof
It is a direct consequence of Corollary \ref{cor1}.
\endproof
\begin{rem}
In equation \eqref{finalsolprod}, a determination of the argument of $-1$ has to be fixed. As the solutions $f_1,\dots,f_k$ are typically multivalued, the solution derived in \eqref{finalsolprod} depends on this choice.
\end{rem}

\section{{\cyr B}-Theorem, determinantal identities, and integral representations}
\label{SEC7}

\subsection{{\cyr B}-Theorem}\label{secbthm1} Recall the definition of the subset $L_\kappa\subseteq\C^n$ and of its ring of holomorphic functions $\mc O_{L_\kappa}$, see Section \ref{secext}. Consider the space $\mathscr S_{H_{k,n}}^{\mc O}$ of sections of $H_{k,n}$ solving the equivariant differential equations \eqref{qde.0}, and  holomorphic in $\bm z\in L_\kappa$. The space $\mathscr S_{H_{k,n}}^{\mc O}$ is a module over $\mc O_{L_\kappa}$.
\vskip1,5mm
By extension of scalars, the topological--enumerative morphism induces a map
\[\mc S\colon H^{(\kappa)}_T(G(k,n))\to\mathscr S_{H_{k,n}}^{\mc O},\quad a\mapsto S(\bm z;\bm p;\kappa)a.
\]
Recall, from Section \ref{secideksol}, the identification $\mu_{k,n}\colon K^T_0(G(k,n))_\C\to\mathscr S_{H_{k,n}}$ defined by $[P]\mapsto \on{Stab}_{k,n}\Psi_P$. Since every element of $\mathscr S_{H_{k,n}}$ is entire in $\bm z$, it is an element of $\mathscr S_{H_{k,n}}^{\mc O}$.

In \cite{CV21}, the following result was proved in the case of projective spaces. Subsequently, in \cite{TV23}, the result was generalized to arbitrary partial flag varieties.

\begin{thm}\label{bth}\cite[\textcyr{B}-Theorem]{CV21,TV23}
The following diagram is commutative,
\beq
\xymatrix{
K^T_0(G(k,n))_\C\ar[rr]^{\text{\textnormal{\cyr{B}}}(-;\kappa)}\ar[dr]_{\mu_{k,n}}&&H^{(\kappa)}_T(G(k,n))\ar[dl]^{\mc S}\\
&\mathscr S^{\mc O}_{H_{k,n}}&
}
\eeq
\end{thm}

In this section we prove that the validity of the \textcyr{B}-theorem for the Grassmannians $G(k,n)$ can be deduced from the analogue statement for $\Pb^{n-1}$, via Satake identification. We first briefly recall the proof as in \cite{CV21,TV23}.

\proof[Proof of Theorem \ref{bth}]
Let $P(\bm\Gm,\bm Z)\in K^T_0(G(k,n))_\C$, and let $a_P(\bm\gm,\bm z)\in H^{(\kappa)}_T(G(k,n))$ be such that 
\[\mu_{k,n}(P)=\mc S(a_P),\quad\mc S\,\,\text{topological--enumerative morphism.}
\]The right-hand side is of the form
\[\left(1+O(p_1^{-1}p_2)\right)p_1^{c_1^{T}(E_1)/\kappa}p_2^{c_1^{T}(E_2)/\kappa}a_P,\quad \text{in the regime }p_1^{-1}p_2\to0.
\]The class $p_1^{c_1^{T}(E_1)/\kappa}p_2^{c_1^{T}(E_2)/\kappa}a_P$ can be then identified with the leading term of the residue series expansions defining the left-hand side, which equals
\begin{multline*}\kappa^{-k(n-k)-k+k}\sum_{J\in\mc I_{k,n}}\sum_{I\in\mc I_{k,n}}\acute{P}(\bm z_{\si_J};\bm z;\kappa)\left(\kappa^{n-k}p_1\right)^{\sum_{a\in J_1}z_a/\kappa}\left(\kappa^{-k}p_2\right)^{\sum_{a\in J_2}z_a/\kappa}\\
\times \prod_{a\in J_1}\prod_{b\in J_2}\Gm\left(\frac{z_a-z_b}{\kappa}\right)W_I(\Si_J;\bm z)\on{Stab}_I.
\end{multline*}
Introduce the polynomial
\[\Hat W_{k,n}(\bm t;\bm \gm):=\prod_{j=1}^k\prod_{i=1}^{n-k}(t_j-\gm_{2,i}),
\]and let $[\Hat W_{k,n}(\bm t;\bm \gm)]$ be its image in $\C[\bm t]\otimes H^\bullet_T(G(k,n),\C)$. We have
\[[\Hat W_{k,n}(\bm t;\bm \gm)]=\sum_{I\in\mc I_{k,n}}W_I(\bm t;\bm z)\on{Stab}_I,
\]and also
\[\Hat W_{k,n}(\Si_J;\bm z_{\si_I})=\dl_{I,J}R(\bm z_{\si_J}),\quad R(\bm z):=\prod_{a=1}^k\prod_{b=1}^{n-k}(z_a-z_b).
\]For any $J\in\mc I_{k,n}$, we then have
\beq\label{apj}
a_P|_{\on{pt}_J}=\sum_{J\in\mc I_{k,n}}\acute{P}(\bm z_{\si_J};\bm z;\kappa)\kappa^{(n-k)\sum_{a\in J_1}z_a/\kappa+k\sum_{b\in J_2}z_b/\kappa}\prod_{a\in J_1}\prod_{b\in J_2}\Gm\left(1+\frac{z_a-z_b}{\kappa}\right).
\eeq
We claim that this equals the restriction \textcyr{B}$(P;\kappa)|_{\on{pt}_J}$. Indeed, we have $TG(k,n)\cong E_1^*\otimes E_2$, so that the equivariant Chern roots of $TG(k,n)$ are $\{-\gm_{1,i}+\gm_{2,j}\colon i=1,\dots,k,\,j=1,\dots,n-k\}$ and also
\[c_1^{T}(G(k,n))=\on{rk}(E_1^*)c_1^{T}(E_2)+\on{rk}(E_2)c_1^{T}(E_1^*)=-(n-k)\sum_{i=1}^k\gm_{1,i}+k\sum_{j=1}^{n-k}\gm_{2,j}.
\]From this, it is easily seen that \textcyr{B}$(P;\kappa)|_{\on{pt}_J}$ equals the right-hand side of \eqref{apj}. Hence $a_P=\textcyr{B}(P;\kappa)$. This proves the theorem.
\endproof

\subsection{{\cyr B}-morphisms of Grassmannians and Satake identification}

For short, denote by $\gm=\gm_{1,1}$ the equivariant first Chern root be of the tautological bundle $\mc O(-1)$ on $\Pb^{n-1}$. The equivariant Chern roots of the tangent bundle ${T}\Pb^{n-1}$ are $(-\gm+z_a)_{a=1}^n$. The \textcyr{B}-morphism is
\begin{multline*}
\text{\textcyr{B}}^{(\kappa)}_\Pb\colon K^T_0(\Pb^{n-1})_\C\to H^{(\kappa)}_T(\Pb^{n-1},\C),\\ V\mapsto \exp\left[\frac{\log\kappa
}{\kappa}\left(n\gm-\sum_{a=1}^nz_a\right)\right]\prod_{a=1}^n\Gm\left(1+\frac{\gm-z_a}{\kappa}\right)\on{Ch}^{(\kappa)}(V).
\end{multline*}
This induces a morphism $\bigwedge^k\text{\textcyr{B}}^{(\kappa)}_\Pb\colon \bigwedge^kK^T_0(\Pb^{n-1})_\C\to\bigwedge^kH^{(\kappa)}_T(\Pb^{n-1})$.

\begin{prop}\label{propsatbmor}
We have
\[\textnormal{\textcyr B}^{(\kappa)}_{G}(\Si_{\bm\la}E_1^*)=\left(\frac{2\pi\sqrt{-1}}{\kappa}\right)^{-\binom{k}{2}}\exp[\frac{\pi\sqrt{-1}}{\kappa}c_1^{T}(\wedge^2 E_1^*)]\,\,\thi_{k,n}\left(\bigwedge_{j=1}^k\textnormal{\textcyr B}^{(\kappa)}_{\Pb}(\mc O(\la_j+k-j))\right).
\]
\end{prop}
\proof
Recall that $TG(k,n)\cong E_1^*\otm E_2$, so that in $K_0^{T}(G(k,n))$ we have $[TG(k,n)]=[E_1^*\otm \underline{\C^n}]-[E_1^*\otm E_1]$. This is easily seen by tensoring the Euler sequence \eqref{Euler} by $E_1^*$. The equivariant Chern roots of $E_1^*\otm\underline{\C^n}$ are $\{-\gm_{1,i}+z_a,\,i=1,\dots,k,\,a=1,\dots,n\}$, and the equivariant Chern roots of $E_1^*\otm E_1$ are $\{-\gm_{1,i}+\gm_{1,j},\,i,j=1,\dots,k\}$. So, we have
\[
\exp\left[-\frac{\log\kappa
}{\kappa}c_1(G(k,n))\right]\Hat\Gm^{(\kappa)}_{G(k,n)}=\frac{\prod_{i=1}^k\prod_{a=1}^n\Gm\left(1+\frac{\gm_{1,i}-z_a}{\kappa}\right)\exp\left[\frac{\log\kappa
}{\kappa}(\gm_{1,i}-z_a)\right]}{\prod_{i,j=1}^k\Gm\left(1+\frac{\gm_{1,i}-\gm_{1,j}}{\kappa}\right)\exp\left[\frac{\log\kappa
}{\kappa}(\gm_{1,i}-\gm_{1,j})\right]}.
\]Notice that
\begin{multline*}
\prod_{i,j=1}^k\Gm\left(1+\frac{\gm_{1,i}-\gm_{1,j}}{\kappa}\right)\exp\left[\frac{\log\kappa
}{\kappa}(\gm_{1,i}-\gm_{1,j})\right]=\prod_{i<j}\Gm\left(1+\frac{\gm_{1,i}-\gm_{1,j}}{\kappa}\right)\Gm\left(1+\frac{\gm_{1,j}-\gm_{1,i}}{\kappa}\right)\\
=\prod_{i<j}\frac{2\pi\sqrt{-1}}{\kappa}\frac{\gm_{1,i}-\gm_{1,j}}{\exp[\frac{\pi\sqrt{-1}}{\kappa}(\gm_{1,i}-\gm_{1,j})]-\exp[\frac{\pi\sqrt{-1}}{\kappa}(\gm_{1,j}-\gm_{1,i})]}\\
=\left(\frac{2\pi\sqrt{-1}}{\kappa}\right)^{\binom{k}{2}}\prod_{i<j}(\gm_{1,i}-\gm_{1,j})\prod_{i<j}\frac{\exp[\frac{\pi\sqrt{-1}}{\kappa}(\gm_{1,i}+\gm_{1,j})]}{\exp[\frac{2\pi\sqrt{-1}}{\kappa}\gm_{1,i}]-\exp[\frac{2\pi\sqrt{-1}}{\kappa}\gm_{1,j}]}\\
=\left(\frac{2\pi\sqrt{-1}}{\kappa}\right)^{\binom{k}{2}}\exp[-\frac{\pi\sqrt{-1}}{\kappa}c_1^{T}(\wedge^2 E_1^*)]\prod_{i<j}\frac{\gm_{1,i}-\gm_{1,j}}{\exp[\frac{2\pi\sqrt{-1}}{\kappa}\gm_{1,i}]-\exp[\frac{2\pi\sqrt{-1}}{\kappa}\gm_{1,j}]}.
\end{multline*}
For the last equality, we used the fact that $\{-\gm_{1,i}-\gm_{1,j}\}_{i<j}$ are the equivariant Chern roots of the bundle $\wedge^2E_1^*$. Then we have
\begin{multline*}\exp\left[-\frac{\log\kappa
}{\kappa}c_1(G(k,n))\right]\Hat\Gm^{(\kappa)}_{G(k,n)}=\\
\left(\frac{2\pi\sqrt{-1}}{\kappa}\right)^{-\binom{k}{2}}\exp[\frac{\pi\sqrt{-1}}{\kappa}c_1^{T}(\wedge^2 E_1^*)]\,\,\prod_{i<j}\frac{\exp[\frac{2\pi\sqrt{-1}}{\kappa}\gm_{1,i}]-\exp[\frac{2\pi\sqrt{-1}}{\kappa}\gm_{1,j}]}{\gm_{1,i}-\gm_{1,j}}\\
\cdot \prod_{i=1}^k\prod_{a=1}^n\Gm\left(1+\frac{\gm_{1,i}-z_a}{\kappa}\right)\exp\left[\frac{\log\kappa
}{\kappa}(\gm_{1,i}-z_a)\right].
\end{multline*}
Since the Chern character $\on{Ch}^{(\kappa)}_T$ defines a morphism of rings, we have
\[\on{Ch}^{(\kappa)}_T(\Si_{\bm\la}(E_1^*))=\on{Ch}^{(\kappa)}_T(s_{\bm\la}(\bm\Gm_1^{-1}))=s_{\bm\la}(\acute{\bm\gm_1}^{-1})=\frac{\det[\exp[\frac{2\pi\sqrt{-1}}{\kappa}\gm_{1,i}(\la_j+k-j)]]_{i,j}}{\prod_{i<j}\left(\exp[\frac{2\pi\sqrt{-1}}{\kappa}\gm_{1,i}]-\exp[\frac{2\pi\sqrt{-1}}{\kappa}\gm_{1,j}]\right)}.
\]
This proves the claim.
\endproof

\subsection{Isometric identification of $K$-theories}
Define the equivariant line bundle $\mathscr L$ on $G(k,n)$ as $\mathscr L:=\det \bigwedge^2E_1^*$. 
\begin{lem}
The equivariant $K$-theoretical class $[\mathscr L]$ equals the Laurent polynomial $\left(\prod_{i=1}^k\Gm_{1,i}\right)^{1-k}$.
\end{lem}
\proof
The equivariant $K$-theoretical Chern roots of $\wedge^2E_1^*$ are $\left(\Gm_{1,i}\Gm_{1,j}\right)^{-1}$ for all $i,j=1,\dots,k$ with $i<j$. The class $[\mathscr L]$ equals the product of these Chern roots.
\endproof
Introduce the 
$\C[\bm Z^{\pm 1}]$-linear map 
$\vartheta^{K
}_{k,n}\colon \bigwedge^k K^T_0(\Pb^{n-1})_\C\to K^T_0(G(k,n))_\C$ by
\begin{align}
\vartheta^{K
}_{k,n}[\mc O(\la_1+k-1)\wedge\mc O(\la_2+k-2)\wedge\dots\wedge\mc O(\la_k)]&=\mathscr L\otimes \Si_{\bm\la}E_1^*,
\end{align}
for any $\bm\la\subseteq k\times(n-k)$. In terms of the polynomials presentation \eqref{kthpres}, the map $\vartheta^{K}_{k,n}$ can be described as follows. For short, set $X:=\Gm_{1,1}^\Pb$ and let $\bm\Gm=(\bm\Gm_1,\bm\Gm_2)$ be the $K$-theoretical Chern roots of $E_1^G,E_2^G$, respectively. The map $\thi^K_{k,n}$ is then given by 
\begin{multline*}
\thi_{k,n}^{K
}\colon \bigwedge\nolimits^k_{\C}\left(\C[X^{\pm 1},\bm Z^{\pm 1}]\Big/\Big\langle\prod_{i=1}^n(X-Z_i)\Big\rangle\right)\\ \to \C[\bm\Gm^{\pm 1}]^{S_{k}\times S_{n-k}}[\bm Z^{\pm 1}]\Big/\Big\langle   e_j(\bm Z)-\sum_{h=0}^j e_h(\bm\Gm_1)e_{j-h}(\bm\Gm_2)\Big\rangle_{j=1}^n,
\end{multline*}
\begin{align}
\label{forthi+}
\vartheta_{k,n}^{K
}(f_1(X;\bm Z)\wedge\dots\wedge f_k(X;\bm Z)):= \frac{\det
\left(\begin{array}{ccc}f_1(\Gm_{1,1};\bm Z)&\dots&f_k(\Gm_{1,1};\bm Z)\\
\vdots&&\vdots\\
f_1(\Gm_{1,k};\bm Z)&\dots&f_k(\Gm_{1,k};\bm Z)
\end{array}\right)}{\left(\prod_{i=1}^k\Gm_{1,i}\right)^{k-1}\frak D_k(\bm\Gm^{-1})},
\end{align}
where $\frak D_k$ denotes the Vandermonde determinant as in \eqref{vandet}.
\vskip2mm
\begin{rem}
The denominator of \eqref{forthi+} equals
\[\left(\prod_{i=1}^k\Gm_{1,i}\right)^{k-1}\frak D_k(\bm\Gm^{-1})=\det\begin{pmatrix}
1&\Gm_{1,1}&\dots&(\Gm_{1,1})^{k-1}\\
1&\Gm_{1,2}&\dots&(\Gm_{1,2})^{k-1}\\
\vdots&&&\vdots\\
1&\Gm_{1,k}&\dots&(\Gm_{1,k})^{k-1}
\end{pmatrix}=(-1)^{\binom{k}{2}}\frak D_k(\bm \Gm).
\]
\end{rem}
The equivariant Grothendieck--Euler--Poincar\'e pairing $\chi^\Pb$ on $K^T_0(\Pb^{n-1})_\C$ induces a pairing $\chi^{\wedge\Pb}$ on the exterior product $\bigwedge\nolimits^k_{\C[\bm Z^{\pm 1}]} K^T_0(\Pb^{n-1})_\C$, defined by
\[\chi^{\wedge\Pb}(V_1\wedge\dots\wedge V_k,W_1\wedge\dots\wedge W_k)=\det\left(\chi^\Pb(V_i,W_j)\right)_{i,j=1}^k.
\]

\begin{prop}\label{propisomkth}
The map $\vartheta_{k,n}^{K}$ is an isomorphism 
of $\C[\bm Z^{\pm 1}]$-modules. Moreover, it preserves 
the Grothendieck--Euler--Poincar\'e pairings, that is 
\beq
\label{kisom}
\chi^{\wedge\Pb}(\al_1,\al_2)=\chi^G\left(\vartheta_{k,n}^{K}(\al_1),\vartheta_{k,n}^{K}(\al_2)\right), 
\quad \al_1,\al_2\in \bigwedge\nolimits^k_{\C[\bm Z^{\pm 1}]} K^T_0(\Pb^{n-1})_\C.
\eeq
\end{prop}
\proof The collections $\left(\bigwedge_{j=1}^k\mc O(\la_j+k-j)\right)_{\bm\la\subseteq k\times(n-k)}$ and $\left(\Si_{\bm\la}E_1^*\right)_{\bm\la\subseteq k\times(n-k)}$ are bases of the $\C[\bm Z^{\pm 1}]$-modules $\bigwedge\nolimits^k_{\C[\bm Z^{\pm 1}]} K^T_0(\Pb^{n-1})_\C$ and $K^T_0(G(k,n))_\C$, respectively. See Corollary \ref{kapfullcoll}. Hence, the morphism $\vartheta_{k,n}^{K}$ is an isomorphism 
 by definition.
Given $\bm\la^{1},\bm\la^2\subseteq k\times (n-k)$, by Theorem \ref{HRR2} and Proposition \ref{isometry1}, we have
\begin{align*}
&\chi^G(\Si_{\bm\la^1}E_1^*,\Si_{\bm\la^2}E_1^*)=\left(-\frac{\kappa}{2\pi\sqrt{-1}}\right)^{k(n-k)}\int_{G(k,n)}^{\rm eq}\textnormal{\textcyr{B}}(\Si_{\bm\la^1}E_1^*;e^{-\pi\sqrt{-1}}\kappa)\textnormal{\textcyr{B}}(\Si_{\bm\la^2}E_1^*;\kappa)\\
&=\left(-\frac{\kappa}{2\pi\sqrt{-1}}\right)^{k(n-k)+2\binom{k}{2}}(-1)^{\binom{k}{2}}\,\,\times\\
&\qquad\qquad\int^{\rm eq}_{G(k,n)}\thi_{k,n}\left(\bigwedge_{j=1}^k\textnormal{\textcyr{B}}\left(\mc O(\la^1_j+k-j);e^{-\pi\sqrt{-1}}\kappa\right)\right)\thi_{k,n}\left(\bigwedge_{h=1}^k\textnormal{\textcyr{B}}\left(\mc O(\la^2_h+k-h);\kappa\right)\right)\\
&=\det\left(\left(-\frac{\kappa}{2\pi\sqrt{-1}}\right)^{n-1}\int^{\rm eq}_{\Pb^{n-1}}\textnormal{\textcyr{B}}\left(\mc O(\la^1_j+k-j);e^{-\pi\sqrt{-1}}\kappa\right)\textnormal{\textcyr{B}}\left(\mc O(\la^2_h+k-h);\kappa\right)\right)_{j,h=1}^k\\
&=\det\left(\chi^\Pb(\mc O(\la^1_j+k-j),\mc O(\la^2_h+k-h))\right)_{j,h=1}^k.
\end{align*}
This proves the claim, since tensoring by a line bundle preserves the Grothendieck--Euler--Poincar\'e pairing.
\endproof

Let $(e_1,\dots e_n)$ be a basis of $K^T_0(\Pb^{n-1})$. For each $I\in\mc I_{k,n}$, with $I=(I_1,I_2)$ and $I_1=\{i_1<\dots<i_k\}$, set $e_I:=e_{i_1}\wedge\dots\wedge e_{i_k}$. The tuple $\left(\thi^K_{k,n}e_I\right)_{I\in\mc I_{k,n}}$ is then a basis of $K^T_0(G(k,n))$.

\begin{cor}\label{corkexcbases}
If $(e_1,\dots e_n)$ is an exceptional basis of $K^T_0(\Pb^{n-1})$, then the basis $\left(\thi^K_{k,n}e_I\right)_{I\in\mc I_{k,n}}$, lexicographically ordered, is an exceptional basis of $K^T_0(G(k,n))$.\qed
\end{cor}

The following result, used later, describes how the braid group orbits of exceptional bases in $K$-theory relate to the Satake identification. Recall also the more general action of the signed braid group $\Z^{\binom{n}{k}}\rtimes\mc B_{\binom{n}{k}}$ on the set of exceptional collections, see Remark \ref{remsignedbraid}.

\begin{prop}[{\cite[Prop.\,7.4]{CDG1}}]\label{propsatbr}
 If $(e_i)_{i=1}^n$ and $(f_i)_{i=1}^n$ are two exceptional bases of $K^T_0(\Pb^{n-1})$ related by the action of a braid in $\mc B_n$, then
the exceptional bases $\left(\thi^K_{k,n}e_I\right)_{I\in\mc I_{k,n}}$ and $\left(\thi^K_{k,n}f_I\right)_{I\in\mc I_{k,n}}$ of $K^T_0(G(k,n))$, lexicographically ordered, are in the same orbit with respect to the action of signed braids in $\Z^{\binom{n}{k}}\rtimes\mc B_{\binom{n}{k}}$.\qed
\end{prop}

\subsection{Topological--enumerative morphism and Satake identification}
By extension of scalars, the topological--enumerative morphism $\mc S_\Pb(\bm z;\bm p;\kappa)$ can act on elements of $H_T^{(\kappa)}(\Pb^{n-1})$. This induces a morphism 
\[\bigwedge\nolimits^k\mc S_\Pb(\bm z;\bm p;\kappa)\colon\bigwedge\nolimits^kH_T^{(\kappa)}(\Pb^{n-1})\to \bigwedge\nolimits^kH_T^{(\kappa)}(\Pb^{n-1}),\quad k=1,\dots,n.
\]
In what follows, set
\beq\label{05.02.24-1}
\mc S_{\Pb,k}(\bm z;\bm p;\kappa):=\mc S_\Pb\left(\bm z; e^{\pi\sqrt{-1}(k-1)}p_1,p_2;\kappa\right).
\eeq
\begin{prop}\label{propsattopsol}
The Satake identification $\thi_{k,n}\colon \bigwedge\nolimits^kH_T^{(\kappa)}(\Pb^{n-1})\to H^{(\kappa)}_T(G(k,n))$ intertwines the operator 
\beq\label{sattop1}
p_2^{(1-k)\sum_{a=1}^nz_a/\kappa}\bigwedge\nolimits^k\left[\mc S_{\Pb,k}(\bm z;\bm p;\kappa)\circ\exp[-\frac{\pi\sqrt{-1}}{\kappa}(k-1)c_1^{T}(E_1^\Pb)]\right]
\eeq
and the topological--enumerative morphism $\mc S_G(\bm z;\bm p;\kappa)$. That is,
\begin{multline}\label{sattop2}
\thi_{k,n}\circ\left(p_2^{(1-k)\sum_{a=1}^nz_a/\kappa}\bigwedge\nolimits^k\left[\mc S_{\Pb,k}(\bm z;\bm p;\kappa)\circ\exp[-\frac{\pi\sqrt{-1}}{\kappa}(k-1)c_1^{T}(E_1^\Pb)]\right]\right)=\mc S_G(\bm z;\bm p;\kappa)\circ\thi_{k,n}.
\end{multline}
\end{prop}
\proof
Via the Satake identification, the operator \eqref{sattop1} defines a solution of the qDE of $G(k,n)$, by Theorem \ref{satsolpg}. Consider the idempotent basis $(\Dl_{[i]})_{i=1}^n$ of $H^{(\kappa)}_T(\Pb^{n-1})$, and introduce the following matrices: let
\begin{itemize}
\item $\mc Z_1^\Pb,\mc Z_2^\Pb$, as in \eqref{Zmatrices}, be the diagonal matrices representing the endomorphisms of $H^{(\kappa)}_T(\Pb^{n-1})$ given by
\[a\mapsto c_1^{T}(E_1^\Pb)a,\quad a\mapsto c_1^{T}(E_2^\Pb)a,\quad\text{respectively;}
\]
\item $R$ 
be the diagonal matrix representing the endomorphism of $H^{(\kappa)}_T(\Pb^{n-1})$ given by
\[a\mapsto \exp[\frac{\pi\sqrt{-1}}{\kappa}(k-1)c_1^{T}(E_1^\Pb)]a
\] 
\end{itemize}
In the idempotent basis $(\Dl_{[i]})_{i=1}^n$, the morphism $\mc S_{\Pb,k}(\bm z;\bm p;\kappa)$ 
is represented by a matrix of the form
\[F(\bm z;\bm p;\kappa)(p_1^{\mc Z_1^\Pb}p_2^{\mc Z_2^\Pb})^\frac{1}{\kappa}R, 
\quad F(\bm z; \bm p;\kappa)={\bf 1}+\sum_{k=1}^\infty F_k(\bm z; \bm p;\kappa)(p_1^{-1}p_2)^k.
\] Then, the matrix representing the morphism \eqref{sattop1} in the basis $(\Dl_{[i_1]}\wedge\dots\wedge\Dl_{[i_k]})_{i_1<\dots<i_k}$ equals
\[p_2^{(1-k)\sum_{a=1}^nz_a/\kappa}\bigwedge\nolimits^kF(\bm z;\bm p;\kappa)\cdot \bigwedge\nolimits^k(p_1^{\mc Z_1^\Pb}p_2^{\mc Z_2^\Pb})^\frac{1}{\kappa}.
\]We have
\[\thi_{k,n}\circ \left[p_2^{(1-k)\sum_{a=1}^nz_a/\kappa}\bigwedge\nolimits^k(p_1^{\mc Z_1^\Pb}p_2^{\mc Z_2^\Pb})^\frac{1}{\kappa}\right]=p_1^{\mc Z_1^G}p_2^{\mc Z_2^G}\circ\thi_{k,n},
\]by Theorem \ref{gSc}. Moreover, the matrix $\bigwedge\nolimits^kF(\bm z;\bm p;\kappa)$ has the form ${\bf 1}+ O(p_1^{-1}p_2)$. 

Let $D$ be the diagonal matrix with entries $D_{II}=\frak D_k(\bm z_{\si_I})$, $I\in\mc I_{k,n}$. Then, the matrix representing the morphism 
\[\thi_{k,n}\circ\left(p_2^{(1-k)\sum_{a=1}^nz_a/\kappa}\bigwedge\nolimits^k\left[\mc S_{\Pb,k}(\bm z;\bm p;\kappa)\circ\exp[-\frac{\pi\sqrt{-1}}{\kappa}(k-1)c_1^{T}(E_1^\Pb)]\right]\right)\circ\thi_{k,n}^{-1}
\]
in the idempotent basis $(\Dl_I)_{I\in\mc I_{k,n}}$ equals
\[D\left(\bigwedge\nolimits^k F(\bm z;\bm p;\kappa)\right)D^{-1}\cdot p_1^{\mc Z_1^G}p_2^{\mc Z_2^G}=\left({\bf 1}+O(p_1^{-1}p_2)\right)p_1^{\mc Z_1^G}p_2^{\mc Z_2^G}.
\]This follows by Proposition \ref{propidemp}. The uniqueness statement asserted in Theorem \ref{thmlev}, together with Theorem \ref{sollev}, 
implies identity \eqref{sattop2}.
\endproof

\subsection{{\cyr{B}}-theorem for Grassmannians and Satake identification}
Define the space $\mathscr S^{\mc O}_{H_{1,n}}(k)$ 
of sections of $H_{1,n}$ by
\beq
\mathscr S^{\mc O}_{H_{1,n}}(k):=\{f(\bm z;\bm p;\kappa)\colon f(\bm z;\bm p;\kappa)=g(\bm z;e^{\pi\sqrt{-1}(k-1)}p_1,p_2;\kappa),\quad g\in\mathscr S^{\mc O}_{H_{1,n}}\},
\eeq
and define the map $p_{k,n}\colon \mathscr S_{H_{1,n}}^{\mc O}\to \mathscr S^{\mc O}_{H_{1,n}}(k)$ 
by 
\beq
p_{k,n}(f)(\bm z;\bm p;\kappa)=f\left(\bm z;e^{\pi\sqrt{-1}(k-1)}p_1,p_2;\kappa\right).
\eeq
We have $p_{k,n}\circ \mc S_\Pb(\bm z;\bm p;\kappa)=\mc S_{\Pb,k}(\bm z;\bm p;\kappa)$. 
From the \textcyr{B}-theorem for $\Pb^{n-1}$, we deduce the commutativity of the diagram
\beq\label{commdiagpn1}
\xymatrix{
K^T_0(\Pb^{n-1})_\C\ar[rr]^{A_1}\ar[dr]_{A_3}&&H^{(\kappa)}_T(\Pb^{n-1})\ar[dl]^{A_2}\\
&\mathscr S^{\mc O}_{H_{1,n}}(k)&
}
\eeq where
\begin{align*}
&A_1(V)=\text{\textnormal{\textcyr{B}}}_\Pb(V;\kappa)\exp[\frac{\pi\sqrt{-1}}{\kappa}(k-1)c_1^{T}(E_1^\Pb)],\quad &V\in K^T_0(\Pb^{n-1})_\C,\\ 
&A_2(v)=\mc S_{\Pb,k}(\bm z;\bm p;\kappa)\exp[-\frac{\pi\sqrt{-1}}{\kappa}(k-1)c_1^{T}(E_1^\Pb)]v,\quad &v\in H^{(\kappa)}_T(\Pb^{n-1}),\\
&A_3(V)=p_{k,n}[\mu_{1,n}(V)],\quad &V\in K^T_0(\Pb^{n-1})_\C.
\end{align*}
Define $\Hat{\mathscr{S}}_{k,n}$ 
as the space of sections of $\bigwedge^kH_{1,n}$ of the form
\[(\bm z;\bm p)\mapsto p_2^{(1-k)\sum_{a=1}^nz_a/\kappa}\,f_1(\bm z;\bm p;\kappa)\wedge\dots\wedge f_k(\bm z;\bm p;\kappa),\quad f_i\in\mathscr S^O_{H_{1,n}}(k).
\]From the commutative diagram \eqref{commdiagpn1}, we obtain the commutativity of the following diagram
\beq\label{wedgecommdiagpn1}
\xymatrix{
\bigwedge\nolimits^kK^T_0(\Pb^{n-1})_\C\ar[rr]^{B_1}\ar[dr]_{B_3}&&\bigwedge\nolimits^kH^{(\kappa)}_T(\Pb^{n-1})\ar[dl]^{B_2}\\
&\Hat{\mathscr{S}}_{k,n}&
}
\eeq
where the maps $B_1,B_2,B_3$ 
 are 
\begin{align*}
&B_1(V_{1}\wedge\dots\wedge V_k)=\left(\frac{2\pi\sqrt{-1}}{\kappa}\right)^{-\binom{k}{2}}\left(A_1(V_1)\wedge\dots\wedge A_1(V_k)\right),\\
&B_2(v_1\wedge\dots\wedge v_k)=p_2^{(1-k)\sum_{a=1}^nz_a/\kappa} \left(A_2(v_1)\wedge\dots\wedge A_2(v_k)\right),\\
&B_3(V_{1}\wedge\dots\wedge V_k)=\left(\frac{2\pi\sqrt{-1}}{\kappa}\right)^{-\binom{k}{2}}p_2^{(1-k)\sum_{a=1}^nz_a/\kappa} \left(A_3(V_1)\wedge\dots\wedge A_3(V_k)\right),
\end{align*}
for any $V_1,\dots, V_k\in K^T_0(\Pb^{n-1})_\C$, and $v_1,\dots,v_k\in H^{(\kappa)}_T(\Pb^{n-1})$.
\vskip2mm

\begin{thm}\label{thmcompatibility}
We have a commutative diagram: 
\beq\label{commprisms}
\xymatrix{
\bigwedge\nolimits^kK^T_0(\Pb^{n-1})_\C\ar[rr]^{B_1}\ar[dr]_{B_3}\ar[dd]_{\thi_{k,n}^{K}}&&\bigwedge\nolimits^kH^{(\kappa)}_T(\Pb^{n-1})\ar[dl]^{B_2}\ar[dd]^{\thi_{k,n}}\\
&\Hat{\mathscr{S}}_{k,n}\ar[dd]^<<<<<<<<{\thi_{k,n}^\diamond\circ}&\\
K^T_0(G(k,n))_\C\ar[dr]_{\mu_{k,n}}\ar[rr]|\hole^<<<<<<<{\textnormal{\textcyr{B}}_{G}^{(\kappa)}}&&H^{(\kappa)}_T(G(k,n))\ar[dl]^{\mc S_G}\\
&\mathscr S^{\mc O}_{H_{k,n}}&
}
\eeq
\end{thm}
\proof
The commutativity of the top base face follows from the \textcyr{B}-theorem for $\Pb^{n-1}$.

The commutativity of the bottom base face is the \textcyr{B}-theorem for $G(k,n)$.

The commutativity of the side face
\[
\xymatrix{
\bigwedge\nolimits^kH^{(\kappa)}_T(\Pb^{n-1})\ar[d]_{\vartheta_{k,n}}\ar[r]^<<<<<{B_2}&\Hat{\mathscr{S}}_{k,n}\ar[d]^{\vartheta^\diamond_{k,n}\circ}\\
H^{(\kappa)}_T(G(k,n))\ar[r]^<<<<<{\mc S_G}&\mathscr S^{\mc O}_{H_{k,n}}
}
\]
is equivalent to Proposition 
\ref{propsattopsol}. The commutativity of the side face
\[
\xymatrix{
\bigwedge\nolimits^kK^T_0(\Pb^{n-1})_\C\ar[r]^<<<<<{B_1}\ar[d]_{\vartheta^{K}_{k,n}}& 
\bigwedge\nolimits^kH^{(\kappa)}_T(\Pb^{n-1})\ar[d]^{\vartheta_{k,n}}\\
K^T_0(G(k,n))_\C\ar[r]^{\textnormal{\textcyr{B}}_{G}^{(\kappa)}}&
H^{(\kappa)}_T(G(k,n))
}\]
is equivalent to Proposition \ref{propsatbmor}. This is easily checked, by invoking the identities $$c_1^{T}(\wedge^2E_1^*)=({\rm rk}\,E_1-1)c_1^{T}(E_1^*)=(1-k)c_1^{T}(E_1),\qquad {\rm Ch}^{(\kappa)}_T(\mathscr L)=\exp\left[\frac{2\pi\sqrt{-1}}{\kappa}(k-1)c_1^{T}(E_1)\right].$$
The commutativity of the remaining side face
\beq\label{finaldiag}
\xymatrix{
\bigwedge\nolimits^kK^T_0(\Pb^{n-1})_\C\ar[r]^<<<<<{B_3}\ar[d]_{\vartheta^{K}_{k,n}}& 
\Hat{\mathscr{S}}_{k,n}\ar[d]^{\vartheta^\diamond_{k,n}\circ}\\
K^T_0(G(k,n))_\C\ar[r]^<<<<<{\mu_{k,n}}&
\mathscr S^{\mc O}_{H_{k,n}}
}
\eeq
is then automatic.
\endproof

\begin{rem}\label{generalell}
The commutative diagram \eqref{commprisms} admits further generalizations. To construct solutions of the joint qDE and qKZ system for $G(k,n)$ starting from corresponding solutions for $\Pb^{n-1}$, we invoke Theorem \ref{satsolpg}. In particular, a determination of $\arg(-1)$ has to be fixed in formula \eqref{finalsolprod}. In formula \eqref{05.02.24-1}, the determination $\arg(-1)=\pi$ has been chosen. Different choices lead to commutative diagrams analog to \eqref{commprisms}, with suitably modified vertices and arrows.

More precisely, the choice $\arg(-1)=(2\ell+1)\pi$, with $\ell\in\Z$, leads to a commutative diagram analog to \eqref{commdiagpn1}, the vertex $\mathscr S^{\mc O}_{H_{1,n}}(k)$ being replaced by the space $\mathscr S^{\mc O}_{H_{1,n}}(k,\ell)$ of sections $f(\bm z;\bm p;\kappa)$ of $H_{1,n}$ of the form
\[f(\bm z;\bm p;\kappa)=g(\bm z;e^{\pi\sqrt{-1}(2\ell+1)(k-1)}p_1,p_2;\kappa),\quad g\in\mathscr S^{\mc O}_{H_{1,n}},
\]and the arrows $A_1,A_2,A_3$ being replaced by
\begin{align*}
&A_{1,\ell}(V)=\text{\textnormal{\textcyr{B}}}_\Pb(V;\kappa)\exp[\frac{\pi\sqrt{-1}}{\kappa}(2\ell+1)(k-1)c_1^{T}(E_1^\Pb)],\quad &V\in K^T_0(\Pb^{n-1})_\C,\\ 
&A_{2,\ell}(v)=\mc S_{\Pb,k}(\bm z;\bm p;\kappa)\exp[-\frac{\pi\sqrt{-1}}{\kappa}(2\ell+1)(k-1)c_1^{T}(E_1^\Pb)]v,\quad &v\in H^{(\kappa)}_T(\Pb^{n-1}),\\
&A_{3,\ell}(V)=p_{k,n,\ell}[\mu_{1,n}(V)],\quad &V\in K^T_0(\Pb^{n-1})_\C,
\end{align*}
where the map $p_{k,n,\ell}\colon \mathscr S_{H_{1,n}}^{\mc O}\to \mathscr S^{\mc O}_{H_{1,n}}(k,\ell)$ is defined 
by 
\beq
p_{k,n,\ell}(f)(\bm z;\bm p;\kappa)=f\left(\bm z;e^{\pi\sqrt{-1}(2\ell+1)(k-1)}p_1,p_2;\kappa\right).
\eeq
If we define $\Hat{\mathscr{S}}_{k,n,\ell}$ 
as the space of sections of $\bigwedge^kH_{1,n}$ of the form
\[(\bm z;\bm p)\mapsto p_2^{(1-k)\sum_{a=1}^nz_a/\kappa}\,f_1(\bm z;\bm p;\kappa)\wedge\dots\wedge f_k(\bm z;\bm p;\kappa),\quad f_i\in\mathscr S^O_{H_{1,n}}(k,\ell),
\]for each $\ell\in\Z$ we have a commutative diagram
\beq\label{commprismsell}
\xymatrix{
\bigwedge\nolimits^kK^T_0(\Pb^{n-1})_\C\ar[rr]^{B_{1,\ell}}\ar[dr]_{B_{3,\ell}}\ar[dd]_{\thi_{k,n,\ell}^{K}}&&\bigwedge\nolimits^kH^{(\kappa)}_T(\Pb^{n-1})\ar[dl]^{B_{2,\ell}}\ar[dd]^{\thi_{k,n}}\\
&\Hat{\mathscr{S}}_{k,n,\ell}\ar[dd]^<<<<<<<<{\thi_{k,n}^\diamond\circ}&\\
K^T_0(G(k,n))_\C\ar[dr]_{\mu_{k,n}}\ar[rr]|\hole^<<<<<<<{\textnormal{\textcyr{B}}_{G}^{(\kappa)}}&&H^{(\kappa)}_T(G(k,n))\ar[dl]^{\mc S_G}\\
&\mathscr S^{\mc O}_{H_{k,n}}&
}
\eeq
where 
\begin{align*}
&B_{1,\ell}(V_{1}\wedge\dots\wedge V_k)=\left(\frac{2\pi\sqrt{-1}}{\kappa}\right)^{-\binom{k}{2}}\left(A_{1,\ell}(V_1)\wedge\dots\wedge A_{1,\ell}(V_k)\right),\\
&B_{2,\ell}(v_1\wedge\dots\wedge v_k)=p_2^{(1-k)\sum_{a=1}^nz_a/\kappa} \left(A_{2,\ell}(v_1)\wedge\dots\wedge A_{2,\ell}(v_k)\right),\\
&B_{3,\ell}(V_{1}\wedge\dots\wedge V_k)=\left(\frac{2\pi\sqrt{-1}}{\kappa}\right)^{-\binom{k}{2}}p_2^{(1-k)\sum_{a=1}^nz_a/\kappa} \left(A_{3,\ell}(V_1)\wedge\dots\wedge A_{3,\ell}(V_k)\right),
\end{align*}
for any $V_1,\dots, V_k\in K^T_0(\Pb^{n-1})_\C$, and $v_1,\dots,v_k\in H^{(\kappa)}_T(\Pb^{n-1})$, and the map $\thi^K_{k,n,\ell}$ is defined by
\[\vartheta^{K
}_{k,n,\ell}[\mc O(\la_1+k-1)\wedge\mc O(\la_2+k-2)\wedge\dots\wedge\mc O(\la_k)]=\mathscr L^{\otimes(\ell+1)}\otimes \Si_{\bm\la}E_1^*,
\]for any $\bm\la\subseteq k\times(n-k)$. In particular, for the choice $\ell=-1$ the map $\thi^K_{k,n,\ell}$ acts as the Schur functors $\Si_{\bm\la}$ only, with no twists by a power of $\mathscr L$. The choice $\ell=0$ in formulating Theorem \ref{thmcompatibility} is due to ``historical'' reasons, to have a perfect match with the computations of \cite{CDG1}, where $\mathscr L$ first appeared.
\end{rem}
\subsection{Determinantal identities for hypergeometric solutions} The commutativity of the diagram \eqref{finaldiag} can be reformulated as a non-trivial determinantal identity (Theorem \ref{detprop} below) relating multi-dimensional hypergeometric solutions of the qDE and qKZ systems associated with $\Pb^{n-1}$ and $G(k,n)$, respectively. 
\vskip2mm
Let $P_1,\dots,P_k\in \C[X^{\pm 1},\bm Z^{\pm 1}]$, and introduce the polynomial $P\in\C[{\bm\Gm}_1^{\pm 1},{\bm Z}^{\pm 1}]$ by
\beq\label{piplus} P(\bm\Gm_1;\bm Z)= \frac{\det
\left(\begin{array}{ccc}P_1(\Gm_{1,1};\bm Z)&\dots&P_k(\Gm_{1,1};\bm Z)\\
\vdots&&\vdots\\
P_1(\Gm_{1,k};\bm Z)&\dots&P_k(\Gm_{1,k};\bm Z)
\end{array}\right)}{\left(\prod_{i=1}^k\Gm_{1,i}\right)^{k-1}\frak D_k(\bm\Gm^{-1})}.
\eeq
The $K$-theoretical classes $[P_1],\dots,[P_k]\in K^T_0(\mathbb P^{n-1})_\C$ and $[P]\in K^T_0(G(k,n))_\C$ are then related by the identity
\[[P]=\thi^K_{k,n}\left([P_1]\wedge\dots\wedge[P_k]\right).
\]Introduce also the functions $\mu_1(\bm z;\bm p;\kappa),\dots,\mu_k(\bm z;\bm p;\kappa)$ and $\mu(\bm z;\bm p;\kappa)$ as Jackson integrals, using the morphisms $\mu_{k,n}$ of Section \ref{secideksol} (see eq. \eqref{defmukn}):
\begin{align}
&\mu_i(\bm z;\bm p;\kappa):=\mu_{1,n}[P_i](\bm z;\bm p;\kappa),\quad i=1,\dots,k,\\
&\mu(\bm z;\bm p;\kappa):=\mu_{k,n}[P](\bm z;\bm p;\kappa).
\end{align}
We have 
\[\mu_i(\bm z;\bm p;\kappa)=\sum_{j=1}^n\mu_{ij}(\bm z;\bm p;\kappa)\on{Stab}_{[j]},\qquad \mu(\bm z;\bm p;\kappa)=\sum_{J\in{\mc I_{k,n}}}\mu_{J}(\bm z;\bm p;\kappa)\on{Stab}_J,
\]where
\begin{multline}
\mu_{ij}(\bm z;\bm p;\kappa)=\kappa^{-n}\sum_{h=1}^n\acute{P}_i(z_h;\bm z;\kappa)\mc M_{[h]}\!\left(\Phi_{1,n}W_{[j]}\right)(\bm z;\bm p;\kappa),\\ i=1,\dots,k,\quad j=1,\dots,n,
\end{multline}
and
\beq
\mu_J(\bm z;\bm p;\kappa)=\kappa^{-k(n-k)-k}\sum_{H\in\mc I_{k,n}}\acute{P}(\bm z_{\si_H};\bm z;\kappa)\mc M_{H}\!\left(\Phi_{k,n}W_{J}\right)(\bm z;\bm p;\kappa).
\eeq
\vskip2mm
The functions $\mu_{ij}(\bm z;\bm p;\kappa),\mu(\bm z;\bm p;\kappa)$ satisfy a remarkable determinantal identity.

\begin{thm}\label{detprop}Let $J=(J_1,J_2)\in\mc I_{k,n}$, with $J_1=\{j_1<\dots<j_k\}$. We have
\[\mu_J(\bm z;\bm p;\kappa)=\left(\frac{2\pi\sqrt{-1}}{\kappa}\right)^{-\binom{k}{2}}p_2^{(1-k)\sum_{a=1}^nz_a/\kappa} \det\left(\mu_{i,j_i}(\bm z;p_1e^{\pi\sqrt{-1}(k-1)},p_2;\kappa)\right)_{i=1}^k.
\]
\end{thm}
\proof
This is equivalent to the commutativity of the diagram \eqref{finaldiag}.
\endproof

\begin{example}\label{exdetid}
Consider the cases $(k,n)=(1,2),(2,2)$. These correspond to the Grassmannians $\Pb^1$ and $G(2,2)=\on{pt}$, respectively. 
\vskip1,5mm
For $(k,n)=(1,2)$, the weight functions are
\[W_{[1]}(t;\bm z)=1,\quad W_{[2]}(t;\bm z)= t-z_1.
\]The master function is 
\[\Phi_{1,2}(t;\bm z;\bm p;\kappa)=(p_2/\kappa)^{(z_1+z_2)/\kappa}(\kappa^2p_1/p_2)^{t/\kappa}\Gm\left(\frac{t-z_1}{\kappa}\right)\Gm\left(\frac{t-z_2}{\kappa}\right).
\]Let $P_1(X;\bm Z), P_2(X;\bm Z)$ be two Laurent polynomials. We have four functions $\mu_{ij}$, for $i,j=1,2$:
\begin{multline*}\mu_{i1}(\bm z;\bm p;\kappa)=\kappa^{-2}\left(\acute{P}_i(z_1;\bm z;\kappa)\mc M_{[1]}\Phi_{1,2}(\bm z;\bm p;\kappa)+\acute{P}_i(z_2;\bm z;\kappa)\mc M_{[2]}\Phi_{1,2}(\bm z;\bm p;\kappa)\right)\\
=\kappa^{-2}\left(\acute{P}_i(z_1;\bm z;\kappa)\kappa^{\frac{z_1-z_2}{\kappa}+1}p_1^{z_1/\kappa}p_2^{z_2/\kappa}\sum_{\ell=0}^\infty\frac{(-1)^\ell}{\ell!}\Gm\left(\frac{z_1-z_2-\ell\kappa}{\kappa}\right)\left(\kappa^{-2}\frac{p_2}{p_1}\right)^\ell\right.\\
\left.+\acute{P}_i(z_2;\bm z;\kappa)\kappa^{\frac{z_2-z_1}{\kappa}+1}p_1^{z_2/\kappa}p_2^{z_1/\kappa}\sum_{h=0}^\infty\frac{(-1)^h}{h!}\Gm\left(\frac{z_2-z_1-h\kappa}{\kappa}\right)\left(\kappa^{-2}\frac{p_2}{p_1}\right)^h\right),
\end{multline*}
\begin{multline*}\mu_{i2}(\bm z;\bm p;\kappa)=\kappa^{-2}\left(\acute{P}_i(z_1;\bm z;\kappa)\mc M_{[1]}\left[\Phi_{1,2}W_{[2]}\right](\bm z;\bm p;\kappa)+\acute{P}_i(z_2;\bm z;\kappa)\mc M_{[2]}\left[\Phi_{1,2}W_{[2]}\right](\bm z;\bm p;\kappa)\right)\\
=\kappa^{-2}\left(\acute{P}_i(z_1;\bm z;\kappa)\kappa^{\frac{z_1-z_2}{\kappa}+2}p_1^{z_1/\kappa}p_2^{z_2/\kappa}\sum_{a=0}^\infty\frac{(-1)^{a+1}a}{a!}\Gm\left(\frac{z_1-z_2-a\kappa}{\kappa}\right)\left(\kappa^{-2}\frac{p_2}{p_1}\right)^a\right.\\
\left.+\acute{P}_i(z_2;\bm z;\kappa)\kappa^{\frac{z_2-z_1}{\kappa}+1}p_1^{z_2/\kappa}p_2^{z_1/\kappa}\sum_{b=0}^\infty\frac{(-1)^b}{b!}(z_2-z_1-b\kappa)\Gm\left(\frac{z_2-z_1-b\kappa}{\kappa}\right)\left(\kappa^{-2}\frac{p_2}{p_1}\right)^b\right),
\end{multline*}
where $i=1,2$.

Let us check the validity of Theorem \ref{detprop}, at least up to the first order. We have
\begin{multline*}
\det\begin{pmatrix}
\mu_{11}(\bm z;\bm p;\kappa)&\mu_{12}(\bm z;\bm p;\kappa)\\
\mu_{21}(\bm z;\bm p;\kappa)&\mu_{22}(\bm z;\bm p;\kappa)
\end{pmatrix}\\=\kappa^{-4}\left(\kappa^2\acute{P}_1(z_1;\bm z;\kappa)\acute{P}_2(z_2;\bm z;\kappa)(p_1p_2)^{(z_1+z_2)/\kappa}(z_2-z_1)\Gm\left(\frac{z_1-z_2}{\kappa}\right)\Gm\left(\frac{z_2-z_1}{\kappa}\right)\right.\\
\left.-\kappa^2\acute{P}_2(z_1;\bm z;\kappa)\acute{P}_1(z_2;\bm z;\kappa)(p_1p_2)^{(z_1+z_2)/\kappa}(z_2-z_1)\Gm\left(\frac{z_1-z_2}{\kappa}\right)\Gm\left(\frac{z_2-z_1}{\kappa}\right)\right)(1+O(p_2/p_1))\\
=\kappa^{-2}\det\begin{pmatrix}
\acute{P}_1(z_1;\bm z;\kappa)&\acute{P}_2(z_1;\bm z;\kappa)\\
\acute{P}_1(z_2;\bm z;\kappa)&\acute{P}_2(z_2;\bm z;\kappa)
\end{pmatrix}(p_1p_2)^{(z_1+z_2)/\kappa}\\
\times(z_2-z_1)\Gm\left(\frac{z_1-z_2}{\kappa}\right)\Gm\left(\frac{z_2-z_1}{\kappa}\right)(1+O(p_2/p_1))\\
=\kappa^{-1}\det\begin{pmatrix}
\acute{P}_1(z_1;\bm z;\kappa)&\acute{P}_2(z_1;\bm z;\kappa)\\
\acute{P}_1(z_2;\bm z;\kappa)&\acute{P}_2(z_2;\bm z;\kappa)
\end{pmatrix}(p_1p_2)^{(z_1+z_2)/\kappa}\\
\times \frac{2\pi\sqrt{-1}}{\exp(\frac{\pi\sqrt{-1}}{\kappa}(z_1-z_2))-\exp(\frac{\pi\sqrt{-1}}{\kappa}(z_2-z_1))}(1+O(p_2/p_1)),
\end{multline*}
where we have invoked the identity
\[t\Gm(t)\Gm(-t)=\frac{2\pi\sqrt{-1}}{e^{-\pi\sqrt{-1}t}-e^{\pi\sqrt{-1}t}}.
\]
We deduce
\begin{multline}\label{detid1}\frac{\kappa}{2\pi\sqrt{-1}}p_2^{-(z_1+z_2)/\kappa}\det\begin{pmatrix}
\mu_{11}(\bm z;e^{\pi\sqrt{-1}}p_1,p_2;\kappa)&\mu_{12}(\bm z;e^{\pi\sqrt{-1}}p_1,p_2;\kappa)\\
\mu_{21}(\bm z;e^{\pi\sqrt{-1}}p_1,p_2;\kappa)&\mu_{22}(\bm z;e^{\pi\sqrt{-1}}p_1,p_2;\kappa)
\end{pmatrix}\\
=\det\begin{pmatrix}
\acute{P}_1(z_1;\bm z;\kappa)&\acute{P}_2(z_1;\bm z;\kappa)\\
\acute{P}_1(z_2;\bm z;\kappa)&\acute{P}_2(z_2;\bm z;\kappa)
\end{pmatrix}\frac{\exp{(\frac{2\pi\sqrt{-1}}{\kappa}(z_1+z_2))}}{\exp{(\frac{2\pi\sqrt{-1}}{\kappa}z_1)}-\exp{(\frac{2\pi\sqrt{-1}}{\kappa}z_2)}}p_1^{(z_1+z_2)/\kappa}(1+O(p_2/p_1)).
\end{multline}
For $(k,n)=(2,2)$, we have a single weight function
\[W_{(\{1,2\},\emptyset)}(\bm t;\bm z)=1.
\]The master function is 
\[\Phi_{2,2}(\bm t;\bm z;\bm p;\kappa)=(\kappa^{-2}p_2)^{(z_1+z_2)/\kappa}\left(\kappa^2p_1/p_2\right)^{(t_1+t_2)/\kappa}\frac{\Gm\left(\frac{t_1-z_1}{\kappa}\right)\Gm\left(\frac{t_1-z_2}{\kappa}\right)\Gm\left(\frac{t_2-z_1}{\kappa}\right)\Gm\left(\frac{t_2-z_2}{\kappa}\right)}{\Gm\left(\frac{t_1-t_2}{\kappa}\right)\Gm\left(\frac{t_2-t_1}{\kappa}\right)}.
\]
We also have 
\[\acute{P}(\bm \gm;\bm z;\kappa)=\det\begin{pmatrix}
\acute{P}_1(\gm_{1,1};\bm z;\kappa)&\acute{P}_2(\gm_{1,1};\bm z;\kappa)\\
\acute{P}_1(\gm_{1,2};\bm z;\kappa)&\acute{P}_2(\gm_{1,2};\bm z;\kappa)
\end{pmatrix}\frac{\exp{(\frac{2\pi\sqrt{-1}}{\kappa}(\gm_{1,1}+\gm_{1,2}))}}{\exp{(\frac{2\pi\sqrt{-1}}{\kappa}\gm_{1,1})}-\exp{(\frac{2\pi\sqrt{-1}}{\kappa}\gm_{1,2})}}.
\]
A simple computation finally shows that the leading terms of the functions $\mc \mu_{(\{1,2\},\emptyset)}(\bm z;\bm p;\kappa)$ equal
\beq\label{detid3}\mc \mu_{(\{1,2\},\emptyset)}(\bm z;\bm p;\kappa)=\acute{P}(\bm z;\bm z;\kappa)p_1^{(z_1+z_2)/\kappa}(1+O(p_2/p_1)).
\eeq
The equality of \eqref{detid1}, 
\eqref{detid3} is in accordance with the determinantal identity of Theorem \ref{detprop}. Higher order terms can be identified similarly. We leave the check to the reader.
\qetr
\end{example}

\begin{rem}
We derived Theorem \ref{detprop} directly from Theorem \ref{thmcompatibility} (in particular, from the commutative diagram \eqref{commprisms}). Theorem \ref{detprop}, however, could be directly proved via some lengthy computations -- in the same spirit of Example \ref{exdetid} above. Consequently, assuming the validity of Proposition \ref{propsatbmor}, Proposition \ref{propsattopsol}, and the \textcyr{B}-theorem specifically for projective spaces, we deduce the logical equivalence among 
\begin{enumerate}
\item Theorem \ref{bth} (the \textcyr{B}-theorem for Grassmannians), 
\item Theorem \ref{thmcompatibility} (the commutativity of the diagram \eqref{commprisms}),
\item Theorem \ref{detprop} (the validity of the determinantal identities).
\end{enumerate}
Once any one of (1),(2),(3) is proved, the other two automatically follow.
\end{rem}

\begin{rem}
The commutative diagram \eqref{commprisms} can be generalized as discussed in Remark \ref{generalell}. For any choice of $\arg(-1)=(2\ell+1)\pi$, with $\ell\in \Z$, we have a corresponding commutative diagram \eqref{commprismsell}. For any such a choice of $\ell\in\Z$, we have a corresponding determinantal identity analogue to Theorem \ref{detprop}. More precisely, given $\ell\in\Z$, replace $P(\bm \Gm_1;\bm Z)$ in \eqref{piplus} by
\[P^{(\ell)}(\bm \Gm_1;\bm Z)=\frac{\det
\left(\begin{array}{ccc}P_1(\Gm_{1,1};\bm Z)&\dots&P_k(\Gm_{1,1};\bm Z)\\
\vdots&&\vdots\\
P_1(\Gm_{1,k};\bm Z)&\dots&P_k(\Gm_{1,k};\bm Z)
\end{array}\right)}{\left(\prod_{i=1}^k\Gm_{1,i}\right)^{(k-1)(\ell+1)}\frak D_k(\bm\Gm^{-1})},
\]
and set 
\[\mu^{(\ell)}(\bm z;\bm p;\kappa):=\mu_{k,n}[P^{(\ell)}](\bm z;\bm p;\kappa)=\sum_{J\in{\mc I_{k,n}}}\mu_{J}^{(\ell)}(\bm z;\bm p;\kappa)\on{Stab}_J.
\]
If $J=(J_1,J_2)\in\mc I_{k,n}$, with $J_1=\{j_1<\dots<j_k\}$, then the function $\mu_J^{(\ell)}(\bm z;\bm p;\kappa)$ satisfy the determinantal identity
\[\mu_J^{(\ell)}(\bm z;\bm p;\kappa)=\left(\frac{2\pi\sqrt{-1}}{\kappa}\right)^{-\binom{k}{2}}p_2^{(1-k)\sum_{a=1}^nz_a/\kappa} \det\left(\mu_{i,j_i}(\bm z;p_1e^{\pi\sqrt{-1}(k-1)(2\ell+1)},p_2;\kappa)\right)_{i=1}^k.
\]
\end{rem}

\subsection{Integral representations for hypergeometric solutions} In this section, as a further output, we show how the determinantal identity of Theorem \ref{detprop} gives rise to integral representations of solutions of the qDE and qKZ equations for Grassmannians.
\vskip2mm
Introduce the modified master function
\begin{multline}
\Tilde\Phi_{k,n}(\bm t;\bm z;\bm p;\kappa):=\prod_{i=1}^k\Phi_{1,n}(t_i;\bm z;p_1e^{\pi\sqrt{-1}(k-1)},p_2;\kappa)\\
=\left(\frac{p_2}{\kappa}\right)^{k\sum_{a=1}^nz_a/\kappa}\left(\frac{\kappa^n p_1}{p_2}\right)^{\sum_{i=1}^kt_i/\kappa}e^{\pi\sqrt{-1}(k-1)\sum_{i=1}^kt_i/\kappa}\prod_{i=1}^k\prod_{a=1}^n\Gamma\left(\frac{t_i-z_a}{\kappa}\right),
\end{multline}
and, for $J=(J_1,J_2)\in\mc I_{k,n}$, with $J_1=\{j_1<\dots<j_k\}$, the modified weight function
\beq
\Tilde W_J(\bm t;\bm z):=\det\begin{pmatrix}
W_{[j_1]}(t_1;\bm z)&W_{[j_2]}(t_1;\bm z)&\dots&W_{[j_k]}(t_1;\bm z)\\
W_{[j_1]}(t_2;\bm z)&W_{[j_2]}(t_2;\bm z)&\dots&W_{[j_k]}(t_2;\bm z)\\
&\ddots&\\
W_{[j_1]}(t_k;\bm z)&W_{[j_2]}(t_k;\bm z)&\dots&W_{[j_k]}(t_k;\bm z)
\end{pmatrix}.
\eeq
Also, collect the modified weight functions into the cohomological valued function
\beq
\Tilde W(\bm t;\bm\gm;\bm z):=\sum_{J\in\mc I_{k,n}}\Tilde W_J(\bm t;\bm z)\on{Stab}_J(\bm\gm;\bm z).
\eeq
\begin{rem}
The modified weight functions $\Tilde W_J(\bm t;\bm z)$ also admit an expression in factorial Schur polynomials, see equation \eqref{facschurpol}. Let $\bm\la=\bm\la^{\si_J}\subseteq k\times(n-k)$ be the Young diagram associated with the double partition $J=(J_1,J_2)$, with $J_1=\{j_1<j_2<\dots<j_k\}$, as described in Section \ref{secpgrpy}. In other words, let $\bm\la=(\la_1\geq\la_2\geq \la_3\geq\dots)$ with 
\[\la_a=j_{k-a+1}-k+a-1,\quad a=1,\dots,k.
\]Consider the factorial Schur polynomial $s_{\bm\la}(\bm x|\bm y)$ as defined in equation \eqref{facschurpol}. We have
\beq
\Tilde W_J(\bm t;\bm z)=\,s_{\bm\la}(\bm t|\bm z)\prod_{1\leq i< j\leq k}(t_j-t_i).
\eeq
\end{rem}
In the space $\C^k$ consider the $k$-dimensional contour $\mc C$ defined by the products of parabolas 
\[\mc C:=C_1\times\dots\times C_k,
\]where 
\[C_i:=\{\kappa(A_i-s^2+s\sqrt{-1})\colon s\in\R\}\subseteq \C,\quad A_i \text{ is a sufficiently large number.}
\]
\vskip2mm
As in the previous section, consider $P_1,\dots,P_k\in\C[X^{\pm 1}, \bm Z^{\pm 1}]$ and introduce the polynomial $P\in \C[{\bm \Gm}_1^{\pm 1},\bm Z^{\pm 1}]$ by \eqref{piplus}, so that the $K$-theoretical classes $[P_1],\dots,[P_k]\in K^T_0(\mathbb P^{n-1})_\C$ and $[P]\in K^T_0(G(k,n))_\C$ are related by the identity
\[[P]=\thi^K_{k,n}\left([P_1]\wedge\dots\wedge[P_k]\right).
\]
Define the $k$-dimensional integral function
\begin{multline}\label{multint}
\frak I[P_1,\dots,P_k](\bm z;\bm p;\kappa):=\\
\left(2\pi\sqrt{-1}\right)^{-\binom{k}{2}-k}\,\kappa^{\binom{k}{2}-nk}\,p_2^{(1-k)\sum_{a=1}^nz_a/\kappa} \int_{\mc C}\prod_{i=1}^k\acute{P_i}(t_i;\bm z;\kappa)\,\Tilde\Phi_{k,n}(\bm t;\bm z;\bm p;\kappa)\,\Tilde W(\bm t;\bm\gm;\bm z)\,d\bm t.
\end{multline}
By construction, for any fixed $\kappa$, the functions $(\bm z;\bm p)\mapsto \frak I[P_1,\dots,P_k](\bm z;\bm p;\kappa)$ define sections of the bundle $H_{k,n}\to\C^n\times\C^2$ (defined in Section \ref{secqdeqkz}).
\begin{thm}\label{intrepthm}
For any $k$-tuple of Laurent polynomials $(P_1,\dots,P_k)$, the section $\frak I[P_1,\dots,P_k]$ is a solution of the joint system of qDE and qKZ equations \eqref{qde.0} and \eqref{qkz.H}. Moreover, it coincides with the solution $\mu_{k,n}[P]$ associated with the $K$-theoretical class $[P]\in K^T_0(G(k,n))_\C$. In other words,
\[\frak I[P_1,\dots,P_k]=\mu_{k,n}[P]=\on{Stab}_{k,n}\Psi_P.
\]
\end{thm}
\proof
The result directly follows from Theorem \ref{detprop} and the integral representation of solutions of the qDE and qKZ system for $\mathbb P^{n-1}$, described in Example \ref{ex:intforpn}.
\endproof

\begin{rem}
For $k=1$, the integral \eqref{multint} exactly coincides with the one of equation \eqref{intpn-1}.
\end{rem}

\section{Monodromy and Stokes phenomenon of qDE and qKZ equations}\label{SEC8}
\subsection{Specialization of parameters} 
By definition the qDE and qKZ joint system of equations for $G(k,n)$ is a system of differential-difference equations for sections of the trivial bundle $H_{k,n}\to\C^n\times\C^2$. Recall that the coordinates on the base space $\C^n\times\C^2$ are $(\bm z;\bm p)$, and that the fiber over $(\bm z_o;\bm p_o)$ is $H^\bullet_T(G(k,n),\C)_{\bm z_o}$, see Section \ref{secqdeqkz}.
In this last section, we will focus on the joint system of qDE and qKZ equations \eqref{qde.0} and \eqref{qkz.H} for $G(k,n)$, under the specialization of parameters 
\beq\label{specpar}
p_1=q^{-1},\quad p_2=1,\quad \kappa=-1.
\eeq
Such a specialization restricts the base space to $\C^n\times \C$, with coordinates $(\bm z; q)$. In particular, we will deal with a single equivariant qDE w.r.t.\,\,$q$, namely
\beq\label{qde.q}
q\frac{\der}{\der q}f(\bm z; q)=c_1^{T}(E_1)*_{q,\bm z} f(\bm z;q).
\eeq
This differential equation has two singularities: a Fuchsian one at \( q = 0 \) and an irregular one at \( q = \infty \). As a result, solutions are typically multi-valued, exhibit {\it monodromy}, and also display a {\it Stokes phenomenon} at \( q = \infty \). In this section, we study these aspects and relate them to the algebraic and geometric features of the derived category \( \mathsf{D}^b_T(G(k,n)) \).
\subsection{Monodromy operator}\label{seccmonop}Following the notation of Section \ref{secext}, let $L_{-1}\subseteq\C^n$ to be the complement of the hyperplanes $z_i-z_j=m,\quad i,j=1,\dots,n,\quad m\in\Z$. 

Similarly to Section \ref{secbthm1}, we denote by $\mathscr S_{q}^{\mc O}$ the space of of sections $f(\bm z;q)$ of $H_{k,n}|_{\C^n\times\C}$ solving the equivariant differential equations \eqref{qde.q}, and  holomorphic in $\bm z\in L_{-1}$. For a fixed $\bm z_o$, we set 
\[\mathscr S_{q,\bm z_o}^{\mc O}:=\{f(\bm z_o;q)\colon\quad f(\bm z;q)\in \mathscr S_{q}^{\mc O}\}.
\]

Fix $(\bm z;q)$ and increase the argument of $q$ by $2\pi$. The analytic continuation of the solutions of \eqref{qde.q} along this curve produces the \emph{monodromy operator} $M_0(\bm z)$ on the space of solutions $\mathscr S_{q,\bm z}^{\mc O}$:
\[f(\bm z;q)\mapsto f(\bm z;e^{2\pi\sqrt{-1}}q).
\] 
\vskip2mm
For a fixed $\bm z\in\C^n$, consider also the morphism $\mu_{k,n,\bm z}\colon K^T_0(G(k,n))_\C\to \mathscr S^{\mc O}_{q,\bm z}$ induced by $\mu_{k,n}$ of Section \ref{secideksol}, that is
\beq
[P]\mapsto \left(\on{Stab}_{k,n}\Psi_P\right)(\bm z;q).
\eeq
\vskip2mm  As in Section \ref{seclev}, we denote by $E_I(\bm z)$, with $I\in\mc I_{k,n}$, the polynomials
\[E_I(\bm z):=\sum_{i\in I_1}z_i,\quad I=(I_1,I_2),
\]and by $D_{I,J}$, with $I,J\in\mc I_{k,n}$, the set of points $\bm z\in\C^n$ such that $E_I(\bm z)-E_J(\bm z)\in\Z_{\geq 0}$ and $E_I(\bm z)\neq E_J(\bm z)$. Set $D_{k,n}:=\bigcup_{I,J\in\mc I_{k,n}}D_{I,J}$, and
\beq\label{eqgenericz}\mc U:=L_{-1}\cap\left(\C^n\setminus D_{k,n}\right).
\eeq
\begin{prop}\label{propmono}
For $\bm z\in \mc U$, the morphism $\mu_{k,n,\bm z}$ intertwines the operator 
\[K^T_0(G(k,n))_\C\to K^T_0(G(k,n))_\C,\quad [V]\mapsto [\det E_1\otimes V],
\]and the monodromy operator $M_0(\bm z)\colon \mathscr S^{\mc O}_{q,\bm z}\to \mathscr S^{\mc O}_{q,\bm z}$. 
\end{prop}

\proof
The result follows from the expression of the Levelt fundamental solution of Section \ref{seclev}. Notice that the operator $X_1(\bm z;\oslash)$ has eigenvalues $E_I(\bm z)$. Under the isomorphism $\on{Stab}_{k,n}^\diamond$ of Section \ref{secqdeqkz}, the operator $X_1(\bm z;\oslash)$ is identified with the operator of classical equivariant multiplication by $c_1^{T}(E_1)$. Consequently, the operator $\exp(2\pi\sqrt{-1} X_1(\bm z;\oslash))$ is identified with the classical equivariant multiplication by the Chern character $\on{Ch}^{(-1)}_T(\det E_1)=e^{2\pi\sqrt{-1}c_1^{T}(E_1)}$.
\endproof

\begin{cor}
For $\bm z\in\mc U$, the monodromy operator $M_0(\bm z)$ has eigenvalues $\exp(2\pi\sqrt{-1}E_I(\bm z))$, 
with $I\in\mc I_{k,n}$.\qed
\end{cor}

\subsection{Formal reduction and formal solutions} Introduce the branched covering  \[\psi\colon \C^n\times\C\to\C^n\times\C,\quad(\bm z,s)\mapsto (\bm z,s^n),\]and consider the pulled-back bundle $\psi^*H_{k,n}$,
\[\xymatrix{
\psi^*H_{k,n}\ar[d]\ar[r]&H_{k,n}\ar[d]\\
\C^n\times\C\ar[r]^{\psi}&\C^n\times\C
}\]By changing the variable $q=s^n$, the joint system of qDE and qKZ equations can be pulled-back to a system of differential-difference equations for sections of the bundle $\psi^*H_{k,n}$. The sections $(\bm z, s)\mapsto \on{Stab}_{I,\bm z}$, with $I\in\mc I_{k,n}$, define a trivialization of $\psi^*H_{k,n}$, which allows us to write the qDE and qKZ equations in matrix form. 
\vskip2mm
More specifically, given a section $f(\bm z;s)$ of $\psi^*H_{k,n}$, we have
\[f(\bm z;s)=\sum_{I\in\mc I_{k,n}}f^I(\bm z;s)\on{Stab}_{I,\bm z},
\]for suitable components $f^I(\bm z;s)$, which can be arranged into a column vector. More generally, if we have a basis of solutions $f_1(\bm z;q),\dots f_{\binom{n}{k}}(\bm z;q)$, we can arrange their components into a $\binom{n}{k}\times \binom{n}{k}$-matrix $Y(\bm z;s)$, with entries
\[Y^I_i(\bm z;q)=f^I_i(\bm z;q),\quad I\in\mc I_{k,n},\quad i=1,\dots,\binom{n}{k}.
\]
 Hence, the qDE and qKZ equations for sections of $\psi^*H_{k,n}$ take the matrix form
\begin{align}\label{js1}
\frac{d}{ds}Y(\bm z;s)&=\mathsf A(\bm z;s)Y(\bm z;s),\\ 
\label{js1.1}Y(z_1,\dots,z_j-1,\dots,z_n;s)&= \mathsf K_j(\bm z;s)Y(\bm z;s),\quad j=1,\dots,n,
\end{align}
for suitable matrix-valued functions $\mathsf A(\bm z;s)$ and $\mathsf K_j(\bm z;s)$, see Remark \ref{akcoeffs} below. A solution $Y(\bm z; s)$ whose columns define a basis of solutions will be called a {\it fundamental solution matrix}\footnote{We will use this terminology also for solutions $Y(\bm z;s)$ of the qDE or qKZ equations only.} of the system \eqref{js1}, \eqref{js2}.
\begin{rem}\label{akcoeffs}
Let $X_1(\bm z;q)$ and $K_j(\bm z;q;\kappa)$, with $j=1,\dots,n$, be the matrices obtained by specializing the dynamical operator $X_1(\bm z;\bm p)$ and the qKZ operators $K_j(\bm z;\bm p;\kappa)$ of equations \eqref{dynop1}, \eqref{qkzop} to $\bm p=(q^{-1},1)$ and $\kappa=-1$. The matrices $\mathsf A(\bm z;s)$ and $\mathsf K_j(\bm z;s)$ are given by
\[\mathsf A(\bm z;s)=\frac{n}{s}X_1(\bm z;s^n),\quad \mathsf K_j(\bm z;s)=K_j(\bm z;s^n),\quad j=1,\dots,n.
\]
\end{rem}
\vskip2mm
Different choices of trivialization of $\psi^*H_{k,n}$ lead to different matrix representations of the joint system of qDE and qKZ equations: one says that all these matrix representations are {\it gauge equivalent}.
Given a different frame  $\left(g_I(\bm z;s)\right)$ of sections of $\psi^*H_{k,n}$, we can similarly expand the section $f_i(\bm z;s)$, $i=1,\dots,\binom{n}{k}$, w.r.t.\,\,the new frame,
\[f_i(\bm z;s)=\sum_{I\in\mc I_{k,n}}Z^I_i(\bm z;s)\,g_I(\bm z; s),\quad i=1,\dots,\binom{n}{k},
\]
 by obtaining a new matrix form for the joint system of differential-difference equations: 
 \begin{align}\label{js2}
\frac{d}{ds}Z(\bm z;s)&=\mathsf A'(\bm z;s)Z(\bm z;s),\\ 
\label{js2.1}Z(z_1,\dots,z_j-1,\dots,z_n;s)&= \mathsf K_j'(\bm z;s)Z(\bm z;s),\quad j=1,\dots,n,
\end{align}for suitable coefficients $\mathsf A'(\bm z;s)$ and $\mathsf K'_j(\bm z;s)$. To describe the relation between \eqref{js1}, \eqref{js1.1} and \eqref{js2}, \eqref{js2.1} introduce the matrix-valued function $G(\bm z;s)=\left(G(\bm z;s)^J_I\right)_{I,J\in\mc I_{k,n}}$ satisfying
\[\on{Stab}_{I,\bm z}=\sum_{J\in\mc I_{k,n}}G(\bm z;s)^J_I\,g_J(\bm z;s).
\]
We have
\begin{align}Z(\bm z;s)&=G(\bm z;s)Y(\bm z;s),\\
\mathsf A'(\bm z;s)&=G(\bm z;s)\mathsf A(\bm z;s)G(\bm z;s)^{-1}+\frac{dG(\bm z;s)}{ds}G(\bm z;s)^{-1},\\
\mathsf K'_j(\bm z;s)&=G(z_1,\dots,z_j-1,\dots,z_n;s)\mathsf K_j(\bm z;s)G(\bm z;s)^{-1},\quad j=1,\dots,n.
\end{align}
\vskip1,5mm
In this section we start by studying the {\it formal gauge equivalence} class, at $s=\infty$, of the joint system of qDE and qKZ equations. More precisely, we allow $G(\bm z;s)$ to be a matrix-valued function admitting only a formal Laurent series expansion in $s^{-1}$, not necessarily convergent. 
\vskip2mm
Set $\zeta_n:=\exp\left(\frac{2\pi\sqrt{-1}}{n}\right)$. For each $I\in\mc I_{k,n}$, with $I=(I_1,I_2)$ and $I_1=\{i_1<\dots<i_k\}$, introduce the sums and products of $n$-th roots of unity
\beq
\frak s_I:=\sum_{\ell=1}^k\zeta_n^{i_\ell-1},\qquad \frak p_I:=\prod_{\ell=1}^k\zeta_n^{1-i_\ell}.
\eeq
\begin{thm}\label{forredthm}
The pulled-back joint system of qDE and qKZ equations for sections of $\psi^*H_{k,n}$ is formally gauge equivalent, at $s=\infty$, to the joint system in matrix form
\begin{align}
\label{fr1}
&\frac{d}{ds}Z(\bm z,s)=UZ(\bm z,s)\\
\label{fr2}
&Z(z_1,\dots,z_j-1,\dots,z_n,s)=\mc K_jZ(\bm z,s),\quad j=1,\dots,n,
\end{align}
where $U=\diag{\left(n\,e^{\frac{\pi\sqrt{-1}(k-1)}{n}}\frak s_I\right)}_{I\in\mc I_{k,n}}$, and $\mc K_j=\diag{\left(e^{-\frac{\pi\sqrt{-1}k(k-1)}{n}}\frak p_I\right)}_{I\in\mc I_{k,n}}$ for each $j=1,\dots,n$. More precisely, there exists a formal gauge transformation matrix $G(\bm z;s)$ relating the system \eqref{js1}, \eqref{js1.1} to the system \eqref{fr1}, \eqref{fr2} of the form
\beq\label{forgau0}
G(\bm z;s)=\left(H(s)F(\bm z;s)s^{\La_{k,n}(\bm z)}\right)^{-1},
\eeq
where 
\beq
H(s)\text{ is a matrix-valued polynomial in }s^{-1},\text{ invertible for $s\neq 0$,}
\eeq
\beq\label{forgau2}
F(\bm z;s)={\bf 1}+\sum_{i=1}^\infty\frac{F_i(\bm z)}{s^i},\quad F_i(\bm z)\text{ polynomials,}
\eeq
\beq
\La_{k,n}(\bm z)=k\left(\frac{n-1}{2}+\sum_{a=1}^nz_a\right)\cdot{\bf 1},\quad\text{called \emph{exponent of formal monodromy}}.
\eeq
\end{thm}

The theorem claims that, after a formal gauge transformation, the system of qDE and qKZ equations becomes a system with constants coefficients and separated variables. Furthermore, the system splits into the direct sum of systems of rank one.

The joint system \eqref{fr1}, \eqref{fr2} admits a basis if solutions
\beq\label{forsol0}
Z_I(\bm z;s)=\exp\left(n\,e^{\frac{\pi\sqrt{-1}(k-1)}{n}}\frak s_I\, s+\left(\frac{\pi\sqrt{-1}k(k-1)}{n}-\log\frak p_I\right)\sum_{a=1}^nz_a\right)\cdot 1_I,\quad I\in\mc I_{k,n},
\eeq
where $1_I$ denotes the column vector with the only non-zero entry $1$ at the $I$-th  position.
All solutions of the system \eqref{fr1}, \eqref{fr2} are linear combinations of these basis solutions, with coefficients $1$-periodic in $z_1,\dots,z_n$.

\proof[Proof of Theorem \ref{forredthm}] The statement has been proved for projective spaces in \cite[Thm.\,9.2, Thm.\,9.6, Cor.\,10.14]{CV21}. For the Grassmannian $G(k,n)$ the result follows, via Satake correspondence, from the analog one for $\Pb^{n-1}$. Indeed, denote by $Z^\Pb(\bm z;s)$ (resp.\,\,$Z^G(\bm z;s)$) the matrix obtained by collecting the column vectors \eqref{forsol0} for $\mathbb P^{n-1}$ (resp.\,\,$G(k,n)$). 
Also, denote by $G^\Pb(\bm z;s)$ the formal gauge transformation matrix \eqref{forgau0} for $\mathbb P^{n-1}$. Then, the qDE and qKZ equations for $\Pb^{n-1}$ admit a formal fundamental solution matrix
\beq\label{forsolp}
Y^{\Pb}_{\rm for}(\bm z;s):=G^{\Pb}(\bm z;s)^{-1}Z^{\Pb}(\bm z;s).
\eeq
By replacing $s\mapsto e^{\frac{\pi\sqrt{-1}(k-1)}{n}}s$ in \eqref{forsolp}, and by taking the $k$-th exterior power of the r.h.s.\,\,(that is the matrix of $k\times k$ minors), we obtain a formal solution $Y^G_{\rm for}(\bm z;s)$ of the joint system of qDE and qKZ equations for $G(k,n)$, by Theorem \ref{satsolpg}. More explicitly, we have
\begin{multline}
Y^G_{\rm for}(\bm z;s):=\bigwedge\nolimits^k\left[H\left(e^{\frac{\pi\sqrt{-1}(k-1)}{n}}s\right)F\left(\bm z;e^{\frac{\pi\sqrt{-1}(k-1)}{n}}s\right)	\left(e^{\frac{\pi\sqrt{-1}(k-1)}{n}}s\right)^{\La_{1,n}(\bm z)}Z^\Pb\left(\bm z;e^{\frac{\pi\sqrt{-1}(k-1)}{n}}s\right)\right]\\
=e^{\frac{\pi\sqrt{-1}k(k-1)(n-1)}{2n}}\Tilde H(s)\Tilde F(\bm z;s)s^{\La_{k,n}(\bm z)}Z^G(\bm z;s),
\end{multline}
where $\Tilde H(s)$ is a matrix-valued polynomial in $s^{-1}$, $\Tilde F(\bm z;s)$ is of the form \eqref{forgau2}. 
This proves that the system \eqref{js1}, \eqref{js1.1} is formally gauge equivalent to the system \eqref{fr1}, \eqref{fr2}, as claimed.
\endproof

\begin{rem}The change of variables \( q = s^n \) is needed to achieve the formal reduction described in Theorem \ref{forredthm}. If, instead, we work with the original system \eqref{qde.q}, it becomes necessary to introduce a formal gauge transformation that admits a formal Puiseux series expansion in \( q \). This is known as the {\it Fabry phenomenon}, first described in \cite{Fab85}. The underlying cause of this phenomenon is the nilpotency of the leading term in the coefficient of equation \eqref{qde.q} as $q\to\infty$. See \cite{wasow,BJL79a,div1,div2}.
\end{rem}

\begin{rem}
The coefficient $U$ in the differential equation \eqref{fr1} may have a simple spectrum (or may not), depending on the pair $(k,n)$. Surprisingly, this phenomenon is related to the prime numbers distribution: the matrix $U$ has not simple spectrum if and only if $\pi_1(n)\leq k\leq n-\pi_1(n)$, where $\pi_1(n)$ is the smallest prime factor of $n$, see \cite{Cot22}. In this case, one says that the singularity $s=\infty$ of the differential equation \eqref{js1} is {\it irregular with coalescing eigenvalues}.
\end{rem}

\begin{rem}
The formal gauge transformations $G(\bm z;s)$ as in Theorem \ref{forredthm} are not unique, since there is a freedom in the choice of $H(s)$ and $F(\bm z;s)$. If, however, $G_1(\bm z;s)$ and $G_2(\bm z;s)$ are of the form \eqref{forgau0} with the same coefficient $H(s)$, then one actually has $G_1(\bm z;s)=G_2(\bm z;s)$. It is possible to prove, indeed, that all the polynomials $F_i(\bm z)$ in \eqref{forgau2} can be iteratively computed. See \cite[Appendix A]{CV21} for the case of simple spectrum of $U$. For the case of coalescing eigenvalues see e.g.\,\,\cite{CDG0,CDG} where the algorithm in the general case is described. In \cite{CV21}, for the case of the projective space $\Pb^{n-1}$, the following choice for the matrix $H(s)=\left(H(s)^i_j\right)_{i,j=1}^n$ was done:
\[H(s)^i_j=\frac{s^{1-i}}{\sqrt{n}}\exp\left(\frac{(j-1)(1-2i)\sqrt{-1}\pi}{n}\right),\quad i,j=1,\dots, n.
\]
\end{rem}

Following the notation of Section \ref{secext}, let $L_{-1}\subseteq\C^n$ to be the complement of the hyperplanes 
\[z_i-z_j=m,\quad i,j=1,\dots,n,\quad m\in\Z.
\]
\begin{cor}\label{corforsol}
Let $C(\bm z)$ a $\binom{n}{k}\times\binom{n}{k}$ diagonal-matrix-valued function meromorphic on $\C^n$, regular on $L_{-1}$, and with non-vanishing entries on $L_{-1}$. If $C(\bm z)$ is a fundamental solution matrix of the difference equations \eqref{fr2}, over the ring of 1-periodic functions in $\bm z$, then the joint system \eqref{js1}, \eqref{js1.1} for $G(k,n)$ admits a formal solution of the form
\beq
Y_{\rm for}(\bm z;s)=G(\bm z;s)^{-1}e^{Us}C(\bm z).
\eeq
\end{cor}
The solution $Y_{\rm for}(\bm z;s)$ is typically divergent. Nevertheless, it contains information about analytic solutions of the differential system \eqref{js1}: namely, it prescribes the Poincar\'e asymptotics of analytic fundamental solution matrices in suitable sectors of the universal cover of the punctured $s$-plane $\C^*$.

\subsection{Stokes sectors, Stokes bases, and Stokes matrices} Following the notations of \cite{TV21,CV21}, introduce coordinates $(r,\phi)$ on the universal cover $\Tilde{\C^*}$ of the punctured $s$-plane $\C^*$:
\beq
q=s^n,\qquad s=re^{-2\pi\sqrt{-1}\phi},\quad r>0,\,\phi\in\R.
\eeq

\begin{defn}
We call {\it Stokes ray} any ray in the universal cover $\Tilde{\C^*}$ of the $s$-plane defined by the equation
\beq
\on{Re}(\frak s_I\,s)=\on{Re}(\frak s_J\,s),\quad I,J\in\mc I_{k,n},\quad\frak s_I\neq\frak s_J.
\eeq
These are the rays along which at least one dominance relation between the growth orders $|\exp(U_{I}s)|$, for $|s|\to\infty$, changes. 
\end{defn}
\begin{defn}We call {\it Stokes sector} any open sector of $\Tilde{\C^*}$ which contains 
exactly one for any ordered pair $(I,J)$, with $I,J\in\mc I_{k,n}$, $I<J$ and $\frak s_I\neq \frak s_J$.
\end{defn}

The following well-known lemma allows us to easily find Stokes sectors, see \cite{painkapaev}.
\begin{lem}
Any open sector $\mathcal V\subset\widetilde{\mathbb C^*}$ of width $\pi+\delta$, i.e. of the form
\beq
\mathcal V=\left\{s\in\widetilde{\mathbb C^*}\colon a-\frac{1}{2}-\delta<\phi<a\right\},\quad a\in\mathbb R,
\eeq is a Stokes sector for $\delta>0$ sufficiently small.\qed
\end{lem}

The following theorem follows from the general theory of differential equations, see \cite{wasow,BJL79,fed,sibook,painkapaev} and also the new paper \cite{CGM23}.

\begin{thm}\label{teostok}
Let $\mathcal V\subseteq\widetilde{\mathbb C^*}$ be a Stokes sector, and let $Y_{\rm for}(\bm z;s)$ denote the formal solution described in Corollary \ref{corforsol}. There exists a unique fundamental solution matrix $Y(\bm z;s)$ of the qDE  \eqref{js1} satisfying the asymptotic condition 
\beq\label{asym} Y(\bm z;s)\sim Y_{\rm for}(\bm z;s),\quad s\to\infty,\quad s\in\mathcal V,
\eeq uniformly on compact subsets of $L_{-1}$.  The asymptotic expansion \eqref{asym} can actually be extended to a sector wider than $\mathcal  V$, up to the nearest Stokes rays.\qed
\end{thm}

\begin{rem}
In the notations of Theorem \ref{forredthm} and Corollary \ref{corforsol}, the precise meaning of the asymptotic relation \eqref{asym} is the following:  
\[
\forall K\Subset L_{-1},
\
 \forall h\in\mathbb N,
 \
  \forall \overline{\mathcal V}\subsetneq \mathcal V,
  \
   \exists C_{K,h,\overline{\mathcal V}}>0\colon \text{ if }s\in\overline{\mathcal V}\setminus\left\{0 \right\}\text{ then }\]
   \[
    \sup_{\bm z\in K}\left\|  H(s)^{-1}\cdot Y(\bm z; s)\cdot C(\bm z)^{-1}\exp(-s U)s^{-\Lambda_{k,n}(\bm z)}-\sum_{m=0}^{h-1}\frac{F_m(\bm z)}{s^m}\right\|
    <\frac{C_{K,h,\overline{\mathcal V}}}{|s|^h}.
\]
Here $\overline{\mathcal V}$ denotes any unbounded closed sector of $\widetilde{\mathbb C^*}$ with vertex at $s=\infty$, and $F_0(\bm z)={\bf 1}$.  Here, for defining $s^{\Lambda_{k,n}(\bm z)}$, the principal branch of $\log s$ is chosen.
\end{rem}

\begin{defn}
Given a sector $\mc V\subseteq\Tilde{\C^*}$, we say that a fundamental solution matrix $Y(\bm z;s)$ of the qDE \eqref{js1} 
 is a {\it Stokes fundamental solution with normalization $C(\bm z)$ on $\mc V$} if it satisfies the asymptotic expansion
\[Y(\bm z;s)\sim Y_{\rm for}(\bm z;s),\quad s\to\infty,\quad s\in\mc V,
\]uniformly on compact subsets of $L_{-1}$, where $Y_{\rm for}(\bm z;s)$ is as in Corollary \ref{corforsol}. We will call the corresponding basis of sections of $\psi^*H_{k,n}$ a {\it Stokes basis} of the qDE. \end{defn}

\begin{rem}
If $\mc V$ is a Stokes sector, and $C(\bm z)$ is a fixed normalization, there exists a unique Stokes fundamental solution for the qDE with normalization $C(\bm z)$ on $\mc V$, by Theorem \ref{teostok}. A priori, the qKZ difference equations \eqref{js1.1} may not be satisfied by this Stokes solution. If also the qKZ difference equations are solved by $Y(\bm z;q)$, then we will say that $Y(\bm z;q)$ is a Stokes fundamental solution of the joint system of qDE and qKZ equations. 

\end{rem}

Notice that if $\mathcal V$ is a Stokes sector, then also $e^{\pi\sqrt{-1}}\mathcal V$ and $e^{2\pi\sqrt{-1}}\mathcal V$ are Stokes sectors.

\begin{defn}\label{stokesmatr}
Let $Y(\bm z;s)$ be the Stokes fundamental solution (with normalization $C(\bm z)$) of the qDE \eqref{js1} on the Stokes sector $\mathcal V$. Let $Y_1(\bm z;s)$, and $Y_2(\bm z;s)$, be the Stokes fundamental solutions on $e^{\pi\sqrt{-1}}\mathcal V$ and on $e^{2\pi\sqrt{-1}}\mathcal V$, respectively. Define the \emph{ Stokes matrices} attached to $\mathcal V$ and $C(\bm z)$ as the matrices $\mathbb S_1,\mathbb S_2$ (depending on $\bm z\in L_{-1}$) such that
\beq\label{stokmatr}
Y_1(\bm z;s)=Y(\bm z;s)\mathbb S_1,\quad Y_2(\bm z;s)=Y_1(\bm z;s)\mathbb S_2,\quad s\in\widetilde{\mathbb C^*},\quad \bm z\in L_{-1}.
\eeq
\end{defn}

\subsection{Stokes bases as ${T}$-full exceptional collections}

The Satake identification allows us to construct the Stokes fundamental solutions for $G(k,n)$ starting with Stokes solutions for $\Pb^{n-1}$.

\begin{thm}\label{satstoksol}
Let $\mc V$ be a Stokes sector for the qDE of $\Pb^{n-1}$, and let $Y^\Pb(\bm z;s)$ be the Stokes fundamental solution of the qDE for $\Pb^{n-1}$ on $\mc V$ with normalization $C(\bm z)$. Under the Satake identification \eqref{sathi}, the $k$-th exterior power $\bigwedge^kY^{\Pb}(\bm z;s)$ is the Stokes fundamental solution of the qDE for $G(k,n)$ on the Stokes sector $e^{-\frac{\pi\sqrt{-1}(k-1)}{n}}\mc V$ with normalization
\[\Tilde C(\bm z)=\bigwedge\nolimits^k\left[\left(e^{\frac{\pi\sqrt{-1}(k-1)}{n}}\right)^{\La_{1,n}(\bm z)}C(\bm z)\right].
\]
\end{thm}

\proof
Take the asymptotic expansion \eqref{asym} for $Y(\bm z;s)=Y^\Pb(\bm z;s)$, replace $s\mapsto e^{\frac{\pi\sqrt{-1}(k-1)}{n}}s$, and take the $k$-th exterior power of both sides. The resulting asymptotic equivalence gives the claim.
\endproof

\begin{cor}\label{corstokmat}
If $\mathbb S_1^\Pb,\mathbb S_2^\Pb$ are the Stokes matrices for $\Pb^{n-1}$ attached with $\mc V$ and $C(\bm z)$, then the $k$-th exterior powers $\bigwedge^k\mathbb S_1^\Pb,\, \bigwedge^k\mathbb S_2^\Pb$ are the Stokes matrices for $G(k,n)$ attached with $e^{-\frac{\pi\sqrt{-1}(k-1)}{n}}\mc V$ and $\Tilde C(\bm z)$.\qed
\end{cor}

Recall the isomorphism $\mu_{k,n}\colon K^T_0(G(k,n))_\C\to\mathscr S_{H_{k,n}}$ of \eqref{defmukn}, which identifies solutions of the joint system of qDE and qKZ equations with $K$-theoretical classes. 

\begin{thm}\label{mainthm3}
There exist  ${T}$-full exceptional collections in ${\mathsf D}^b_T(G(k,n))$ whose $K$-classes are mapped via $\mu_{k,n}$ to Stokes fundamental solutions on the Stokes sectors of the qDE of $G(k,n)$. In particular, these Stokes solutions are solutions of the joint system of the qDE and qKZ equations.
\end{thm}

\proof
The statement has been proved for $\Pb^{n-1}$ in \cite[Section 10]{CV21}. Consequently, for $G(k,n)$, the Stokes fundamental solutions on Stokes sectors are the images via $\mu_{k,n}$ of suitable exceptional bases in $K_0^{T}(G(k,n))_\C$. This follows from Theorem \ref{satstoksol}, Theorem \ref{thmcompatibility} (in particular the commutativity of the diagram \eqref{finaldiag}), and also Corollary \ref{corkexcbases}. These exceptional bases are $K$-theoretical classes of ${T}$-full exceptional collections obtained by braid mutations of the twisted Kapranov collection $\left(\mathscr L\otimes\Si_{\bm\la}E_1^*\right)_{\bm\la\subseteq k\times(n-k)}$. This follows from the results established in \cite[Section 10]{CV21} for $\Pb^{n-1}$, and from Proposition \ref{propsatbr}.
\endproof

Consider the sectors
\beq\label{stsec1}\mathcal V_{k,\ell}':=\left\{s\in\widetilde{\mathbb C^*}\colon \frac{\ell}{n}-\frac{1}{2}-\frac{1}{2n}+\frac{k-1}{2n}<\phi<\frac{\ell}{n}+\frac{k-1}{2n}\right\},
\eeq
\beq\label{stsec2}
\mathcal V_{k,\ell}'':=\left\{s\in\widetilde{\mathbb C^*}\colon \frac{\ell}{n}-\frac{1}{2}-\frac{1}{n}+\frac{k-1}{2n}<\phi<\frac{\ell}{n}-\frac{1}{2n}+\frac{k-1}{2n}\right\},
\eeq
for $\ell\in\mathbb Z$.

\begin{lem}\label{stoksec}
The sectors $\mathcal V_{k,\ell}',\mathcal V_{k,\ell}''$, defined by \eqref{stsec1}-\eqref{stsec2}, are maximal Stokes sectors for $G(k,n)$ w.r.t.\,\,to the inclusion, i.e.
\begin{enumerate} 
\item they are Stokes sectors,
\item any Stokes sector $\mathcal V$ is contained in one (and only one) $\mathcal V_{k,\ell}'$ or $\mathcal V_{k,\ell}''$.
\end{enumerate}
\qed
\end{lem}

Let us make more explicit the nature of the exceptional bases corresponding to Stokes solutions on each Stokes sector $\mc V_{k,\ell}', \mc V''_{k,\ell}$.
\vskip2mm  
To that purpose, let us first introduce four families $Q'_\ell,Q''_\ell,\Tilde Q'_\ell, \Tilde Q''_\ell$, with $\ell\in\Z$, of ${T}$-full exceptional collections in ${\mathsf D}^b_T(\Pb^{n-1})$, as in \cite[Section 7]{CV21}.
In what follows $\mc {T}$ denotes the tangent sheaf of $\Pb^{n-1}$, and we also set 
\[\bigwedge\nolimits^h\mc {T}(m):=\left(\bigwedge\nolimits^h\mc {T}\right)\otimes\mc O(m),\quad h=1,\dots,n-1,\quad m\in\Z.
\]
\vskip2mm
First, let us define the collections $Q'_\ell,Q''_\ell$ for a given $\ell\in\Z$:
\vskip1,5mm
\begin{enumerate}
\item If $n$ is \emph{odd}, the collection $Q_\ell'$ is
{\footnotesize
\begin{align*}
\mathcal O\left(-\ell-\frac{n-1}{2}\right),\ &\mathcal {T}\left(-\ell-\frac{n-1}{2}-1\right),\ \mathcal O\left(-\ell-\frac{n-1}{2}+1\right),\ \bigwedge\nolimits^3\mathcal {T}\left(-\ell-\frac{n-1}{2}-2\right),\dots,\\
\dots&\ ,\ \bigwedge\nolimits^{n-4}\mathcal {T}\left(-\ell-n+2\right),\ \mathcal O(-\ell-1),\ \bigwedge\nolimits^{n-2}\mathcal {T}\left(-\ell-n+1\right),\ \mathcal O(-\ell),
\end{align*}
}
and the collection $Q_\ell''$ is
{\footnotesize
\begin{align*}
\mathcal O\left(-\ell-\frac{n-1}{2}\right),\ &\mathcal O\left(-\ell-\frac{n-1}{2}+1\right),\ \bigwedge\nolimits^2\mathcal {T}\left(-\ell-\frac{n-1}{2}-1\right),\ \mathcal O\left(-\ell-\frac{n-1}{2}+2\right),\dots,\\
\dots\ &,\ \mathcal O(-\ell-1),\ \bigwedge\nolimits^{n-3}\mathcal {T}\left(-\ell-n+2\right),\ \mathcal O(-\ell),\ \bigwedge\nolimits^{n-1}\mathcal {T}\left(-\ell-n+1\right).
\end{align*}
}
\item If $n$ is \emph{even}, the collection $Q_\ell'$ is
{\footnotesize
\begin{align*}
\mathcal O\left(-\ell-\frac{n}{2}\right),\ &\mathcal O\left(-\ell-\frac{n}{2}+1\right),\ \bigwedge\nolimits^2\mathcal {T}\left(-\ell-\frac{n}{2}-1\right),\ \mathcal O\left(-\ell-\frac{n}{2}+2\right),\dots,\\
\dots\ &,\ \bigwedge\nolimits^{n-4}\mathcal {T}\left(-\ell-n+2\right),\ \mathcal O(-\ell-1),\ \bigwedge\nolimits^{n-2}\mathcal {T}\left(-\ell-n+1\right),\ \mathcal O(-\ell),
\end{align*}
}and the collection $Q_\ell''$ is
{\footnotesize
\begin{align*}
\mathcal O\left(-\ell-\frac{n}{2}+1\right),\ &\mathcal {T}\left(-\ell-\frac{n}{2}\right),\ \mathcal O\left(-\ell-\frac{n}{2}+2\right),\ \bigwedge\nolimits^3\mathcal {T}\left(-\ell-\frac{n}{2}-1\right),\dots,\\
\dots&\ ,\ \mathcal O(-\ell-1),\ \bigwedge\nolimits^{n-3}\mathcal {T}\left(-\ell-n+2\right),\ \mathcal O(-\ell),\ \bigwedge\nolimits^{n-1}\mathcal {T}\left(-\ell-n+1\right).
\end{align*}
}
\end{enumerate}
In these exceptional collections, each of the objects 
$\mathcal O(m)$  sits in degree $0$ and each of the objects
  $\bigwedge^h\mathcal {T}(m)$ sits in degree $-h$.
\vskip2mm
Now, we can introduce the collections $\Tilde Q'_\ell, \Tilde Q''_\ell$. The objects of the 
collections $\widetilde{Q}_\ell',\widetilde{Q}_\ell''$ are obtained 
from the objects of
$Q_\ell',Q_\ell''$ by twisting their ${T}$-equivariant structures.
More precisely, for $a\in\mathbb Z$ define the ${T}$-characters 
\beq 
\underbrace{\bigwedge\nolimits ^n V\otimes\dots\otimes \bigwedge\nolimits ^n V}_{a\text{ times}},\quad \text{if}\
\ a\geq 0,
\eeq
\beq
 \underbrace{\bigwedge\nolimits ^n V^*\otimes\dots\otimes \bigwedge\nolimits ^n V^*}_{-a\text{ times}},\quad 
 \quad \text{if}\
\ a< 0,
\eeq
where $V\cong\mathbb C^n$ is the diagonal representation of ${T}$.
Given $m\in \mathbb Z$, define $a\in \mathbb Z$ from 
\beq
0\leq m+an\leq n-1.
\eeq
To construct the collections $\Tilde Q'_\ell, \Tilde Q''_\ell$: the ${T}$-equivariant structure of any
object $\mathcal O(-m)$ or $\bigwedge^{m-h}\mathcal {T}(-m)$ corresponding to bases
$Q_\ell',Q_\ell''$  must be tensored with the corresponding character defined above. See also \cite[Sec.\,7.6]{CV21}.
\vskip2mm
Let us describe the $K$-classes of the collections $\Tilde Q_\ell',\Tilde Q_\ell''$ as elements of $$K^T_0(\Pb^{n-1})\cong \Z[X^{\pm 1},\bm Z^{\pm 1}]/\Big\langle\prod_{j=1}^n(X-Z_j)\Big\rangle,\qquad X=\Gm_{1,1}^\Pb.$$
For any $h,m\in\Z$ such that $0\leq m-h\leq n$, we denote by $X^m(h)$ the polynomial
\[X^m(h):=X^m-e_1(\bm Z)X^{m-1}+\dots+(-1)^{m-h}e_{m-h}(\bm Z)X^h,
\]where $e_i(\bm Z)$ is the $i$-th elementary symmetric polynomial in $\bm Z$.
\vskip2mm
Denote by $[Q_\ell'],[Q''_\ell]$ the exceptional bases of $K^T_0(\Pb^{n-1})$ corresponding to $Q_\ell',Q''_\ell$.
The polynomials in $K^T_0(\Pb^{n-1})$ corresponding to the objects of $[Q_\ell'],[Q''_\ell]$ can be reconstructed by application of the following rule:
\begin{center}{\it
each object $\bigwedge^{m-h}\mc {T}(-m)$ (in degree $h-m$) appearing in the exceptional collections $Q_\ell',Q''_\ell$ is represented by the polynomial $X^m(h).$}
\end{center}
Finally, the basis $\Tilde Q'_\ell$ (resp.\,\,$\Tilde Q''_\ell$) is obtained from the basis $Q'_\ell$ (resp.\,\,$Q''_\ell$) by substituting any polynomial $X^m(h)$ with 
\[((-1)^{n+1}e_n(\bm Z))^aX^m(h),
\]where $a\in\Z$ is such that $0\leq m+an\leq n-1$.

\begin{example}
Let $n=5$, and $\ell=-1$. We have
\begin{align*}
Q'_{-1}&=(X^1,X^2(1),X^0,X^3,X^{-1}),\\
\Tilde Q'_{-1}&=(X^1,X^2(1),X^0,X^3,e_5(\bm Z)X^{-1}),\\
Q''_{-1}&=(X^1,X^0,X^2,X^{-1},X^3(-1)),\\
\Tilde Q''_{-1}&=(X^1,X^0,X^2,e_5(\bm Z)X^{-1},X^3(-1)).
\end{align*}
\qetr
\end{example}

The exceptional bases \( [\Tilde{Q}'_\ell ]\) and \( [\Tilde{Q}''_\ell] \) of \( K^T_0(\mathbb{P}^{n-1}) \) induce, by Corollary \ref{corkexcbases}, two exceptional bases \( [^{(k)}\!\mathcal{Q}'_\ell ]\) and \( [^{(k)}\!\mathcal{Q}''_\ell] \) of \( K^T_0(G(k,n)) \). These are the $K$-classes of two ${T}$-full exceptional collections  \( ^{(k)}\!\mathcal{Q}'_\ell \) and \( ^{(k)}\!\mathcal{Q}''_\ell \) of \( {\mathsf D}^b_T(G(k,n))\), suitable mutations (up to some shifts of objects) of the twisted Kapranov collection $\left(\mathscr L\otimes\Si_{\bm\la}E_1^*\right)_{\bm\la\subseteq k\times(n-k)}$. 
\vskip2mm
The polynomials representing the basis \( [^{(k)}\!\mathcal{Q}'_\ell ]\) (resp.\,\,\( [^{(k)}\!\mathcal{Q}''_\ell ]\)) are obtained via formula \eqref{forthi+} by specializing \( f_1(X; \bm{Z}), \dots, f_k(X; \bm{Z}) \) to the polynomials of \( [\Tilde{Q}'_\ell] \) (resp. \( [\Tilde{Q}''_\ell ]\)).

\begin{example}
Let $n=5, k=4$ and $\ell=-1$. The basis $[^{(4)}\!\mc Q'_{-1}]$ can be computed from the basis $[\Tilde Q'_{-1}]$ of the previous example. Set
\[f_1(X;\bm Z)=X,\quad f_2(X;\bm Z)=X^2-e_1(\bm Z)X,\quad f_3(X;\bm Z)=1,\]\[ f_4(X;\bm Z)=X^3,\quad f_5(X;\bm Z)=e_5(\bm Z)X^{-1}.
\]The polynomials of the basis $[^{(4)}\!\mc Q'_{-1}]$ are
\[\thi_{4,5}^K(f_1\wedge f_2\wedge f_3\wedge f_4)=1,
\]
\[\thi_{4,5}^K(f_1\wedge f_2\wedge f_3\wedge f_5)=-\frac{e_5(\bm Z)}{\Gm_{1,1} \Gm_{1,2} \Gm_{1,3} \Gm_{1,4}},
\]
\[\thi_{4,5}^K(f_1\wedge f_2\wedge f_4\wedge f_5)=e_5(\bm Z) \left(-\frac{1}{\Gm_{1,1}}-\frac{1}{\Gm_{1,2}}-\frac{\Gm_{1,3}+\Gm_{1,4}}{\Gm_{1,3} \Gm_{1,4}}\right),
\]
\[\thi_{4,5}^K(f_1\wedge f_3\wedge f_4\wedge f_5)=\frac{e_5(\bm Z) (\Gm_{1,1}+\Gm_{1,2}+\Gm_{1,3}+\Gm_{1,4})}{\Gm_{1,1} \Gm_{1,2} \Gm_{1,3} \Gm_{1,4}},
\]
\begin{multline*}\thi_{4,5}^K(f_2\wedge f_3\wedge f_4\wedge f_5)=\\ \frac{e_5(\bm Z) (-e_1(\bm Z)(\Gm_{1,1}+\Gm_{1,2}+\Gm_{1,3}+\Gm_{1,4})+\Gm_{1,1} (\Gm_{1,2}+\Gm_{1,3}+\Gm_{1,4})+\Gm_{1,4} (\Gm_{1,2}+\Gm_{1,3})+\Gm_{1,2} \Gm_{1,3})}{\Gm_{1,1} \Gm_{1,2} \Gm_{1,3} \Gm_{1,4}}.\end{multline*}
\qetr
\end{example}

\begin{thm}\label{thmqstok}
Via the map $\mu_{k,n}$, the exceptional basis \( [^{(k)}\!\mathcal{Q}'_\ell ]\) (resp.\,\,\( [^{(k)}\!\mathcal{Q}''_\ell] \)) is mapped to the Stokes basis of the joint system of qDE and qKZ equations for $G(k,n)$ on the Stokes sector $\mc V'_{k,\ell}$ (resp.\,\,$\mc V''_{k,\ell}$), with normalization
\[C(\bm z)=\mathsf c_{k,n}\cdot\diag\left(\mathsf a_{I}\cdot(e^{\frac{\pi\sqrt{-1}k(k-1)}{n}}\frak p_I^{-1})^{\sum_{a=1}^nz_a+\frac{n-1}{2}}\right)_{I\in\mc I_{k,n}},
\]where $\mathsf c_{k,n}:=(2\pi)^{\frac{k(n-1)}{2}}e^{-\pi\sqrt{-1}\frac{k(n-1)}{2}}$, and $\mathsf a_I:=\prod_{i\in I_1}e^{\frac{(i-1)\pi\sqrt{-1}}{n}}$.
\end{thm}

\proof 
It follows from the results for $\Pb^{n-1}$ in \cite[Section 10]{CV21}, and Theorem \ref{satstoksol}.
\endproof

\subsection{Dual exceptional collections, Serre functor, helices} Let $\frak E$ be a ${T}$-full exceptional collection in ${\mathsf D}^b(G(k,n))$. Introduce the braid 
\beq
\beta:=\tau_1(\tau_{2}\tau_{1})\dots(\tau_{\binom{n}{k}-2}\dots \tau_{1})(\tau_{\binom{n}{k}-1}\tau_{\binom{n}{k}-2}\dots \tau_{1})\in\mc B_{\binom{n}{k}},
\eeq
and define the {\it left} and {\it right dual exceptional collections}, $^\vee\frak E$ and $\frak E^\vee$, as the mutations
\beq
^\vee\frak E:=\bt\,\frak E,\qquad \frak E^\vee:=\bt^{-1}\,\frak E.
\eeq
Similarly, given an exceptional basis $\eps:=(e_1,\dots, e_{\binom{n}{k}})$ of $K^T_0(G(k,n))$, we can introduce the dual exceptional bases $^\vee\eps:=(^\vee e_1,\dots,^\vee e_{\binom{n}{k}})$ and $\eps^\vee:=(e_1^\vee,\dots,e_{\binom{n}{k}}^\vee)$. 
\vskip2mm
Given a matrix $A\in M_n(R({T})_\C)=M_h(\C[\bm Z^{\pm 1}])$, we define two matrices $A^*,A^\dag\in M_h(R({T})_\C)$ as follows: the matrix $A^*$ is obtained by applying the duality involution \eqref{dualrap} entrywise, and the matrix $A^\dag$ is defined by 
\beq
(A^\dag)_{a,b}:=A^*_{b,a},\quad a,b=1,\dots, h.
\eeq
\begin{prop}[\cite{CV21,CDG1}]$\quad$
\begin{enumerate}
\item 
The following orthogonality conditions hold:
\beq
\chi^{T}(^\vee e_i,e_j)=\dl_{i+j,\binom{n}{k}+1},\quad \chi^{T}(e_i,e_j^\vee)=\dl_{i+j,\binom{n}{k}+1},\quad i,j=1,\dots,\binom{n}{k}.
\eeq
\item For any $v\in K^T_0(G(k,n))$, we have
\beq
\label{dual1}v=\sum_{h=1}^{\binom{n}{k}}\chi^G(v,e_h^\vee)^*\ e_{\binom{n}{k}+1-h},\quad
v=\sum_{h=1}^{\binom{n}{k}}\chi^G(^\vee e_h, v)\ e_{{\binom{n}{k}}+1-h}.
\eeq
\item Let $\mc G$ be the Gram matrix of $\chi^{T}$ w.r.t. the exceptional basis $\eps$ of $K^T_0(G(k,n))$. Then the Gram matrix of $\chi^{T}$ w.r.t.\,\,both $^\vee\eps$ and 
$\eps^\vee$ 
equals
\beq
J\cdot(\mc G^\dag)^{-1}\cdot J,\quad \text{where }J_{a,b}=\dl_{a+b,\binom{n}{k}+1},\quad a,b=1,\dots,\binom{n}{k}.
\eeq
\end{enumerate}\qed
\end{prop}

Let $\om_{G(k,n)}^{T}$ be the ${T}$-equivariant canonical sheaf of $G(k,n)$. The {\it Serre functor} is the functor $\mc K\colon {\mathsf D}^b_T(G(k,n))\to {\mathsf D}^b_T(G(k,n))$ defined by
\beq
\mc K(E):=(\om^T_{G(k,n)}\otimes E)[\dim_CG(k,n)],\quad E\in\on{Ob}({\mathsf D}^b_T(G(k,n))).
\eeq The Serre functor is uniquely characterized (up to canonical isomorphism) by the condition
\beq\label{serre1}
\Hom_T^\bullet(E,F)^*\cong\Hom_T(F,\mc K(E)),\quad E,F\in\on{Ob}({\mathsf D}^b_T(G(k,n))).
\eeq
From this it also follows the {\it Serre periodicity}
\beq\label{serre2}
\Hom^\bullet_T(E,F)\cong \Hom^\bullet_T(\mc K(E),\mc K(F)),\quad E,F\in\on{Ob}({\mathsf D}^b_T(G(k,n))).
\eeq

\begin{prop}\label{serreexc}\cite{CV21}
Let $\frak E$ 
be a ${T}$-full exceptional collection of ${\mathsf D}^b_G(X)$. The following operations are equivalent, i.e. produce the same exceptional collection when applied to $\frak E$:
\begin{enumerate}
\item to act on $\frak E$ with the braid $(\tau_1\dots\tau_{\binom{n}{k}-1})^{-\binom{n}{k}}$,
\item to take the double right-dual exceptional collection $(\frak E^\vee)^\vee$,
\item to apply the Serre functor to each object of $\frak E$. 
\end{enumerate}
\qed
\end{prop}

The $K$-theoretical analog of the Serre functor is the {\it canonical operator} 
\[{\mathsf k}\colon K^T_0(G(k,n))\to  K^T_0(G(k,n)),\quad [E]\mapsto [\mc K(E)].
\]The $K$-theoretical analogs of the isomorphisms \eqref{serre1}, \eqref{serre2} are
\beq\label{cano}
\chi^{T}(e,f)^*=\chi^{T}(f,{\mathsf k}(e)),\quad \chi^{T}(e,f)=\chi^{T}({\mathsf k}(e),{\mathsf k}(f)),\quad e,f\in  K^T_0(G(k,n)).
\eeq In particular, the canonical operator ${\mathsf k}$ is a $\chi^{T}$-isometry.

\begin{prop}\label{matrk}$\quad$
\begin{enumerate}
\item The canonical operator ${\mathsf k}\colon K^T_0(G(k,n))\to  K^T_0(G(k,n))$ acts as follows on polynomials
\beq\label{polykappa}
f(\bm\Gm;\bm Z)\mapsto (-1)^{k(n-k)}\frac{(\Gm_{1,1}\dots\Gm_{1,k})^n}{(Z_1\dots Z_n)^k}f(\bm \Gm;\bm Z).
\eeq
\item Let $\varepsilon$ be a basis of $K_0^{T}(G(k,n))$, and  $\mathcal G$  the Gram matrix of $\chi^{T}$ w.r.t.\,\,$\varepsilon$. Then the matrix of the canonical operator ${\mathsf k}\colon K_0^{T}(G(k,n))\to K_0^{T}(G(k,n))$ w.r.t.\,\,the basis $\varepsilon$ equals
\beq\label{eqgramk}
\mathcal G^{-1}\mathcal G^\dag.
\eeq
\end{enumerate}
\end{prop}
\proof 
Point (1) follows from the identity $[\om^T_{G(k,n)}]=(\prod_{i=1}^k\Gm_{1,i})^n(\prod_{j=1}^nZ_j)^{-k}$ in $K$-theory. This holds because
\[\om^T_{G(k,n)}=(\det TG(k,n))^*=\det(E_1\otimes E_2^*)=(\det E_1)^{\otimes(n-k)}\otimes(\det E_2^*)^{\otimes k},
\]and from the short exact sequence \eqref{Euler}, we deduce 
\begin{align*}0\to \det E_1\to\det \underline{\C^n}\to \det E_2\to 0\quad&\Longrightarrow\quad \det E_1\otimes \det E_2=\det \underline{\C^n}\\
&\Longrightarrow\quad [\det E_2^*]=\frac{[\det E_1]}{Z_1\dots Z_n}\\
&\Longrightarrow\quad [\om^T_{G(k,n)}]=\frac{[(\det E_1)^{\otimes n}]}{(Z_1\dots Z_n)^k}.
\end{align*}
Point (2) follows from the first of identities \eqref{cano}, written in matrix notation.
\endproof

Denote by $\bigwedge^k{\mathsf k}^\Pb\in \on{End}_\C\left( \bigwedge^kK^T_0(\Pb^{n-1})_\C\right)$ the linear map induced by the canonical operator ${\mathsf k}^\Pb$ of $\Pb^{n-1}$, and by ${\mathsf k}^G$ the canonical operator of $G(k,n)$.
\begin{prop}\label{eigenk}
The Satake identification $\thi^K_{k,n}\colon \bigwedge^kK^T_0(\Pb^{n-1})\to K^T_0(G(k,n))$ intertwines the operator $\bigwedge^k{\mathsf k}^\Pb$ and the canonical operator ${\mathsf k}^G$. Consequently, the eigenvalues of ${\mathsf k}^G$ are
\beq\label{eigenvaluescanop}
\nu_I(\bm Z)=
(-1)^{k(n-k)}\frac{(\prod_{i\in I_1}Z_{i})^n}{(Z_1\dots Z_n)^k},\quad I=(I_1,I_2)\in\mc I_{k,n}.
\eeq
\end{prop}

\proof
The result directly follows from the polynomial representation of ${\mathsf k}^\Pb$, and formula \eqref{forthi+}.
\endproof

Collect the polynomials $\nu_I(\bm Z)$ into the tuple $\bm\nu(\bm Z):=\left(\nu_I(\bm Z)\right)_{I\in\mc I_{k,n}}$, and set
\beq\label{nujz}
\tilde\nu_{k,n,j}(\bm Z):=e_j\left(\bm \nu(\bm Z)\right),\quad e_j\text{ is $j$-th elementary symmetric polynomial},\quad j=0,\dots,\binom{n}{k}.
\eeq

\begin{cor}\label{diophantusg}
Let $\varepsilon$ be a basis of $K_0^{T}(G(k,n))$, and  $\mathcal G$  the Gram matrix of $\chi^{T}$ w.r.t.\,\,$\varepsilon$. We have
\beq
\det({\mathsf t}\cdot{\bf 1}-\mc G^{-1}\mc G^\dag)=\sum_{j=0}^{\binom{n}{k}}(-1)^j\,\tilde\nu_{k,n,j}(\bm Z)\,{\mathsf t}^{\binom{n}{k}-j}.
\eeq
\end{cor}
\proof
A direct consequence of Proposition \ref{matrk} and Proposition \ref{eigenk}.
\endproof

Given a ${T}$-full exceptional collection $\frak E=(E_1,\dots,E_{\binom{n}{k}})$ of ${\mathsf D}^b_T(G(k,n))$, the {\it helix} generated by $\frak E$ is the infinte family of objects $(E_i)_{i\in\Z}$ obtained by the iterated mutations
\[E_{i+\binom{n}{k}}:=\mathbb R_{E_{i+\binom{n}{k}-1}}\dots\mathbb R_{E_{i+1}}E_i,\quad E_{i-\binom{n}{k}}:=\mathbb L_{E_{i-\binom{n}{k}+1}}\dots\mathbb L_{E_{i-1}}E_i,\quad i\in\mathbb Z.
\]Such a helix is said to be of \emph{period $\binom{n}{k}$}. Any family of $\binom{n}{k}$ consecutive objects $(E_i,\dots, E_{i+\binom{n}{k}})$ is called a \emph{foundation} of the helix.

\begin{prop}\label{prophelix}
Let $(E_i)_{i\in\Z}$ be the helix generated by $\frak E=(E_1,\dots,E_{\binom{n}{k}})$. For any $a\in\Z$, the foundation $\frak E_a=(E_{a\binom{n}{k}+1},\dots,E_{a\binom{n}{k}+\binom{n}{k}})$ is obtained by applying $-a$ times the Serre functor to $\frak E$, that is
$\frak E_a=\mc K^{-a}\frak E.$\qed
\end{prop}

\subsection{Stokes matrices as Gram matrices}In this section we relate the Stokes phenomenon of Stokes solutions of the joint system of qDE and qKZ equations with the Euler--Poincar\'e--Grothendieck pairing in $K^T_0(G(k,n))$, defined in \eqref{epgm}. More precisely, we prove that the Stokes matrices equal the Gram matrices of $\chi^{T}$ w.r.t.\,\,the exceptional bases \( [^{(k)}\!\mathcal{Q}'_\ell] \) and \( [^{(k)}\!\mathcal{Q}''_\ell] \).
\vskip2mm
Let $\mc V$ be a Stokes sector for the qDE of $G(k,n)$, and let $Y(\bm z;s)$ be the Stokes fundamental solution matrix on $\mc V$ associated with the ${T}$-full exceptional collection $\frak E$. The sector $\mc V$ is contained in one (and only one) of the sectors $\mc V'_{k,\ell}$ or $\mc V''_{k,\ell}$, $\ell\in\Z$, by Lemma \ref{stoksec}. Consequently, by Theorem \ref{thmqstok}, the collection $\frak E$ projects in $K$-theory either on the exceptional basis  \( [^{(k)}\!\mathcal{Q}'_\ell ]\) or \( [^{(k)}\!\mathcal{Q}''_\ell] \).
\vskip2mm
Let $Y_1(\bm z;s)$ and $Y_2(\bm z; s)$ be the Stokes fundamental solutions matrices on the Stokes sectors $e^{\pi\sqrt{-1}}\mc V$ and $e^{2\pi\sqrt{-1}}\mc V$, respectively. Let $\frak E_1$ and $\frak E_2$ be the corresponding ${T}$-full exceptional collections. Recall that $Y(\bm z;s), Y_1(\bm z;s)$ and $Y_2(\bm z;s)$ are related by the Stokes matrices $\mathbb S_1,\mathbb S_2$ as in equation \eqref{stokmatr}. 

\begin{thm}\label{thmleftdual}
We have $[\frak E_1]=\,[^\vee\frak E]$ and $[\frak E_2]=\,[^\vee(^\vee\frak E)]$.
\end{thm}
\proof
The statement holds true for $\Pb^{n-1}$, see \cite[Section 11]{CV21}. The result for $G(k,n)$ then follows by the Satake correspondence: the left/right duality relations between the exceptional bases $\eps,\eps^\vee,\,^\vee\eps$ are preserved on the induced exceptional bases in the $K$-theory of $G(k,n)$. This follows by the definition of dual exceptional bases and Proposition \ref{propisomkth}.
\endproof

\begin{thm}\label{corstokmatr}
Let $J$ be the anti-diagonal matrix 
\[J_{ab}=\dl_{a+b,\binom{n}{k}+1},\quad a,b=1,\dots,\binom{n}{k},
\]and let $\mc G$ be the Gram matrix of $\chi^{T}$ w.r.t.\,\,the exceptional basis $[\frak E]$.
\begin{enumerate}
\item The Stokes matrix $\mathbb S_1$ equals the Gram matrix of $\chi^{T}$ w.r.t.\,\,the left dual exceptional basis $[\frak E_1]=[^\vee\frak E]$, that is
\beq
\mathbb S_1=J\left(\mc G^\dag\right)^{-1} J.
\eeq
\item The matrix $J\mathbb S_2J$ equals the Gram matrix of $\chi^{T}$ w.r.t.\,\,the exceptional basis $[\frak E]$, that is
\beq
\mathbb S_2=J\mc GJ.
\eeq
\item Both $\mathbb S_1$ and $\mathbb S_2$ have entries in the ring of symmetric Laurent polynomials with integer coefficients, that is $\mathbb S_1,\mathbb S_2\in M_n(\Z[\bm Z^{\pm 1}]^{ S_n})$. We also have
\beq
\mathbb S_2=\left(\mathbb S_1^\dag\right)^{-1}.
\eeq
\end{enumerate}
\end{thm}

\proof All the points hold for $\Pb^{n-1}$, and Corollary \ref{corstokmat} implies the results for $G(k,n)$. Alternatively, we can argue as follows.
Points (1) and (2) follow from Theorem \ref{thmleftdual} and the second identity \eqref{dual1}. For point (3), notice that the Gram matrix of $\chi^{T}$ w.r.t.\,\,the Kapranov exceptional basis $\left([\Si_{\bm\la}E_1^*]\right)_{\bm\la}=\left(s_{\bm\la}(\bm\Gm_1^{-1})\right)_{\bm\la}$ is with integer symmetric Laurent polynomial entries, by Corollary \ref{corchikapr}. The braid group action preserves this property.
\endproof

Recall the polynomials $\tilde\nu_{k,n,j} (\bm Z)$ defined in equation \eqref{nujz}. The following result gives some non-trivial diophantine constraints on the Stokes matrices. 
\begin{cor}
For both $\mathbb S=\mathbb S_1,\mathbb S_2$, we have
\beq\label{constrS}
\on{Tr} \bigwedge\nolimits^{\ell}\left(\mathbb S^\dag\mathbb S^{-1}\right)=\tilde\nu_{k,n,\ell}(\bm Z),\quad \ell=1,\dots,\binom{n}{k}.
\eeq
\end{cor}

\proof For an arbitrary matrix $A$, we have $\det(\mathsf t\cdot{\bf 1}-A)=\sum_{i\geq 0}(-1)^i\on{Tr}(\bigwedge^iA)\,\mathsf t^i$. 
The result then follows from Theorem \ref{corstokmatr} and Corollary \ref{diophantusg}. 
\endproof

\begin{example}
Consider the case $(k,n)=(1,3)$, i.e. the case of the projective plane $\Pb^2$. Set
\[\mathbb S=\begin{pmatrix}
1&a&b\\
0&1&c\\
0&0&1
\end{pmatrix},\quad a,b,c\in\Z[\bm Z^{\pm 1}]^{ S_3}.
\]We have two non-trivial diophantine constraints \eqref{constrS} (for $\ell=1,2$):
\begin{align}\label{markov1}
aa^*+bb^*+cc^*-ab^*c&=3-\frac{Z_1^3+Z_2^3+Z_3^3}{Z_1Z_2Z_3},\\
\label{markov2}
aa^*+bb^*+cc^*-a^*bc^*&=3-\frac{Z_2^3Z_3^3+Z_1^3Z_3^3+Z_2^3Z_3^3}{Z_1^2Z_2^2Z_3^2}.
\end{align}
Notice that \eqref{markov2} is equivalent to \eqref{markov1}, by taking the $*$-dual of both sides. The equation \eqref{markov1} coincides with the $*$-Markov equation for Laurent polynomials\footnote{In \cite{CVmarkov} we changed variables $(a,b,c)\mapsto (a,b^*,c)$. The resulting $*$-Markov equation takes the more symmetric form $aa^*+bb^*+cc^*-abc=3-\frac{Z_1^3+Z_2^3+Z_3^3}{Z_1Z_2Z_3}$}, extensively studied in \cite{CVmarkov}.
\vskip1,5mm
Consider now the case $(k,n)=(2,3)$, i.e. the case of the Grassmannian $G(2,3)$ of $2$-planes in $\C^3$. In this case as well, we have the two diophantine constraints \eqref{constrS} (for $\ell=1,2$):
\begin{align}\label{23markov1}
aa^*+bb^*+cc^*-ab^*c&=3-\frac{Z_2^3Z_3^3+Z_1^3Z_3^3+Z_2^3Z_3^3}{Z_1^2Z_2^2Z_3^2},\\
\label{23markov2}
aa^*+bb^*+cc^*-a^*bc^*&=3-\frac{Z_1^3+Z_2^3+Z_3^3}{Z_1Z_2Z_3}.
\end{align}
Once again, we recognize the $*$-Markov equation: if $(a,b,c)$ is a solution of \eqref{markov1}, then $(a^*,b^*,c^*)$ is a solution of \eqref{23markov1} and \eqref{23markov2}. This is coherent with the fact that $G(2,3)$ actually is isomorphic to $\Pb^2$, being $G(2,3)$ the projective plane for the dual space $(\C^3)^*$.
\vskip1,5mm
The Satake correspondence also gives us a further way to produce solutions of \eqref{23markov1} and \eqref{23markov2} starting from solutions of \eqref{markov1}: if $(a,b,c)$ is a solution of \eqref{markov1}, then $(c,ac-b,a)$ is a solution of \eqref{23markov1} and \eqref{23markov2}. This follows from Corollary \ref{corstokmat}.
\qetr
\end{example}

\begin{rem} Let us briefly discuss two general features observed in the previous example.
For fixed $n$, set $R_n:=\Z[\bm Z^{\pm 1}]^{ S_n}$, $N_{k,n}:=\frac{1}{2}\binom{n}{k}\left[\binom{n}{k}-1\right]$, and introduce the $R_n$-module $\mathsf A_{k,n}:=R_n^{N_{k,n}}$. We can identify $\mathsf A_{k,n}$ with the space of all possible Stokes matrices $\mathbb S$ above: each $N_{k,n}$-tuple of polynomials in $R_n$ gives the non-trivial (i.e.\,\,$\neq 0,1$) entries of $\mathbb S$. 

Denote by $\mathsf V_{k,n}\subseteq \mathsf A_{k,n}$ the zero locus defined by the equations \eqref{constrS}. For example, ${\mathsf V_{1,3}}$ is the set of triples $(a,b,c)$ solving equations \eqref{markov1}, \eqref{markov2}, and ${\mathsf V_{2,3}}$ is the set of triples $(a,b,c)$ solving equations \eqref{23markov1}, \eqref{23markov2}. 

We have $N_{k,n}=N_{n-k,n}$, so that the duality map $$g\colon A_{k,n}\to A_{n-k,n},\quad (f_1(\bm Z),\dots, f_{N_{k,n}}(\bm Z))\mapsto (f_1(\bm Z^{-1}),\dots, f_{N_{k,n}}(\bm Z^{-1}))$$ is well-defined. It is easy to see that $g(\mathsf V_{k,n})=\mathsf V_{n-k,n}$. This follows from the identities
\[\tilde\nu_{k,n,j}(\bm Z)=\tilde\nu_{k,n,\binom{n}{k}-j}(\bm Z)^*=\tilde\nu_{n-k,n,\binom{n}{k}-j}(\bm Z),\qquad j=0,\dots,\binom{n}{k}.
\]The identification of $\mathsf V_{k,n}$ with $\mathsf V_{n-k,n}$ via $g$ is consistent with the fact that $G(n-k,n)$ can be identified with the Grassmannian of $k$-planes in the dual space $(\C^*)^n$.

Furthermore, the Satake correspondence (more precisely, Corollary \ref{corstokmat}) gives a polynomial map $f\colon \mathsf A_{1,n}\to\mathsf A_{k,n}$ of degree $k$ such that $f(\mathsf V_{1,n})\subset \mathsf V_{k,n}$.
\end{rem}

To conclude, let us prove a structural constraint between the monodromy operator $M_0(\bm z)$, the formal monodromy $\La_{k,n}(\bm z)$, and the Stokes matrices $\mathbb S_1,\mathbb S_2$. 

\begin{rem}Notice that the monodromy operator $M_0(\bm z)$, as originally defined in Section \ref{seccmonop}, acts as $Y(\bm z;s)\mapsto Y(\bm z;\zeta_ns)$ on the space of solutions of the equation \eqref{js1}. The operator $Y(\bm z;s)\mapsto Y(\bm z;e^{2\pi\sqrt{-1}}s)$ equals $M_0(\bm z)^n$.
\end{rem}

\begin{prop}
For generic $\bm z$, we have 
\beq\label{constraint}
M_0(\bm z)^n=\exp\left(2\pi\sqrt{-1}\La_{k,n}(\bm z)\right)(\mathbb S_1\mathbb S_2)^{-1}.
\eeq
\end{prop}
\proof
Recall $\mc U$ as in \eqref{eqgenericz}.
In the  notations of Definition \ref{stokesmatr}, the following identities hold true for any $s\in\widetilde{\mathbb C^*}$ and $\bm z\in\mc U$:
\begin{enumerate}
\item $Y_{2}\left(\bm z; e^{2\pi\sqrt{-1}}s\right)=Y(\bm z;s)\cdot \exp\left(2\pi\sqrt{-1}\Lambda_{k,n} (\bm z)\right)$,
\item $Y_2(\bm z; s)=Y(\bm z;s)\cdot \mathbb S_1\mathbb S_2$,
\item $Y\left(\bm z;e^{2\pi\sqrt{-1}}s\right)=Y(\bm z;s)\cdot \exp\left(2\pi\sqrt{-1}\Lambda_{k,n} (\bm z)\right)\cdot \left(\mathbb S_1\mathbb S_2\right)^{-1}$.
\end{enumerate}
For (1), notice that
\[Y_{2}\left(\bm z;e^{2\pi\sqrt{-1}}s\right)\cdot \exp\left(-2\pi\sqrt{-1}\Lambda_{k,n} (\bm z)\right)
\]is a solution of \eqref{js1} with asymptotic expansion $Y_{\rm for}(\bm z;s)$ on the Stokes sector $\mathcal V$. Hence it must coincide with $Y(\bm z;s)$. Point (2) is a direct consequence of the definition of Stokes matrices. Point (3) follows from points (1) and (2).
\endproof

The constraint \eqref{constraint} has a clear $K$-theoretical interpretation: written in the form
\[(\mathbb S_1\mathbb S_2)^{-1}=\exp\left(-2\pi\sqrt{-1}\La_{k,n}(\bm z)\right)M_0(\bm z)^n,
\]the constraint is equivalent to Proposition \ref{matrk}. Indeed, 
under the identification of solutions and $K$-theoretical classes, both sides of the equation above can be identified with the canonical operator $\mathsf k$:
\begin{itemize} 
\item 
The exceptional collections $\frak E,\frak E_2$, associated with the Stokes solutions $Y(\bm z;s)$ and $Y_2(\bm z;s)$, respectively, are related by the identity $\mc K\frak E_2=\frak E$, by Proposition \ref{prophelix}. Consequently, the matrix $(\mathbb S_1\mathbb S_2)^{-1}$ represents the canonical morphism $\mathsf k$. Alternatively, the product $(\mathbb S_1\mathbb S_2)^{-1}$ is of the form \eqref{eqgramk}, by Theorem \ref{corstokmatr}.
\item The r.h.s.\,\,exactly is the polynomial representation \eqref{polykappa} of the canonical operator $\mathsf k$. This directly follows from Proposition \ref{propmono}.
\end{itemize}

\end{document}